\documentclass{article}

\PassOptionsToPackage{numbers, compress}{natbib}

\usepackage[final]{neurips_2023}

\usepackage[utf8]{inputenc} 
\usepackage[T1]{fontenc}    
\usepackage{hyperref}       
\usepackage{url}            
\usepackage{booktabs}       
\usepackage{amsfonts}       
\usepackage{nicefrac}       
\usepackage{microtype}      
\usepackage{xcolor}         

\usepackage{graphicx}
\usepackage{epstopdf} 
\usepackage{subfigure}
\usepackage{caption}
\usepackage{algorithmic,algorithm}

\newcommand{\X}{\mathcal{X}}
\newcommand{\blambda}{\boldsymbol{\lambda}}
\newcommand{\bzt}{\boldsymbol{\zeta}}
\newcommand{\bxi}{\boldsymbol{\xi}}

\newcommand{\reals}{\mathbb{R}}

\newcommand{\argmin}{\mathop{\mathrm{arg\,min}{}}}

\newcommand{\ba}{{\mathbf a}}

\newcommand{\bx}{{\mathbf x}}
\newcommand{\bz}{{\mathbf z}}
\newcommand{\by}{{\mathbf y}}
\newcommand{\bu}{{\mathbf u}}

\newcommand{\bb}{{\mathbf b}}

\newcommand{\bA}{{\mathbf A}}

\newcommand{\vertiii}[1]{{\left\vert\kern-0.20ex\left\vert\kern-0.20ex\left\vert #1
		\right\vert\kern-0.20ex\right\vert\kern-0.20ex\right\vert}}

\usepackage{amsmath}
\usepackage{amssymb}
\usepackage{mathtools}
\usepackage{amsthm}

\theoremstyle{plain}
\newtheorem{theorem}{Theorem}[section]
\newtheorem{proposition}[theorem]{Proposition}
\newtheorem{lemma}[theorem]{Lemma}

\theoremstyle{definition}
\newtheorem{definition}[theorem]{Definition}
\newtheorem{assumption}[theorem]{Assumption}
\theoremstyle{remark}
\newtheorem{remark}[theorem]{Remark}

\title{Oracle Complexity of Single-Loop Switching Subgradient Methods for Non-Smooth Weakly Convex Functional Constrained Optimization}

\author{%
  Yankun Huang \\
  Department of Business Analytics \\
  University of Iowa \\
  Iowa City, IA 52242 \\
  \texttt{yankun-huang@uiowa.edu} \\
  \And
  Qihang Lin \\
  Department of Business Analytics \\
  University of Iowa \\
  Iowa City, IA 52242 \\
  \texttt{qihang-lin@uiowa.edu} \\
}

\begin{document}

\maketitle

\begin{abstract}
We consider a non-convex constrained optimization problem, where the objective function is weakly convex and the constraint function is either convex or weakly convex. To solve this problem, we consider the classical switching subgradient method, which is an intuitive and easily implementable first-order method whose oracle complexity was only known for convex problems. This paper provides the first analysis on the oracle complexity of the switching subgradient method for finding a nearly stationary point of non-convex problems. Our results are derived separately for convex and weakly convex constraints. Compared to existing approaches, especially the double-loop methods, the switching gradient method can be applied to non-smooth problems and achieves the same complexity using only a single loop, which saves the effort on tuning the number of inner iterations. 
\end{abstract}

\section{Introduction}
Continuous optimization with nonlinear constraints arises from many applications of machine learning and statistics. Examples include 
Neyman-Pearson classification~\citep{rigollet2011neyman} and learning with fairness constraints~\citep{zafar2019fairness}. In this paper, we consider the following general nonlinear constrained optimization  
	\begin{eqnarray}
	\label{eq:gco}
	f^*\equiv \min_{\bx\in\mathcal{X}}f(\bx) \quad \text{s.t.} \quad g(\bx)\leq 0,
	\end{eqnarray}
where $\mathcal{X}\subset\mathbb{R}^d$ is a closed convex set that allows a computationally easy projection operator, $f$ is weakly-convex, $g$ is either convex or weakly convex, and functions $f$ and $g$ are not necessarily smooth.
When $g(\bx)\equiv\max_{i=1,\dots,m}g_i(\bx)$, \eqref{eq:gco} is equivalent to an optimization problem with multiple nonlinear constraints $g_i(\bx)\leq0$ for $i=1,\dots,m$.


A weakly convex function can be non-convex, so computing an optimal solution of \eqref{eq:gco} is challenging in general even without constraints. For this reason, theoretical analysis for gradient-based algorithms for non-convex problems mostly focuses on an algorithm's (oracle) complexity for finding an $\epsilon$-stationary solution for \eqref{eq:gco}. When a problem is non-smooth, finding an $\epsilon$-stationary solution is generally difficult even if the problem is convex~\cite{kornowski2022complexity}. Hence, in this paper, we consider finding a nearly $\epsilon$-stationary solution for \eqref{eq:gco}, whose definition will be stated later in Definition~\ref{def:stationary}. 

In the past decade, there have been many studies on non-convex constrained optimization. However, most of the existing algorithms and their theoretical complexity analysis are developed by assuming $f$ and $g_i$'s are all smooth or can be written as the sum of a smooth and a simple non-smooth functions. A non-exhaustive list of the works under such an assumption includes \cite{facchinei2021ghost, zhang2022global, zhang2022iteration, melo2020iteration, kong2023iteration, kong2023accelerated, kong2022iteration, li2021augmented, sahin2019inexact, lin2022complexity, li2021rate, boob2022level,  curtis2021inexact, berahas2021sequential,curtis2021worst,jin2022stochastic,li2022zeroth}. Their results cannot be applied to \eqref{eq:gco} due to non-smoothness in the problem.

Non-smoothness is common in optimization in machine learning, e.g., when a non-smooth loss function is applied, but there are much fewer studies on non-smooth non-convex constrained optimization. Under the weak-convexity assumption, an effective approach for solving a non-smooth non-convex problem with theoretical guarantees is the (inexact) proximal point method, where a quadratic proximal term is added to objective and constraint functions to construct a strongly convex constrained subproblem and then a sequence of solutions can be generated by solving this subproblem inexactly and updating the center of the proximal term. Oracle complexity for this method to find a nearly $\epsilon$-stationary has been established by \cite{boob2022stochastic, ma2020quadratically, jia2022first} under different constraint qualifications. 

The inexact proximal point method is a \emph{double-loop} algorithm where the inner loop is another optimization algorithm for solving the aforementioned strongly convex subproblems. The complexity results in \cite{boob2022stochastic, ma2020quadratically, jia2022first} require the inner loop solves each subproblem to a targeted optimality gap. However, the optimality gap is hard to evaluate and thus cannot be used to terminate the inner loop. Although the number of inner iterations needed to achieve the targeted gap can be bounded theoretically, the bound usually involves some constants that are unknown or can only be estimated conservatively. Hence, using the theoretical iteration bound to stop the inner loop usually leads to significantly more inner iterations than needed, making the whole algorithm inefficient. In practices, users often need to tune the number of inner iterations to improve algorithm's efficiency, which is an inconvenience common to most double-loop methods. 

On the contrary, a \emph{single-loop} algorithm is usually easier to implement as it does not require tuning the number of inner iterations. Therefore, the \textbf{main contribution} of this paper is showing that a single-loop first-order algorithm can find a nearly $\epsilon$-stationary point of \eqref{eq:gco} with oracle complexity $O(1/\epsilon^4)$, which matches state-of-the-art results obtained only by the double-loop methods~\cite{boob2022stochastic, ma2020quadratically, jia2022first}. The algorithm we study is the classical \emph{switching subgradient} (SSG) method proposed by Polyak \cite{swtichgradientpolyak} in 1967, which is intuitive and easy to implement but has only been analyzed before in the convex case. We first show that the SSG method has complexity  complexity $O(1/\epsilon^4)$ when $g$ is convex. Then, we show that the same complexity result holds when $g$ is weakly convex and appropriate constraint qualification holds. We also present some practical examples where all of our assumptions hold, including a fair classification problem subject to the demographic parity constraint~\cite{agarwal2018reductions}. Our \textbf{technical novelty} is inventing a \emph{switching stepzize rule} to accompany the switching subgradient. In particular, we use a fixed stepsize when the solution is updated by the objective's subgradient while use an adaptive Polyak's stepsize \cite{eremin1965relaxation,polyak1969minimization} when the solution is updated by the constraint's subgradient. This allows us to leverage the local error bound of the constraint function to keep the solution nearly feasible during the algorithm, which prevents the solution from being trapped at an infeasible stationry point. To the best of our knowledge, this paper is the first to establish the complexity of a single-loop first-order method for weakly convex non-smooth nonlinear constrained optimization. In the appendix, we also extend our algorithms and complexity analysis to the stochastic setting. 


\section{Related work}
 \label{sec:relatedwork}
Non-convex constrained optimization has a long history~\cite{fletcher19851,burke1989sequential,di1988exactness,cartis2011evaluation,di1989exact,auslender2013extended,cartis2019optimality} and the interest on this subject has been still growing in the machine learning community because of its new applications such as learning with fairness constraints (see e.g., \cite{zafar2019fairness}). 

In recent literature, the prevalent classes of algorithms for non-convex constrained optimization include the augmented Lagrangian method (ALM), the penalty method~\cite{gao2020admm, zhang2020proximal, zhang2022global, zhang2022iteration, hong2016convergence, kong2019complexity, melo2020iteration, kong2023iteration, kong2023accelerated, sujanani2023adaptive, kong2022iteration, li2021augmented, sahin2019inexact, lin2022complexity, li2021rate,jin2022stochastic,li2022zeroth}, and the sequential quadratic programming method~\cite{bolte2016majorization,facchinei2021ghost, boob2022level, curtis2021inexact, berahas2023sequential,berahas2021sequential,berahas2021stochastic,curtis2021worst}. Besides, an inexact projected gradient method is developed by \cite{boroun2021inexact} and a level conditional gradient method is developed by \cite{cheng2022functional}. However, these works all focus on the case where $g$ is smooth and $f$ is either smooth or equals $f_1+f_2$, where $f_1$ is smooth and non-convex while $f_2$ is non-smooth and convex and has a simple structure that allows a closed-form solution for the proximal mapping $\argmin_{\by}\{f_2(\by)+\frac{\rho}{2}\|\by-\bx\|^2\}$ for any $\bx$. There are relatively fewer works on non-convex non-smooth constrained problems. An alternating direction method of multipliers (ADMM) and an ALM are studied by \cite{wang2019global} and \cite{zeng2022moreau}, respectively, for non-convex non-smooth problems with linear constraints while our study considers nonlinear non-smooth constraints. The methods by \cite{cheng2022functional} and \cite{boob2022level} can be extended to a structured non-smooth case where $f=f_1+f_2$ with $f_1$ being smooth non-convex and $f_2=\max_{\by}\{\by^\top A\bx-\phi(\by)\}$ with a convex $\phi$, and $g$ has a similar structure. The method by \cite{boob2020feasible} can handle a specific non-smooth non-convex constraint, i.e., $g(\bx)=\lambda \|\bx\|_1-h(\bx)$ where $h$ is a convex and smooth. Compared to these works, our results apply to a more general non-smooth problem without those structure assumptions. 



When $f$ and $g$ in \eqref{eq:gco} are weakly convex and non-smooth, the inexact proximal point method has been studied by \cite{boob2022stochastic, ma2020quadratically, jia2022first} under different constraint qualifications and different notions of stationarity. Their complexity analysis utilizes the relationship between the gradient of the Moreau envelope of \eqref{eq:gco} and the near stationarity of a solution, which is originally used to analyze complexity of subgradient methods for weakly convex non-smooth unconstrained problems~\cite{Davis2018,davis2019stochastic,davis2019proximally,alacaoglu2021convergence,deng2021minibatch,rafique2022weakly,zhang2018convergence}. Our analysis utilizes a similar framework. The methods \cite{boob2022stochastic, ma2020quadratically, jia2022first} are double-loop while our algorithm only uses a single loop and achieves the same complexity of $O(1/\epsilon^4)$ as them under similar assumptions. 

The SSG algorithm is first proposed by Polyak \cite{swtichgradientpolyak}. It has been well-studied for convex problems~\cite{nesterov2018lectures,bayandina2018mirror,lan2020algorithms,stonyakin2019mirror,titov2020analogues,titov2019mirror,stonyakin2019adaptive,stonyakin2018one,stonyakin2019some,alkousa2020modification} and quasi-convex problems~\cite{stonyakin2019mirror}.  This paper provides the first complexity analysis for the SSG method under weak convexity assumption. Non-smooth non-convex optimization has also been studied without weak convexity assumption by~\cite{zhang2020complexity,kornowski2021oracle,shamir2020can,kornowski2022complexity,chen2023faster,tian2021computing,tian2022finite}. These works analyze the complexity of first-order methods for computing an $(\epsilon,\delta)$-Goldstein approximate stationary point, which is a more general stationarity notation than what we consider here. However, these works only focus on unconstrained problems. 

\section{Algorithm and near stationarity}\label{sec:pre}
Let $\|\cdot\|$ be the $\ell_2$-norm. For $h:\reals^d\rightarrow \reals\cup\{+\infty\}$, the subdifferential of $h$ at $\bx$ is
\small
\begin{align*}
	\partial h(\bx)=
	&~\big\{\bzt\in\mathbb{R}^d~\big|~ 
	h(\bx')
	\geq h(\bx)+\langle\bzt,\bx'-\bx\rangle +o(\|\bx'-\bx\|), ~\bx'\rightarrow\bx
	\big\},
\end{align*}
\normalsize
and $\bzt\in\partial h(\bx)$ is a subgradient of $h$ at $\bx$. 
We say  $h$ is \emph{$\mu$-strongly convex ($\mu\geq0$)} on $\X$ if 
\begin{align*}
	h(\bx)\geq h(\bx')+\langle\bzt,\bx-\bx'\rangle+\frac{\mu }{2}\|\bx-\bx'\|^2
\end{align*}
for any $(\bx,\bx')\in\X\times\X$ and any $\bzt\in\partial h(\bx')$. We say $h$ is \emph{$\rho$-weakly convex} ($\rho\geq0$) on $\X$ if
\begin{align*}
	h(\bx)\geq h(\bx')+\langle\bzt,\bx-\bx'\rangle-\frac{\rho}{2}\|\bx-\bx'\|^2
\end{align*}
for any $(\bx,\bx')\in\X\times\X$ and any $\bzt\in\partial h(\bx')$. 
We denote the normal cone of $\X$ at $\bx$ by $\mathcal{N}_\X(\bx)$ and the relative interior of $\X$ by $\text{relint}(\X)$.  We say a point $\bx$ is $\epsilon$-feasible if $\bx\in\mathcal{X}$ and $g(\bx)\leq\epsilon$. Let $\delta_{\X}(\bx)$ be the zero-infinity characteristic function of set $\X$, $\text{proj}_{\X}(\cdot)$ be the projection mapping to $\X$, and $\text{dist}(\bx,\mathcal{A}):=\min_{\by\in \mathcal{A}}\|\bx-\by\|$ for set $\mathcal{A}$.
 
	We make the following assumptions on \eqref{eq:gco} throughout the paper.
	\begin{assumption}
		\label{assume:allpaper}
  Functions $f$ and $g$ are continuous with $\partial f(\bx)\neq\emptyset$ and $\partial g(\bx)\neq\emptyset$ on $\X$, and there exists $M$ such that $\|\bzt_f\|\leq M$ and $\|\bzt_g\|\leq M$ for any $\bx\in\X$, $\bzt_f\in\partial f(\bx)$ and $\bzt_g\in\partial g(\bx)$.
	\end{assumption}

	Since \eqref{eq:gco} is non-convex, finding an $\epsilon$-optimal solution is intractable in general. For a non-convex problem, the target is typically to find a \emph{stationary} point of \eqref{eq:gco}, which is a point $\bx^*\in\mathcal{X}$ that satisfies the following Karush-Kuhn-Tucker (KKT) conditions:
	\begin{equation*}
	\bzt_f^*+\lambda^*\bzt_g^*\in-\mathcal{N}_\X(\bx^*), 
	\quad \lambda^*g(\bx^*)=0, \quad g(\bx^*)\leq0,\quad\lambda^*\geq0,
	\end{equation*}
	where $\lambda^*\in\mathbb{R}$ is a Lagrangian multiplier, $\bzt_f^*\in\partial f(\bx^*)$ and $\bzt_g^*\in\partial g(\bx^*)$.  Typically, an exact stationary point can only be approached by an algorithm at full convergence, which may require infinitely many iterations. Within a finite number of iterations, an algorithm can only generate an \emph{$\epsilon$-stationary} point~\cite{boob2022stochastic}, which is a point $\widehat\bx\in\mathcal{X}$ satisfying
	\begin{equation}
	\label{eq:eKKT}
	\text{dist}\left(\widehat\bzt_f+ \widehat\lambda \widehat\bzt_g,-\mathcal{N}_\X(\widehat\bx)\right)\leq \epsilon, \quad
	 |\widehat\lambda g(\widehat\bx)|\leq \epsilon^2,\quad g(\widehat\bx)\leq\epsilon^2,\quad\widehat\lambda\geq0,
	\end{equation}
	where $\widehat\lambda\in\mathbb{R}$ is a Lagrangian multiplier, $\widehat\bzt_f\in\partial f(\widehat\bx)$ and  $\widehat\bzt_g\in\partial g(\widehat\bx)$. However, because $f$ and $g$ are non-smooth, computing an $\epsilon$-stationary point is still challenging even for an unconstrained problem. Nevertheless, under the weak convexity assumption, it is possible to compute a \emph{nearly $\epsilon$-stationary point}, which we will introduce next.

 Given $\hat\rho\geq0$, $\tilde\rho\geq0$ and $\bx\in\X$, we define a quadratically regularized problem of \eqref{eq:gco} as
	\small
	\vspace{-0.3em}
	\begin{equation}
	\begin{aligned}
	\label{eq:phi}
 \varphi(\bx)
 \equiv \min_{\by\in\X} \Big\{ f(\by) + \frac{\hat{\rho}}{2}\|\by - \bx\|^2, 
	~{s.t.}~ g(\by)+\frac{\tilde{\rho}}{2}\|\by - \bx\|^2\leq 0 \Big\} ,\\
	\end{aligned}
	\end{equation}
	\vspace{-0.5em}
	\begin{equation}
	\begin{aligned}
	\label{eq:phix}
 \widehat\bx(\bx)
 \equiv \argmin_{\by\in\X} \Big\{ f(\by) + \frac{\hat{\rho}}{2}\|\by - \bx\|^2, 
	~{s.t.}~ g(\by)+\frac{\tilde{\rho}}{2}\|\by - \bx\|^2\leq 0 \Big\}.
	\end{aligned}
	\end{equation}
	\normalsize
Following the literature on weakly convex optimization~\cite{Davis2018,davis2019proximally,davis2018complexity,boob2022stochastic, ma2020quadratically, jia2022first}, we use the value of $\|\widehat\bx(\bx)-\bx\|$ as a measure of the quality of a solution $\bx$ because it can be interpreted as a stationarity measure. For the purpose of illustration, we assume for now that $\widehat\bx(\bx)$ is uniquely defined and there exists a Lagrangian multiplier $\widehat\lambda\in\mathbb{R}$ such that the following KKT conditions of \eqref{eq:phix} holds.  
	\begin{equation}
	\begin{aligned}
	\label{eq:KKTprox}
	&\widehat\bzt_f  +\hat{\rho}(\widehat\bx (\bx) - \bx)+\widehat\lambda\left(\widehat\bzt_g+\tilde{\rho}(\widehat\bx (\bx) - \bx)\right)\in-\mathcal{N}_\X(\widehat\bx (\bx)),\\
	& \widehat\lambda \left(g(\widehat\bx(\bx))+\frac{\tilde{\rho}}{2}\|\widehat\bx(\bx) - \bx\|^2\right)=0, \quad g(\widehat\bx (\bx))+\frac{\tilde{\rho}}{2}\|\widehat\bx(\bx) - \bx\|^2\leq0, \quad \widehat\lambda\geq0,
	\end{aligned}
	\end{equation}
	where $\widehat\bzt_f\in\partial f(\widehat\bx(\bx))$ and $\widehat\bzt_g\in\partial g(\widehat\bx(\bx))$.	
Therefore, as long as $\|\widehat\bx (\bx) - \bx\|\leq \epsilon$ , we have
\begin{align}
\label{eq:nearKKTbyhat}
\text{dist}\left(\widehat\bzt_f+\widehat\lambda \widehat\bzt_g, -\mathcal{N}_\X(\widehat\bx (\bx)) \right)\leq (\hat{\rho}+\widehat\lambda\tilde{\rho})\epsilon,\quad |\widehat\lambda g(\widehat\bx (\bx))|=\widehat\lambda\tilde\rho\epsilon^2/2,
\quad g(\widehat\bx (\bx))\leq 0.
\end{align}
This means $\widehat\bx (\bx)$ is an $O(\epsilon)$-stationary point of the original problem \eqref{eq:gco} in the sense of \eqref{eq:eKKT}. Since $\bx$ is within an $\epsilon$-distance from $\widehat\bx (\bx)$, we call such an $\bx$ a \emph{nearly $\epsilon$-stationary point} of \eqref{eq:gco}. We formalize this definition as follows.


\begin{definition}
\label{def:stationary}
Suppose $\widehat\bx(\bx)$ is defined in \eqref{eq:phix} with $\epsilon\geq0$. 
A (stochastic) point $\bx\in\X$ is a (stochastic) nearly $\epsilon$-stationary point of \eqref{eq:gco} if 
$\mathbb{E}[\|\widehat\bx(\bx)  - \bx\|]\leq\epsilon$. 
\end{definition}
Of course, we can claim $\widehat\bx (\bx)$ is an $O(\epsilon)$-stationary point of \eqref{eq:gco} based on \eqref{eq:nearKKTbyhat}  only when $\widehat\lambda$ in \eqref{eq:KKTprox} exists and does not go to infinity as $\epsilon$ approaches zero. Fortunately, we can show in Lemmas~\ref{thm:boundlambda} and~\ref{thm:boundlambda_wc} that this is true under some constraint qualifications, which justifies Definition~\ref{def:stationary}. 

\begin{algorithm}[tb]
   \caption{Switching Subgradient (SSG) Method by Polyak~\cite{swtichgradientpolyak}}
   \label{alg:dsgm}
\begin{algorithmic}[1]
   \STATE {\bfseries Input:} $\bx^{(0)}\in \X$, total number of iterations $T$, tolerance of infeasibility $\epsilon_t\geq0$, stepsize $\eta_t$, and a starting index $S$ for generating outputs.
   \STATE {\bfseries Initialization:} $I=\emptyset$ and $J=\emptyset$.
   \FOR{iteration $t=0,1,\cdots,T-1$}
   \IF{$g (\bx^{(t)} )\leq\epsilon_t$}
    \STATE Set $\bx^{(t+1)}= \text{proj}_{\X} (\bx^{(t)}-\eta_{t}\bzt_f^{(t)} )$ for any $\bzt_f^{(t)}\in\partial f(\bx^{(t)})$
    and, if $t\geq S$,  $I= I\cup\{t\}$.
    \ELSE
    \STATE Set $\bx^{(t+1)}= \text{proj}_{\X} (\bx^{(t)}-\eta_{t}\bzt_g^{(t)} )$ for any $\bzt_g^{(t)}\in\partial g(\bx^{(t)})$
    \ and, if $t\geq S$,   $J= J\cup\{t\}$.
   \ENDIF
   \ENDFOR
\STATE {\bfseries Output I:} $\bx^{(\tau)}$ with $\tau$ randomly sampled from $I$ using $\text{Prob}(\tau=t)=\eta_t/\sum_{s\in I}\eta_s$.
\STATE {\bfseries Output II:} $\bx^{(\tau)}$ with $\tau$ randomly sampled from  $I\cup J$ using $\text{Prob}(\tau=t)=\eta_t/\sum_{s\in I\cup J}\eta_s$.
\end{algorithmic}
\end{algorithm}

We present the SSG method in Algorithm~\ref{alg:dsgm} for finding a nearly $\epsilon$-stationary point of \eqref{eq:gco}. At iteration $t$, we check if the current solution $\bx^{(t)}$ is nearly feasible in the sense that $g(\bx^{(t)})\leq \epsilon_t$ for a pre-determined tolerance of infeasibility $\epsilon_t$. If yes, the algorithm performs a subgradient step along the subgradient of $f$. Otherwise, the algorithm switches the updating direction to the subgradient of $g$. The algorithm records the iteration indexes of the nearly feasible solutions in set $I$ and other indexes in set $J$. The final output is randomly sampled from the iterates in $I$ or $I\cup J$ with a distribution weighted by the stepsizes $\eta_t$'s. An index $S$ is set so the algorithm only starts to record $I$ and $J$ when $t\geq S$.  We study the theoretical \emph{oracle complexity} of Algorithm~\ref{alg:dsgm} for finding a nearly $\epsilon$-stationary point, which is defined as the total number of times for which the algorithm queries the subgradient or function value of $f$ or $g$. Our results are presented separately when $g$ is convex and when $g$ is weakly convex. 


\section{Convex constraints}
\label{sec:deterministic}
In this section, we first consider a relatively easy case where $f$ is weakly convex but $g$ is convex. In particular, we make the following assumptions in addition to Assumption~\ref{assume:allpaper} in this section.
\begin{assumption}
		\label{assume:convex}
		The following statements hold:
		\begin{itemize}
			\item[A.] $f(\bx)$ is $\rho$-weakly convex on $\X$ and $g(\bx)$ is convex on $\X$. 
			\item[B.] (Slater's condition) There exists $\bx_{\text{feas}}\in\text{relint}(\X)$ such that $g(\bx_{\text{feas}})<0$.
                \item[C.] There exists $D\in\mathbb{R}$ such that $\|\bx-\bx'\|\leq D$ for any $\bx$ and $\bx'$ in $\X$.
		\end{itemize}
	\end{assumption}

In this section, we  choose parameters in \eqref{eq:phi} such that 
\begin{equation}
\label{eq:parameter1}
\hat\rho>\rho\text{ and }\tilde\rho=0.
\end{equation}
Under Assumption~\ref{assume:convex}, \eqref{eq:parameter1} guarantees that \eqref{eq:phi} is strictly feasible, its  objective function is $(\hat\rho-\rho)$-strongly convex and its constraint function is convex, so $\widehat\bx(\bx)$ in \eqref{eq:phix} is unique and $\widehat\lambda$ in \eqref{eq:KKTprox} exists. We first present an upper bound of $\widehat\lambda$ that is independent of $\bx$. The proof is in Section~\ref{sec:boundlambda}.
\begin{lemma}
\label{thm:boundlambda}
Suppose Assumptions~\ref{assume:allpaper} and \ref{assume:convex} hold. Given any $\bx\in\X$, let $\widehat\bx(\bx)$ be defined as in \eqref{eq:phix} with $(\hat\rho,\tilde\rho)$ satisfying \eqref{eq:parameter1} and $\widehat\lambda$ be the associated Lagrangian multiplier satisfying \eqref{eq:KKTprox}. We have 
\begin{eqnarray}
\label{eq:Lambdabound}
\widehat\lambda \leq  \Lambda:=(MD+\hat\rho D^2)/(-g(\bx_{\textup{feas}})).
\end{eqnarray}
\end{lemma}

For simplicity of notation, we denote $\widehat\bx(\bx^{(t)})$ defined in \eqref{eq:phix} by $\widehat\bx^{(t)}$.
Let $\mathbb{E}_\tau[\cdot]$ be the expectation taken only over the random index $\tau$ when the algorithms stop.  
We present the convergence properties of Algorithm~\ref{alg:dsgm} when $\epsilon_t$ and $\eta_t$ are static and diminishing. 
The proof is provided in Section~\ref{sec:mdconverge_deterministic}.
\begin{theorem}
\label{thm:mdconverge_deterministic}
Suppose Assumptions~\ref{assume:allpaper} and \ref{assume:convex} hold and $\Lambda$ is as in \eqref{eq:Lambdabound}.  Let $\widehat\bx(\bx^{(t)})$ be defined as in \eqref{eq:phix} with $(\hat\rho,\tilde\rho)$ satisfying \eqref{eq:parameter1}
and $\bx^{(\tau)}$ is generated by Output I. Algorithm~\ref{alg:dsgm} guarantees 
$\mathbb{E}_\tau[\|\widehat\bx^{(\tau)} - \bx^{(\tau)}\|]\leq\epsilon$ and $\mathbb{E}_\tau [g(\bx^{(\tau)})]\leq \frac{\epsilon^2(\hat\rho-\rho)}{1+\Lambda}$ in either of the following cases. 

Case I: $S=0$, $\epsilon_t=\frac{\epsilon^2(\hat\rho-\rho)}{1+\Lambda}$, $\eta_t=\frac{2\epsilon^2(\hat\rho-\rho)}{5(1+\Lambda)M^2}$ and 
$T\geq \frac{25M^2D^2(1+\Lambda)^2}{4\epsilon^4(\hat\rho-\rho)^2}
=O(1/\epsilon^4)$.

Case II: $S=T/2$, $\epsilon_t=\frac{5MD}{\sqrt{t+1}}$, $\eta_t=\frac{D}{M\sqrt{t+1}}$ and 
$T\geq  \frac{50M^2D^2(1+\Lambda)^2}{\epsilon^4(\hat\rho-\rho)^2}   
=O(1/\epsilon^4)$.
\end{theorem}

Algorithm~\ref{alg:dsgm} is single-loop with $O(1)$ oracle complexity per iteration, so its total complexity is just $T=O(1/\epsilon^4)$, which matches the start-of-the-art complexity by~\cite{boob2022stochastic, ma2020quadratically, jia2022first}. This result can be generalized to the case where there exist additional linear equality constraints $\bA\bx=\bb$. See the extension of Lemma~\ref{thm:boundlambda} in Section~\ref{sec:boundlambda} and Remark~\ref{remark:complexity_ec}. 

Assumption~\ref{assume:convex}C can be relaxed to only require that the feasible set, i.e., $\mathcal{S}=\{\bx\in\X~|~g(\bx)\leq 0\}$, is bounded instead of $\X$. The same complexity can be achieved by Algorithm~\ref{alg:dsgm} by using a switching step size rule similar to $\eta_t$ in Proposition~\ref{thm:feasiblesubproblem}. This result is provided in Section~\ref{sec:convexgsharpness} but we recommend interested readers to read Section~\ref{sec:weaklyconvexconstraints} first to get introduced to this special step size.

\begin{remark}
\label{remark:feasibility}
Property $\mathbb{E}_\tau [g(\bx^{(\tau)})]\leq \frac{\epsilon^2(\hat\rho-\rho)}{1+\Lambda}$ in the theorems above is not required by Definition~\ref{def:stationary}. By Assumption~\ref{assume:allpaper}A, 
$
\mathbb{E}_\tau [g(\bx^{(\tau)})]\leq \mathbb{E}_\tau[ g(\widehat\bx^{(\tau)})]+M\mathbb{E}_\tau[\|\widehat\bx^{(\tau)}-\bx^{(\tau)}\|]
\leq M\epsilon,
$
which means a nearly $\epsilon$-stationary point must be $O(\epsilon)$-feasible by definition. Property $\mathbb{E}_\tau [g(\bx^{(\tau)})]\leq \frac{\epsilon^2(\hat\rho-\rho)}{1+\Lambda}$ implies $O(\epsilon^2)$-feasibility for the output, which is even better. 
\end{remark}

When $g$ is $\mu$-strongly convex with $\mu>0$, we can show that the complexity of Algorithm~\ref{alg:dsgm} is still $O(1/\epsilon^4)$ but one can simply set $\epsilon_t=0$, which makes $\eta_t$ the only tuning parameter. This makes this single-loop method even more attractive. Due to space limit, we include this result in Section~\ref{sec:strongconvexg}. We also extend our result to the stochastic case in Section~\ref{sec:stochastic} and \ref{sec:mdconverge_prob}.

\section{Weakly convex constraints}
\label{sec:weaklyconvexconstraints}
Next we consider the case where both $f$ and $g$ are weakly convex but not necessarily convex. Let 
\begin{align*}
    g_+(\bx)=\max\{g(\bx),0\},\quad \mathcal{L}=\{\bx\in\X~|~g(\bx)=0\}~\text{ and }~\mathcal{S}=\{\bx\in\X~|~g(\bx)\leq 0\}.  
\end{align*}
We make the following assumptions in addition to Assumption~\ref{assume:allpaper} in this section.
\begin{assumption}
		\label{assume:weaklyconvex}
		The following statements hold:
		\begin{itemize}
			\item[A.] $f(\bx)$ and $g(\bx)$  are $\rho$-weakly convex on $\X$. 
			\item[B.] There exist $\bar\epsilon>0$, $\theta>0$ and $\bar\rho>\rho$ such that, for any $\bar\epsilon^2$-feasible solution $\bx$, there exists $\by\in\text{relint}(\X)$ such that $g(\by)+\frac{\bar\rho}{2}\|\by-\bx\|^2\leq-\theta$. 
            \item[C.] $\underline f:=\inf_{\bx\in\X}f(\bx)>-\infty$.
		\end{itemize}
	\end{assumption}

Assumption~\ref{assume:weaklyconvex}B is called the uniform Slater's condition by \cite{ma2020quadratically}. We will present two real-world examples in Section~\ref{sec:slaterexample} that satisfy this assumption, including  a fair classification problem under demographic parity constraint, which is one of the applications in our numerical experiments in Section~\ref{sec:classification_demographic_parity}. In this section, we choose parameters in \eqref{eq:phi} such that
\begin{equation}
\label{eq:parameter2}
\bar\rho\geq\hat\rho=\tilde\rho>\rho.
\end{equation}
Under Assumption~\ref{assume:weaklyconvex}, \eqref{eq:parameter2} guarantees that \eqref{eq:phi} is uniformly strictly feasible for any $\epsilon^2$-feasible $\bx$, and the objective and constraint functions of  \eqref{eq:phi} are both $(\hat\rho-\rho)$-strongly convex, so $\widehat\bx(\bx)$ is uniquely defined and $\widehat\lambda$ in \eqref{eq:KKTprox} exists. In addition, Assumption~\ref{assume:weaklyconvex} has the following three implications that play important roles in our complexity analysis. 

First, $\widehat\lambda$ in \eqref{eq:KKTprox} can be bounded by a constant independent of $\bx$ and $\epsilon$ as long as $\bx$ is $\epsilon^2$-feasible with $\epsilon\leq \bar\epsilon$. This result is similar to Lemma 1 by Ma et al. \cite{ma2020quadratically} except that they require $\X$ to be bounded but we do not. The proof is given in Section~\ref{sec:boundlambda_wc}.
\begin{lemma}
\label{thm:boundlambda_wc}
Suppose Assumptions~\ref{assume:allpaper} and \ref{assume:weaklyconvex} hold. Given any $\epsilon^2$-feasible $\bx$ with any $\epsilon\leq \bar\epsilon$, let $\widehat\bx(\bx)$ defined as in \eqref{eq:phix} satisfying \eqref{eq:parameter2} and $\widehat\lambda$ is the associated Lagrangian multiplier satisfying \eqref{eq:KKTprox}. We have 
\begin{eqnarray}
\label{eq:Lambdabounda}
\|\widehat\bx(\bx)-\bx\|\leq M/\hat\rho\quad\text{ and }\quad
\widehat\lambda\leq  \Lambda':=2M/\sqrt{2 \theta (\hat{\rho} - \rho)}.
\end{eqnarray}
\end{lemma}

Second, the subgradient of the constraint function $g(\bx)+\delta_{\X}(\bx)$ on $\mathcal{L}$ is uniformly away from the origin. The proof is provided in Section~\ref{sec:boundlambda_wc}.
\begin{lemma}
\label{thm:sharpsubgradient}
Suppose Assumptions~\ref{assume:allpaper} and \ref{assume:weaklyconvex} hold. It holds for any $\bx\in \mathcal{L}$ that
\begin{align}
\label{eq:sharpsubgradient}
\min\limits_{\bzt_g\in\partial g(\bx),\bu\in\mathcal{N}_{\X}(\bx)}\|\bzt_g+\bu\|\geq\nu:=\sqrt{2 \theta (\hat{\rho} - \rho)}.
\end{align}
\end{lemma}
Lastly, note that $\mathcal{S}$ is the optimal set of $\min_{\bx\in\X}g_+(\bx)$, which is a $\rho$-weakly convex non-smooth optimization problem with an optimal value of zero. Lemma~\ref{thm:sharpsubgradient} implies that $g_+(\bx)$ is sharp near the boundary of $\mathcal{S}$, meaning that $g_+(\bx)$ satisfies a linear error bound in an $O(1)$-neighborhood of $\mathcal{S}$. 
A similar result for a convex $g_+$ is given in Lemma 1 in \cite{yang2018rsg}. In the lemma below, we extend their result for a weakly convex $g_+$. The proof is in Section~\ref{sec:boundlambda_wc} and the second conclusion is directly from~\cite{davis2018subgradient}. 


\begin{lemma}
\label{thm:errorbound_wc}
Suppose Assumptions~\ref{assume:allpaper} and \ref{assume:weaklyconvex} hold. It holds
for any $\bx$ satisfying $\textup{dist}(\bx,\mathcal{S})\leq \frac{\nu}{\rho}$  that
\begin{eqnarray}
\label{eq:errorbound_wc}
(\nu/2)\textup{dist}(\bx,\mathcal{S})\leq g_+(\bx).
\end{eqnarray}
Moreover, $\nu\leq 2M$ and $\min_{\bx\in\X}g_+(\bx)$ has no stationary point satisfying $0<\textup{dist}(\bx,\mathcal{S})<\frac{\nu}{\rho}$.
\end{lemma}

Since $g$ is non-convex, Algorithm~\ref{alg:dsgm} may not even find a nearly feasible solution if $\bx^{(t)}$ is trapped at a stationary point of $g$ with $g(\bx)>0$, that is, a sub-optimal stationary point of $\min_{\bx\in\X}g_+(\bx)$. Fortunately, the second conclusion of Lemma~\ref{thm:errorbound_wc} indicates that this situation can be avoided by keeping $\textup{dist}(\bx^{(t)},\mathcal{S})=O(\epsilon^2)$ during the algorithm. To do so, we start with $\bx^{(0)}\in\mathcal{S}$ and use $\epsilon_t=O(\epsilon^2)$ in Algorithm~\ref{alg:dsgm}. Moreover, we apply a switching stepsize rule that sets $\eta_t=O(\epsilon^2)$ when $t\in I$ and $\eta_t=g(\bx^{(t)})/\|\bzt_g^{(t)}\|^2$ when $t\in J$, the latter of which is known by the Polyak's stepsize~\cite{eremin1965relaxation,polyak1969minimization}. This way, when $g (\bx^{(t)} )\leq\epsilon_t$, \eqref{eq:errorbound_wc} ensures $\textup{dist}(\bx^{(t)},\mathcal{S})=O(\epsilon^2)$. When $g (\bx^{(t)} )>\epsilon_t$, \eqref{eq:errorbound_wc} ensures $\text{dist}(\bx^{(t)},\mathcal{S})$ Q-linearly converges to zero \cite{davis2018subgradient}, which also guarantees $\textup{dist}(\bx^{(t)},\mathcal{S})=O(\epsilon^2)$. As a result, we have $g(\bx^{(t)})\leq\epsilon^2$ for any $t$, so problem \eqref{eq:phi} with $\bx=\bx^{(t)}$ and $(\hat\rho,\tilde\rho)$ satisfying \eqref{eq:parameter2} will be strictly feasible according to Assumption~\ref{assume:weaklyconvex}B. This finding is given in the proposition below with its proof in Section~\ref{sec:proof_mainresult_wc}.



\begin{proposition}
\label{thm:feasiblesubproblem}
Suppose Assumptions~\ref{assume:allpaper} and \ref{assume:weaklyconvex} hold and $\epsilon\leq \bar\epsilon$. Also, suppose $\bx^{(t)}$ is generated by Algorithm~\ref{alg:dsgm} using 
$\bx^{(0)}\in\mathcal{S}$, $\epsilon_{t}=\frac{\nu}{4}\min\left\{\epsilon^2/M,\nu/(4\rho)\right\}$ and 
\begin{align*}
\eta_{t}=
\left\{
\begin{array}{ll}
\frac{\nu}{4M^2}\min\left\{\epsilon^2/M,\nu/(4\rho)\right\}&\text{ if }t\in I\\
g(\bx^{(t)})/\|\bzt_g^{(t)}\|^2&\text{ if }t\in J.
\end{array}
\right.
\end{align*}
Then $\textup{dist}(\bx^{(t)},\mathcal{S})\leq \min\left\{\epsilon^2/M,\nu/(4\rho)\right\}$ and $g(\bx^{(t)})\leq \epsilon^2$ for any $t\geq0$. As a consequence, $\bx^{(t)}$ is $\epsilon^2$-feasible to \eqref{eq:phi} where $\bx=\bx^{(t)}$ and $(\hat\rho,\tilde\rho)$ satisfies \eqref{eq:parameter2}. 
\end{proposition}

Let $\widehat\bx^{(t)}$ be $\widehat\bx(\bx^{(t)})$ defined in \eqref{eq:phix} with $(\hat\rho,\tilde\rho)$ satisfying \eqref{eq:parameter2}. The complexity result for the case of weakly convex constraints is as follows. The proof can be found in Section~\ref{sec:proof_mainresult_wc}. 
\begin{theorem}
\label{thm:mdconverge_deterministic_weaklyconvex}
Under the same assumptions as Proposition~\ref{thm:feasiblesubproblem}, Algorithm~\ref{alg:dsgm} guarantees 
$\mathbb{E}_\tau[\|\widehat\bx^{(\tau)} - \bx^{(\tau)}\|]\leq C\epsilon$ and $\mathbb{E}_\tau [g(\bx^{(\tau)})]\leq \epsilon^2$, where $C:= 2\sqrt{\frac{1+\Lambda'}{\hat{\rho}-\rho}}$, if we use Output II and set 
\small
\begin{align*}
    S=0\text{ and }T\geq \frac{8M^2\left(f(\bx^{(0)})-\underline f+3M^2/(2\hat\rho)\right)}{\hat\rho(1+\Lambda')\nu\epsilon^2\min\left\{\epsilon^2/M,\nu/(4\rho)\right\}}=O(1/\epsilon^4).
\end{align*}
\normalsize
\end{theorem}
This theorem indicates that Algorithm~\ref{alg:dsgm} finds a nearly $(C\epsilon)$-stationary point with complexity $O(1/\epsilon^4)$. To obtain a nearly $\epsilon$-stationary point with $\epsilon\leq \bar\epsilon$, one only needs to replace $\epsilon$ in $\eta_t$, $\epsilon_t$ and $T$ in this theorem above by $\epsilon/\max\{C, 1\}$ . This will only change the constant factor in the $O(1/\epsilon^4)$ complexity. This complexity matches the start-of-the-art complexity by~\cite{boob2022stochastic, ma2020quadratically, jia2022first}.

\section{Numerical experiments}
\label{sec:exp}
We demonstrate the performance of the SSG method on two fairness-aware classification problems, which are instances of \eqref{eq:gco} with convex and weakly convex $g$'s, respectively.  We compare with two different double-loop inexact proximal point (IPP) methods~\cite{boob2022stochastic, ma2020quadratically, jia2022first}. The IPP method approximately solves a strongly convex constrained subproblem in each outer iteration, and we use the SSG method in this paper and the ConEx method in \cite{boob2022stochastic} as the solvers (inner loop) because they both have the best theoretical complexity for that subproblem. We use IPP-SSG and IPP-ConEx to denote these two implementations of the IPP method. All experiments are conducted using MATLAB 2022b on a computer with the CPU 3.20GHz x Intel Core i7-8700 and 16GB memory. 


\subsection{Classification problem with ROC-based fairness}
Given a feature vector $\ba\in\mathbb{R}^d$ and a class label $b\in\{1,-1\}$, the goal in a binary linear classification problem is to learn a model $\bx\in\mathbb{R}^d$ to predict $b$ based on the score $\bx^\top\ba$. Let $\mathcal{D}=\{(\ba_i,b_i)\}_{i=1}^n$ be a training set and $\ell(\cdot)$ be a convex non-increasing loss function. Model $\bx$ can be learned by solving 
\begin{equation}
\label{eq:ermL}
L^*=\min_{\bx\in\mathcal{W}}\Big\{L(\bx):=\frac{1}{n}\sum_{i=1}^n \ell(b_i\bx^\top\ba_i)\Big\},
\end{equation}
where $\mathcal{W}=\{\bx\in\mathbb{R}^d~|~\|\bx\|\leq r\}$. Solving \eqref{eq:ermL} may ensure good classification performance of $\bx$ but not its fairness. Suppose there exist two additional datasets. One contains the feature vectors of a protected group, denoted by $\mathcal{D}_p=\{\ba_i^p\}_{i=1}^{n_p}$, and the other one contains the feature vectors of an unprotected group, denoted by $\mathcal{D}_u=\{\ba_i^u\}_{i=1}^{n_u}$. We want to enhance the fairness of $\bx$ between these two groups using the ROC-based fairness metric proposed by \cite{vogel2021learning}. Suppose we set a threshold $\theta$ and 
classify data $\ba$ as positive if $\bx^\top \ba\geq \theta$ and as negative otherwise. 
The ROC-based fairness measure and its continuous approximation are defined as
\begin{align}
\nonumber
&\max_{\theta\in\Theta}\Big|
    \frac{1}{n_p}\sum_{i=1}^{n_p} \mathbb{I}(\bx^\top\ba_i^p\geq\theta)-\frac{1}{n_u}\sum_{i=1}^{n_u} \mathbb{I}(\bx^\top\ba_i^u\geq\theta) \Big|\\\label{eq:ROC}
   \approx R(\bx):=&\max_{\theta\in\Theta}\Big|
\frac{1}{n_p}\sum_{i=1}^{n_p} \sigma(\bx^\top\ba_i^p-\theta)-\frac{1}{n_u}\sum_{i=1}^{n_u} \sigma(\bx^\top\ba_i^u-\theta)
\Big|,
\end{align}
where $\sigma(z)=\exp(z)/(1+\exp(z))$ is the sigmoid function and $\Theta$ is a finite set of thresholds. If the value of this measure is small, model $\bx$ produces similar predicted positive rates for the protected and unprotected groups on various $\theta$'s, indicating the fairness of the model. To obtain a fair $\bx$, we balance \eqref{eq:ermL} and \eqref{eq:ROC} by solving 
\begin{equation}
\label{eq:ROCfairnessClassification}
\min_{\bx\in\mathcal{W}} R(\bx)\text{ s.t. } L(\bx)\leq L^*+\kappa,
\end{equation}
where $\kappa$ is a slackness parameter indicating how much we are willing to increase the classification loss in order to reduce $R(\bx)$ to obtain a more fair model. Problem \eqref{eq:ROCfairnessClassification} is an instance of \eqref{eq:gco} satisfying Assumptions~\ref{assume:allpaper} and~\ref{assume:convex} with $\rho=\beta$ where $\beta$ is defined in \eqref{eq:hLip} in Section~\ref{sec:slaterexample}.

We solve \eqref{eq:ROCfairnessClassification} on three datasets: \textit{a9a}~\cite{kohavi1996scaling}, \textit{bank}~\cite{moro2014data} and \textit{COMPAS}~\cite{angwin2016compas}. The information of
these datasets is given in Table \ref{tbl:data}. We split each dataset into two subsets with a ratio of $2:1$. The larger set is used as $\mathcal{D}$ in the constraint and the smaller set is further split into $\mathcal{D}_p$ and $\mathcal{D}_u$ based on a binary group variable listed in Table \ref{tbl:data}.

\begin{table}[t]
\caption{Information of the datasets. Groups are males VS females in a9a, users with age within $[25,60]$ VS outside $[25,60]$ in bank, and caucasian VS non-caucasian in COMPAS.}
\label{tbl:data}
\vskip 0.1in
\begin{center}
\begin{small}
\begin{sc}
\begin{tabular}{c|ccccc}
\toprule
 Datasets & $n$ & $d$ & Label & Groups \\
\midrule
a9a & 48,842 & 123 & Income & Gender\\
Bank & 41,188 & 54 & Subscription & Age\\
COMPAS & 6,172 & 16 & Recidivism & Race\\
\bottomrule
\end{tabular}
\end{sc}
\end{small}
\end{center}
\vskip -0.1in
\end{table}

\begin{figure*}
     \begin{tabular}[h]{@{}c|ccc@{}}
      & a9a & bank & COMPAS \\
		\hline \vspace*{-0.1in}\\
		\raisebox{12ex}{\small{\rotatebox[origin=c]{90}{Objective}}}
		& \hspace*{-0.06in}\includegraphics[width=0.30\textwidth]{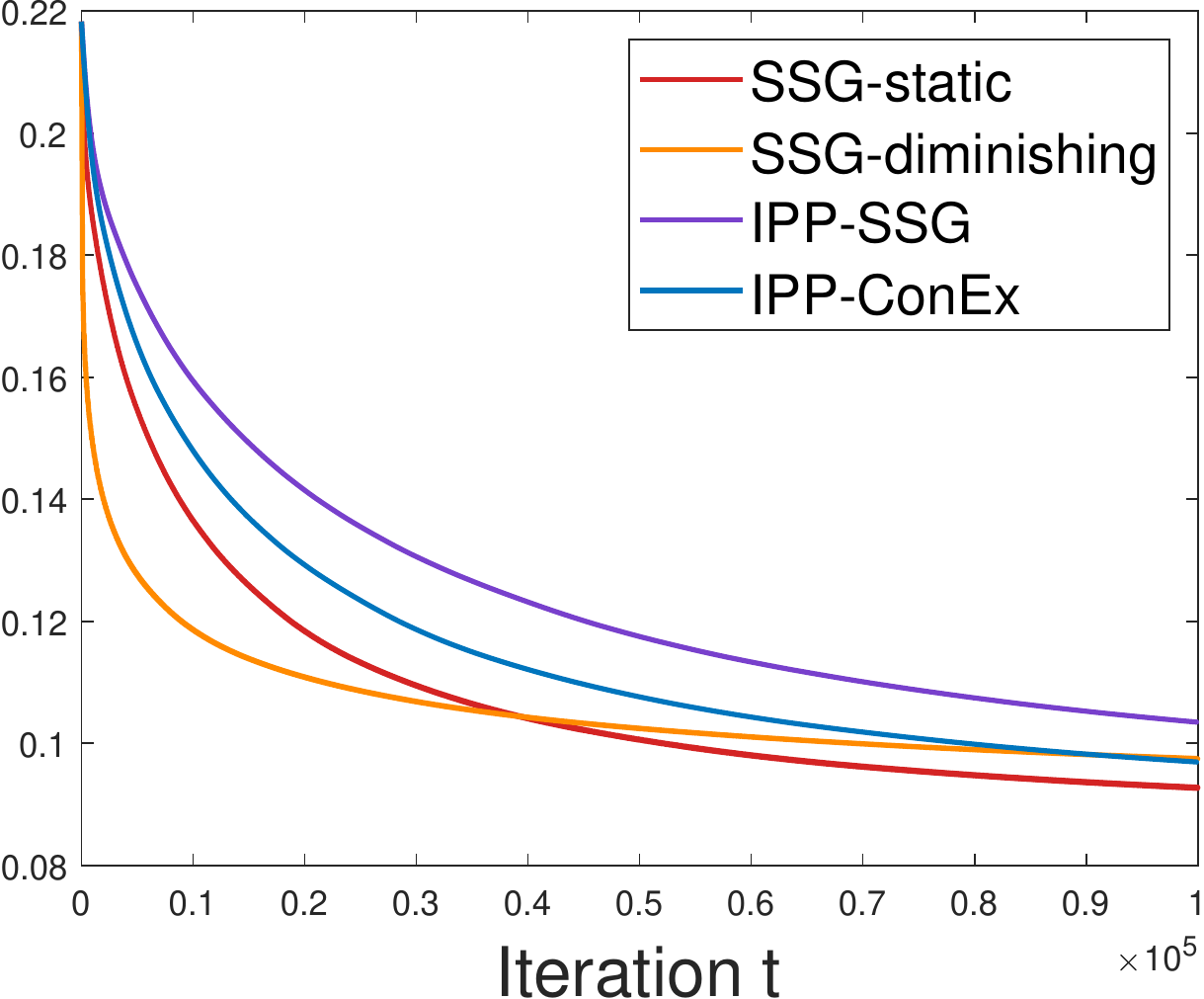}
		& \hspace*{-0.06in}\includegraphics[width=0.30\textwidth]{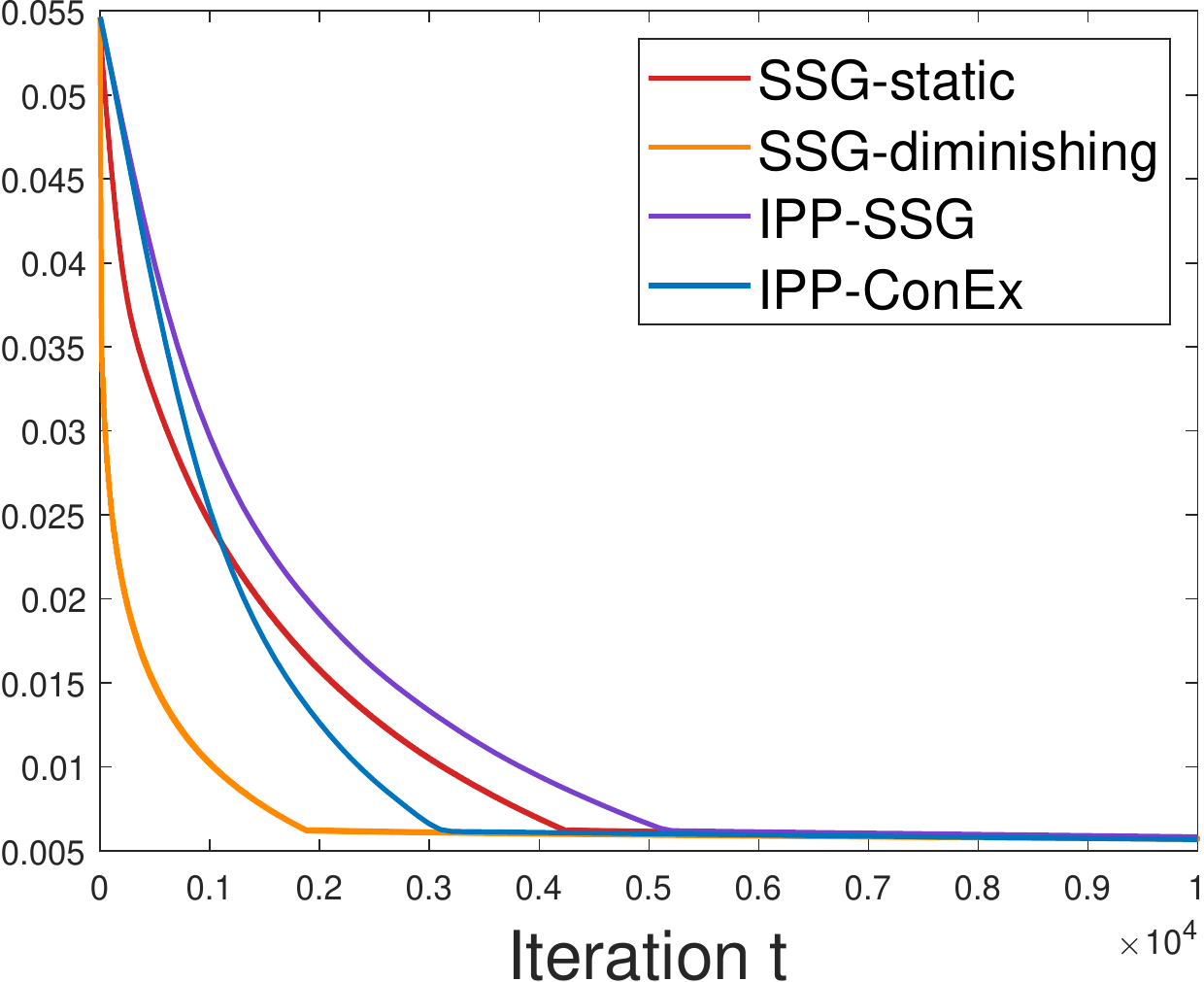}
		& \hspace*{-0.06in}\includegraphics[width=0.30\textwidth]{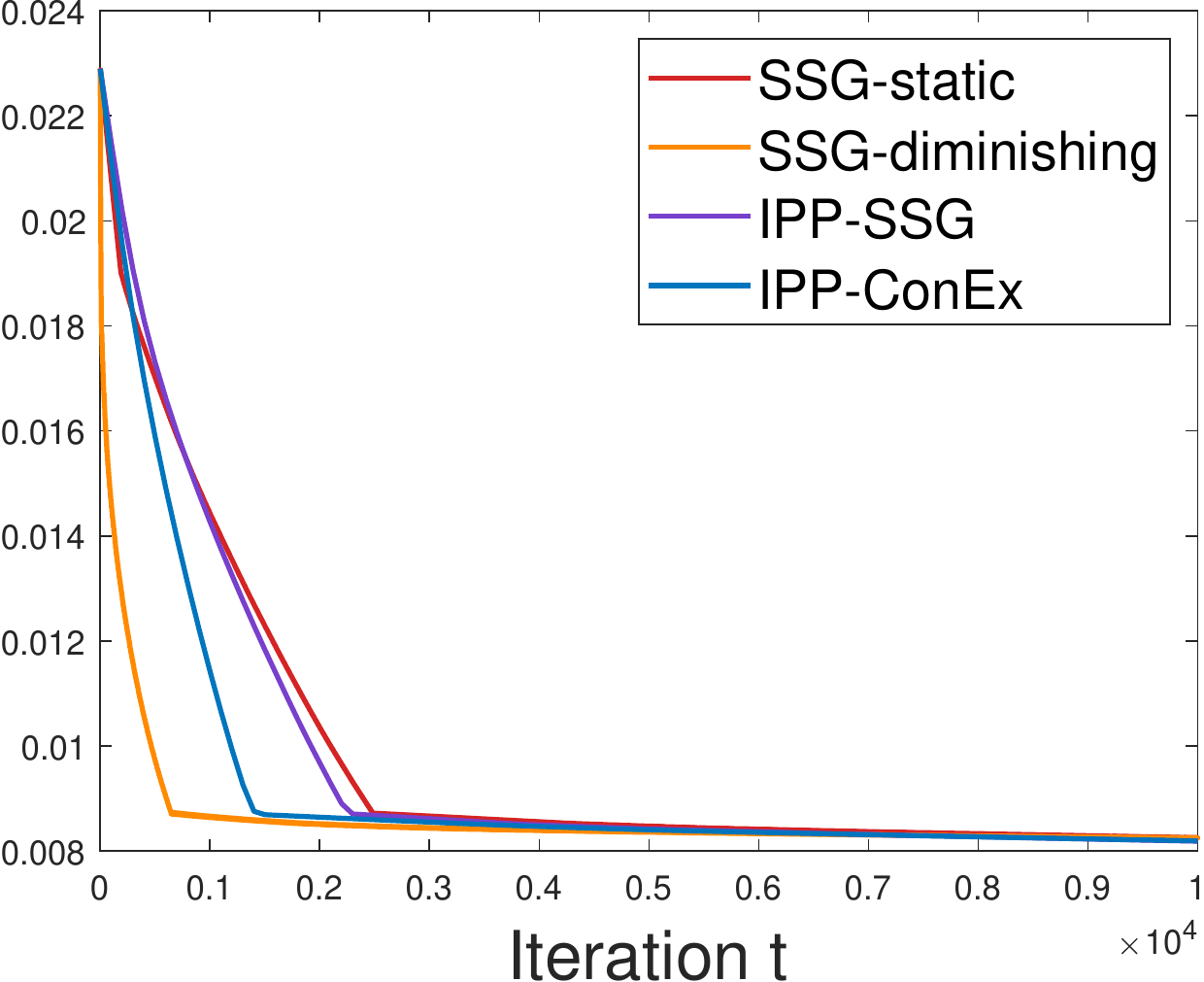}
          \\
		\raisebox{12ex}{\small{\rotatebox[origin=c]{90}{Infeasibility}}}
		& \hspace*{-0.06in}\includegraphics[width=0.30\textwidth]{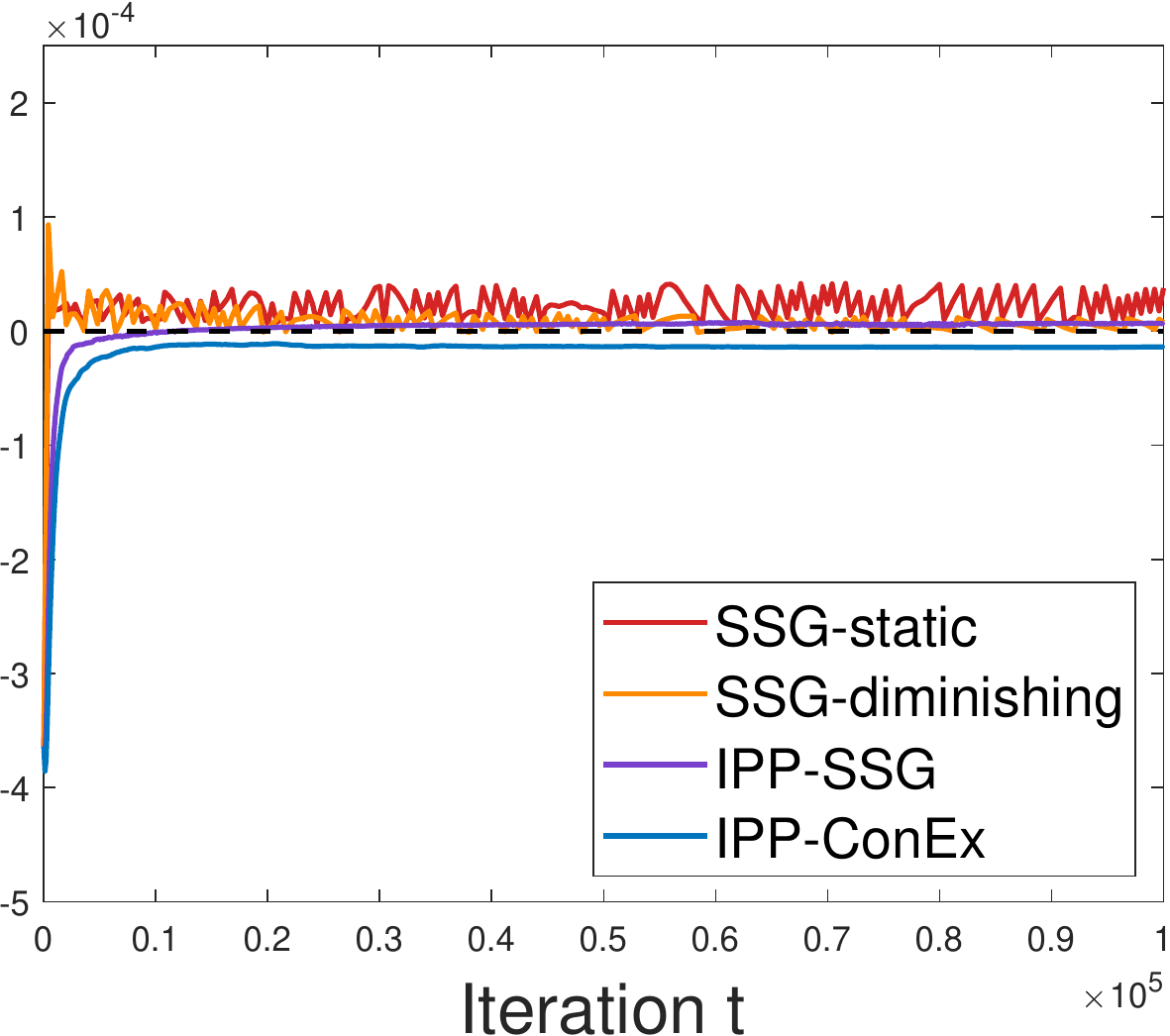}
		& \hspace*{-0.06in}\includegraphics[width=0.30\textwidth]{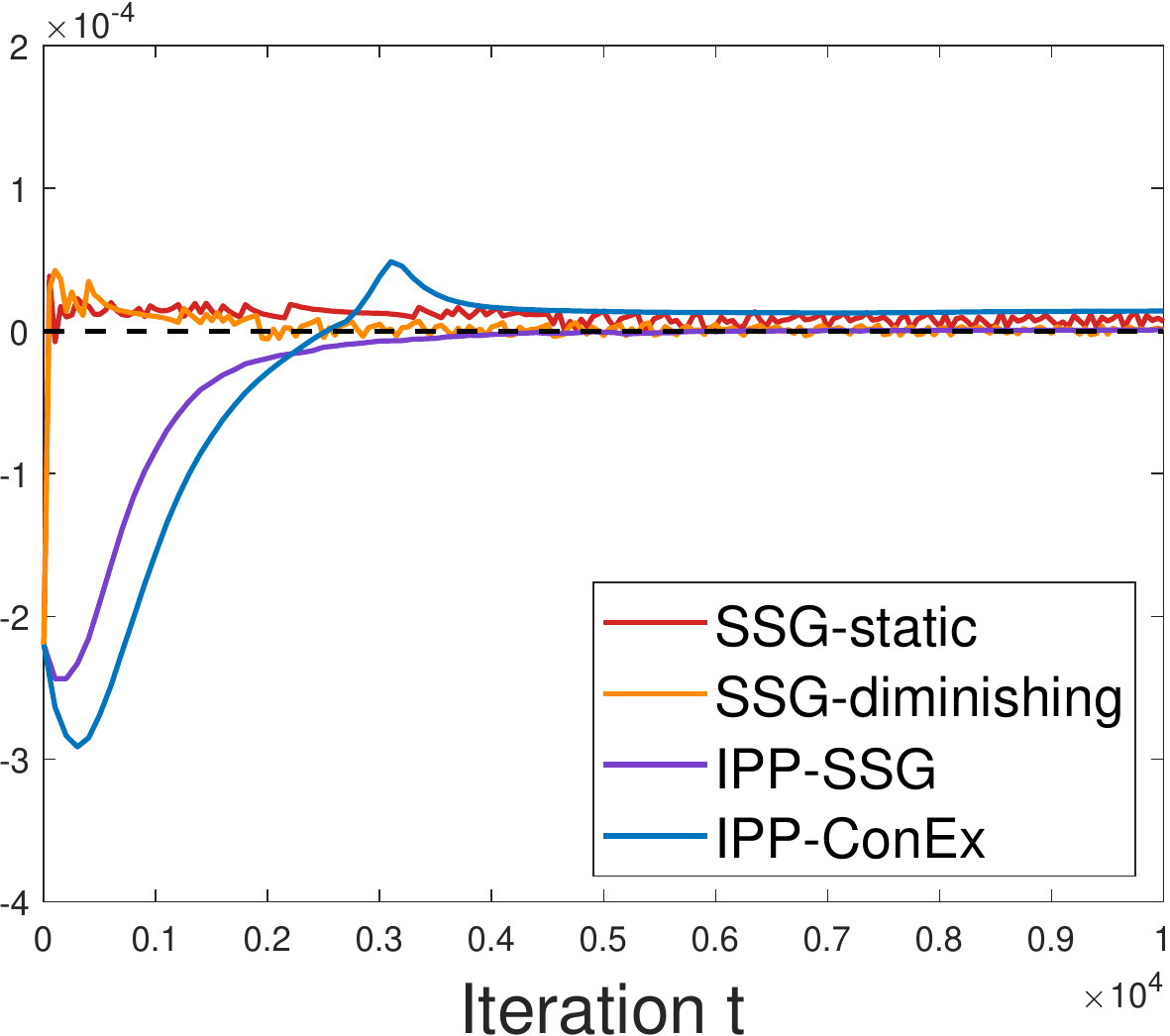}
		& \hspace*{-0.06in}\includegraphics[width=0.30\textwidth]{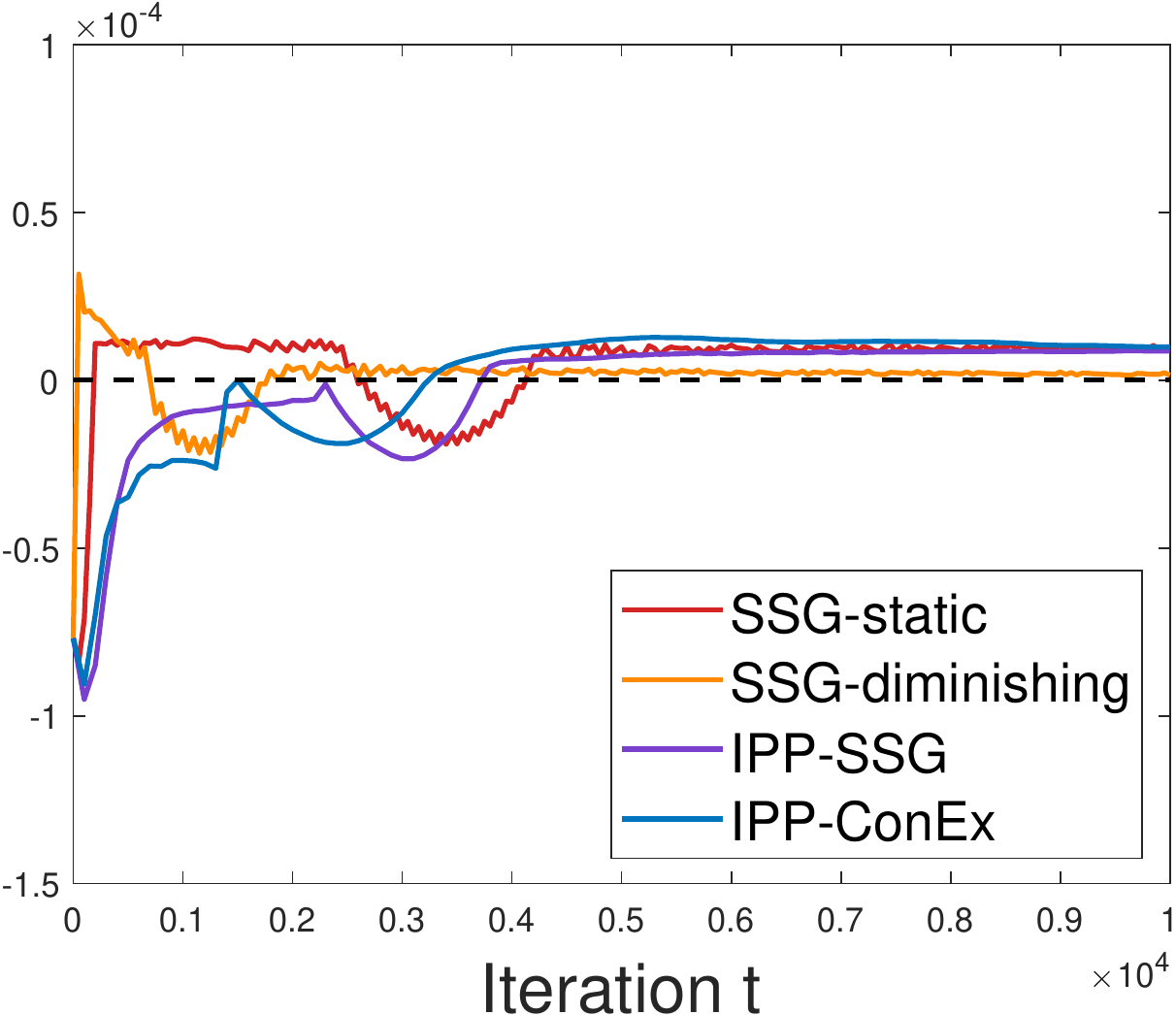}
           \\
        \raisebox{12ex}{\small{\rotatebox[origin=c]{90}{Near Stationarity}}}
		& \hspace*{-0.06in}\includegraphics[width=0.30\textwidth]{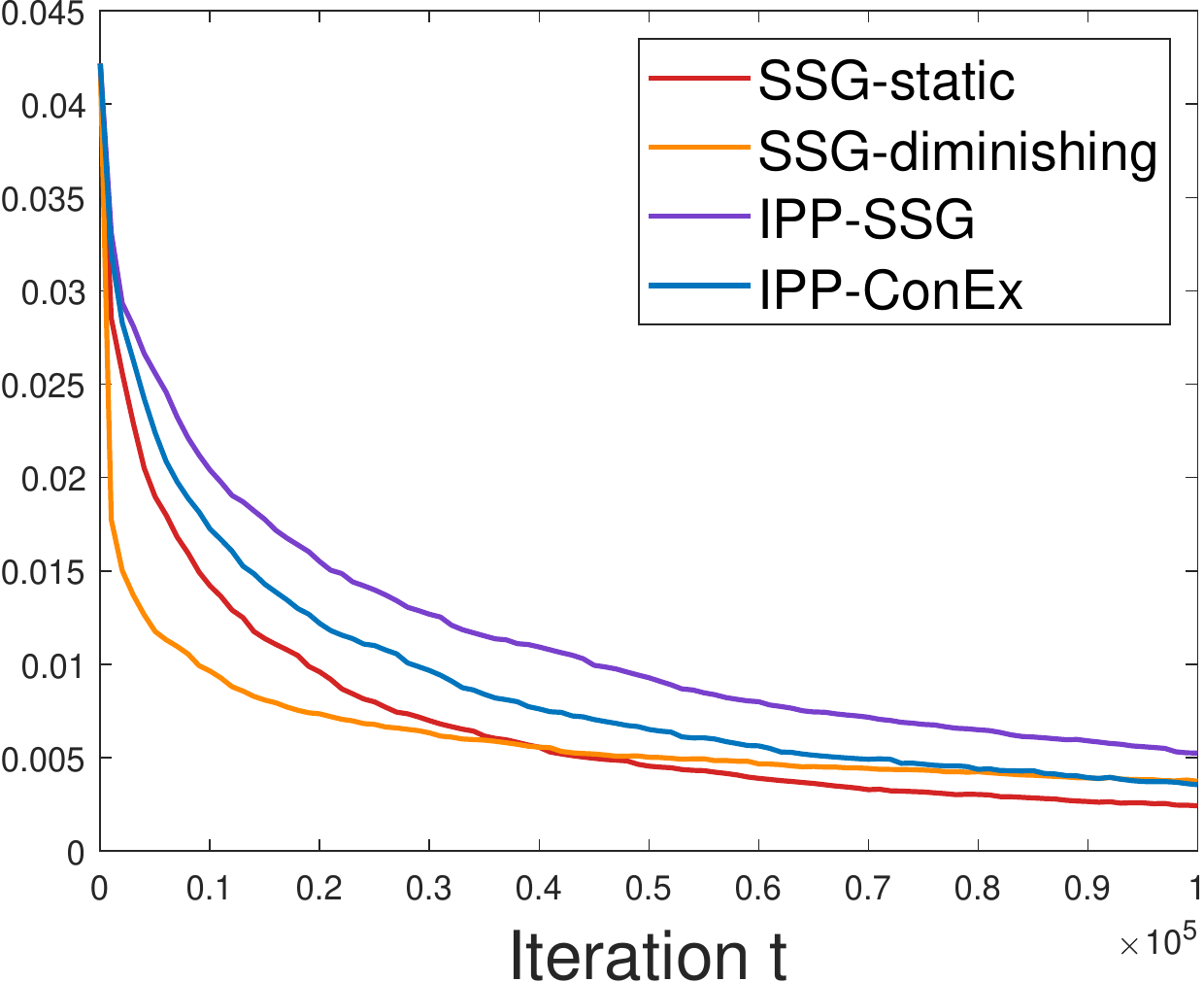}
		& \hspace*{-0.06in}\includegraphics[width=0.30\textwidth]{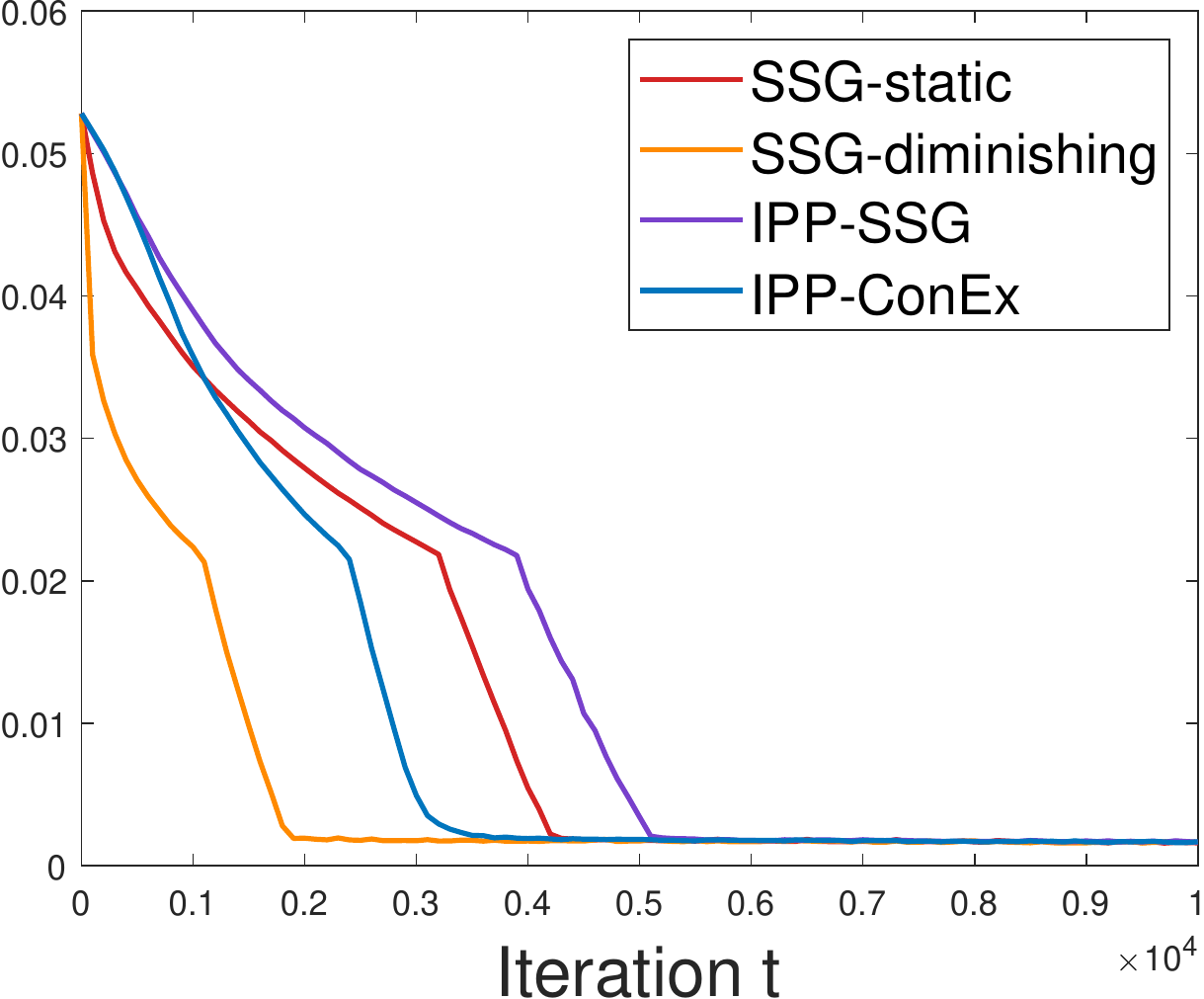}
		& \hspace*{-0.06in}\includegraphics[width=0.30\textwidth]{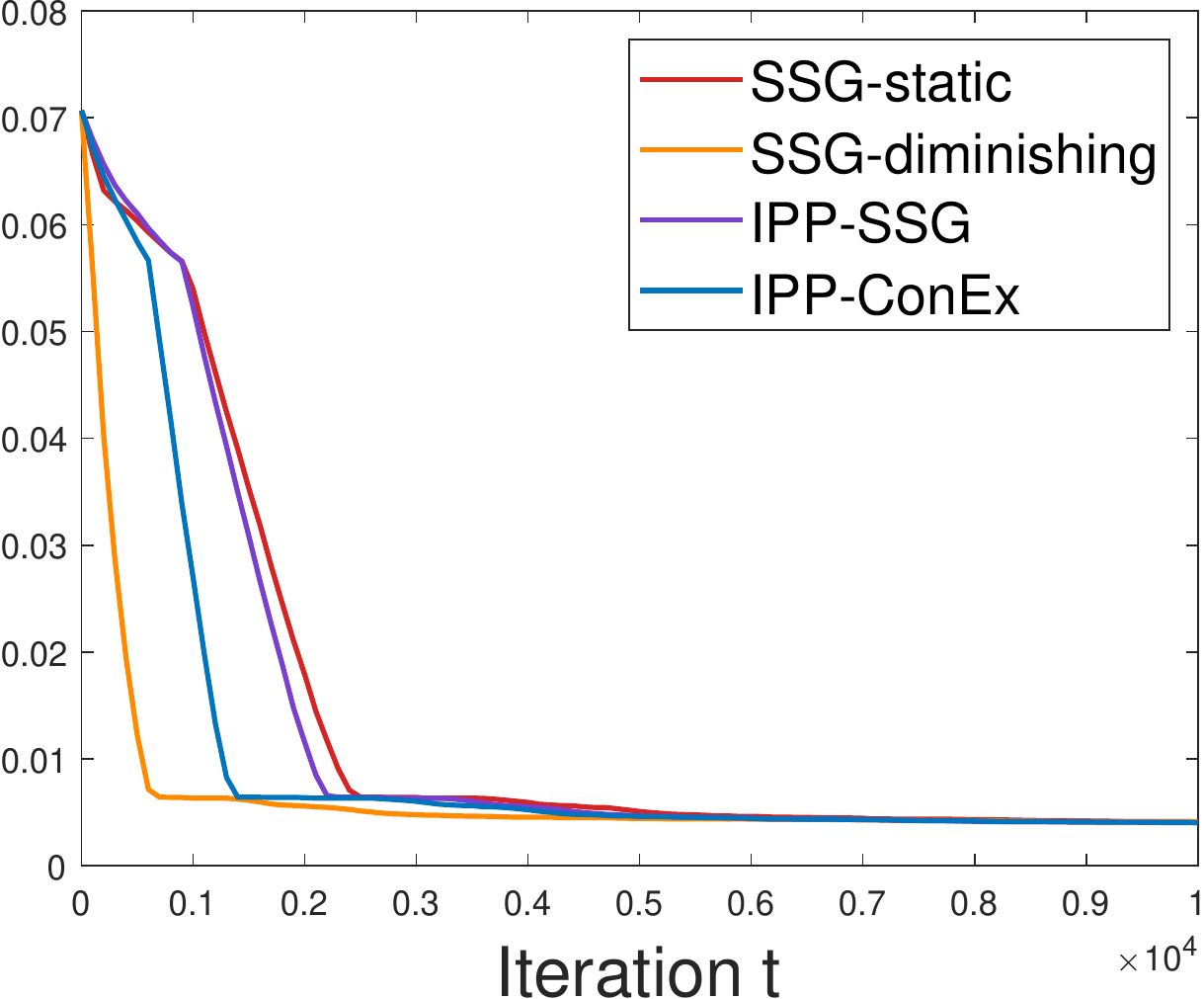}
        \end{tabular}
	\caption{Performances vs number of iterations on classification problems with ROC-based fairness.}  \label{fig:figure_convex_experiment_iteration}
	\vspace{-0.1in}
\end{figure*}

In our experiments, we choose $\ell(z)=(1-z)_+$ and first solve \eqref{eq:ermL} using the subgradient method with a large enough number of iterations to obtain $L^*$ and a solution $\bx_{\text{ERM}}$. Then we set $\kappa=0.001L^*$ and $r=5\|\bx_{\text{ERM}}\|$, and let $\Theta$ consist of $400$ points equally spaced between 
$
\min_{i}\bx_{\text{ERM}}^\top\ba_i-0.5(\max_{i}\bx_{\text{ERM}}^\top\ba_i-\min_{i}\bx_{\text{ERM}}^\top\ba_i)
$
and 
$
\max_{i}\bx_{\text{ERM}}^\top\ba_i+0.5(\max_{i}\bx_{\text{ERM}}^\top\ba_i-\min_{i}\bx_{\text{ERM}}^\top\ba_i).
$

All methods are initialized at $\bx_{\text{ERM}}$. We implemented the SSG method with both static and diminishing stepsizes. For the static stepsize, we select $\epsilon_t$ from $\{10^{-6},2\times10^{-6},5\times10^{-6},10^{-5}\}$ and $\eta_t$ from $\{2\times10^{-4},5\times10^{-4},10^{-3},2\times10^{-3}\}$. For the diminishing stepsize, we set $\epsilon_t=\frac{E_1}{\sqrt{t+1}}$ and $\eta_t=\frac{E_2}{\sqrt{t+1}}$ and select $E_1$ from $\{5\times10^{-5},10^{-4},2\times10^{-4},5\times10^{-4}\}$ and $E_2$ from $\{0.02,0.05,0.1,0.2\}$. We select the best set of parameters that produces the smallest objective value after 5000 iterations. For IPP, we select $\hat\rho$ from $\max\{\rho,1\}\times\{1,1.5,2\}$ for all three datasets, and the proximal point subproblem is approximately solved by SSG and ConEx both with $100$ iterations. For IPP-SSG, we apply a static stepsize with $\epsilon_t$ and $\eta_t$ tuned in the same way as SSG. For IPP-ConEx, following the notation in~\cite{boob2022stochastic}, we set $\theta_t=\frac{t}{t+1}$, $\eta_t=c_1(t+1)$ and $\tau_t=\frac{c_2}{t+1}$, and select $c_1$ from $\{20,50,100,200\}$ and $c_2$ from $\{0.002,0.005,0.01,0.02\}$ by the same procedure adopted by SSG.


We report the performances versus number of iterations of all methods on each dataset in Figure~\ref{fig:figure_convex_experiment_iteration}. The performances versus the used CPU time will be given in Section~\ref{sec:add_exp} in the Appendix. For SSG, the $x$-axis represents the total number of iterations while, for IPP, it represents the total number of inner iterations across all outer iterations. The $y$-axis in each row represents the objective value, infeasibility and near stationarity achieved at each iteration, respectively. To measure near stationarity, we solve \eqref{eq:phi} with $\bx=\bx^{(t)}$ and parameters \eqref{eq:parameter1} using the SSG method with $2500$ iterations and use the last iterate as an approximation of $\widehat\bx(\bx^{(t)})$. We make sure that the change of $\|\widehat\bx(\bx^{(t)})-\bx^{(t)}\|$ is less than 1\% if the number of iterations is increased to $5000$. Then we plot $\|\widehat\bx(\bx^{(t)})-\bx^{(t)}\|$ as near stationarity in Figure~\ref{fig:figure_convex_experiment_iteration}. Since computing $\widehat\bx(\bx^{(t)})$ with a high precision for each $t$ is time-consuming, we only report near stationarity at $100$ equally spaced iterations.

According to Figure~\ref{fig:figure_convex_experiment_iteration}, the SSG method with a diminishing stepsize has the best performance on all three datasets in the sense that it reduces the objective value and the (approximate) near stationarity measure faster than others while keeping the solutions nearly feasible. However, the SSG method with a static stepsize is not always better than 
the IPP methods. This is consistent with our theoretical finding that the SSG and IPP methods have the same oracle complexity. 

\subsection{Classification problem with demographic parity}
\label{sec:classification_demographic_parity}
Following the notation in the previous subsection, we consider a binary classification problem with a constraint enforcing demographic parity~\cite{agarwal2018reductions}. The measure of demographic parity and its continuous approximation are 
\small
\begin{equation}
\label{eq:demographicparity}
\Big|\frac{1}{n_p}\sum_{i=1}^{n_p} \mathbb{I}(\bx^\top\ba_i^p\geq0)-\frac{1}{n_u}\sum_{i=1}^{n_u} \mathbb{I}(\bx^\top\ba_i^u\geq0)
\Big|\approx R_0(\bx):=\Big|
\frac{1}{n_p}\sum_{i=1}^{n_p} \sigma(\bx^\top\ba_i^p)-\frac{1}{n_u}\sum_{i=1}^{n_u} \sigma(\bx^\top\ba_i^u)
\Big|.
\end{equation}
\normalsize
Fairness measure $R_0(\bx)$ is a special case of $R(\bx)$ with $\Theta=\{0\}$. If $R_0(\bx)$ is small, model $\bx$ produces similar predicted positive rates for the protected and unprotected groups. To obtain a fair $\bx$, we balance \eqref{eq:ermL} and \eqref{eq:demographicparity} by solving 
\begin{equation}
\label{eq:demographicparity_fairnessClassification}
\min L(\bx)+\lambda\text{SCAD}(\bx)\text{ s.t. } R_0(\bx)\leq \kappa.
\end{equation}
Different from \eqref{eq:ROCfairnessClassification}, the fairness measure is used as the constraint in \eqref{eq:demographicparity_fairnessClassification} while the objective function becomes the empirical hinge loss plus the smoothly clipped absolute deviation (SCAD) regularizer \cite{fan2001variable} for promoting a sparse solution. Here, $\lambda$ is a regularization parameter and 
\begin{align}
\label{eq:SCAD}
\text{SCAD}(\bx):=\sum_{i=1}^d s(x_i), \quad 
s(x_i)=\left\{
\begin{array}{ll}
2|x_i| & 0\leq |x_i|\leq 1\\
-x_i^2+4|x_i|+1 & 1<|x_i|\leq 2\\
3 & 2<|x_i|.
\end{array}
\right.
\end{align}
We prove in Section~\ref{sec:slaterexample} that \eqref{eq:demographicparity_fairnessClassification} is an instance of \eqref{eq:gco} satisfying Assumptions~\ref{assume:allpaper} and \ref{assume:weaklyconvex} with $\rho=\max\{2\lambda,\beta\}$ where $\beta$ is defined in \eqref{eq:hLip} in Section~\ref{sec:slaterexample}.

We set $\lambda=0.2$ for all datasets and set $\kappa=0.005$, $0.02$  and $0.02$ for \textit{a9a}, \textit{bank} and \textit{COMPAS}, respectively. For SSG, we select $\epsilon_t$ from $\{10^{-6},2\times10^{-6},5\times10^{-6},10^{-5}\}$ and select $\eta_t$ from $\{10^{-4},2\times10^{-4},5\times10^{-4},7.5\times10^{-4}\}$ for $t\in I$ while set $\eta_t=g(\bx^{(t)})/\|\bzt_g^{(t)}\|^2$ for $t\in J$. We select the best set of parameters that produces the smallest objective value after $50000$ iterations. For IPP, we select $\hat\rho$ from $\max\{\rho,1\}\times\{1,1.5,2\}$ and 
the proximal point subproblem is approximately solved by SSG and ConEx both with $600$ iterations. For IPP-SSG, we apply a static stepsize with $\epsilon_t$ and $\eta_t$ tuned in the same way as SSG. For IPP-ConEx, following the notation in~\cite{boob2022stochastic}, we set $\theta_t=\frac{t}{t+1}$, $\eta_t=c_1(t+1)$ and $\tau_t=\frac{c_2}{t+1}$ and select $c_1$ from $\{20,50,100,200\}$ and $c_2$ from $\{0.002,0.005,0.01,0.02\}$ by the same procedure adopted by SSG.  Due to the limit of space, we present the performances vs number of iterations and CPU runtime of all methods in Figure~\ref{fig:figure_weakly_convex_experiment_iteration} and Figure~\ref{fig:figure_weakly_convex_experiment_cputime} respectively in Section~\ref{sec:add_exp} in the Appendix. 

\section{Conclusion}
We study the oracle complexity of the switching subgradient (SSG) method for finding a nearly $\epsilon$-stationary point of a non-smooth weakly convex constrained optimization problem. We show that the complexity of the SSG method matches the best result in literature that is achieved only by double-loop methods. On the contrary, the SSG method is single-loop and easier to implement with reduced tuning effort. This is the first complexity result for a single-loop first-order method for a  weakly-convex non-smooth constrained problem.

\section*{Acknowledgements}
This work was supported by the NSF Grant No. 2147253.

\medskip

\bibliography{SIP,references,iclr2019_conference}
\bibliographystyle{plain}


\newpage
\appendix
\section{Complexity analysis for convex constraints}
In this section, we present the theoretical complexity analysis for the SSG method when $f$ is non-convex but $g$ is convex. Our analysis will first focus on a more general setting where the oracles of the subgradients and function values of $f$ and $g$ are stochastic. Then the complexity for the case of deterministic oracles is derived as a special case.

\subsection{Proof of Lemma~\ref{thm:boundlambda} and extension with linear equality constraints}
\label{sec:boundlambda}
\begin{proof}[Proof of Lemma~\ref{thm:boundlambda}]
For simplicity of notation, we denote $\widehat\bx(\bx)$ in \eqref{eq:phix} by $\widehat\bx$. By Assumption~\ref{assume:convex}B,  there exists a strictly feasible solution $\bx_{\text{feas}}\in\text{relint}(\X)$ with $g(\bx_{\text{feas}})<0$. As a result, the Lagrangian multiplier $\widehat\lambda\geq 0$  corresponding to $\widehat\bx$ is well-defined and satisfies \eqref{eq:KKTprox}, which means $\widehat\lambda g(\widehat{\bx}) = 0$ and 
 \begin{align}
\label{eq:optcond1}
\widehat\bzt_f  +\hat{\rho}(\widehat\bx - \bx)+\widehat\lambda \widehat\bzt_g+\widehat\bu=\mathbf{0},
\end{align}
where $\widehat\bzt_f\in\partial f(\widehat\bx)$, $\widehat\bzt_g\in\partial g(\widehat\bx)$, $\widehat\bu\in \mathcal{N}_\X(\widehat\bx)$ and $\mathcal{N}_\X(\widehat\bx)$ is the normal cone of $\X$ at $\widehat\bx$. 
 
If $\widehat\lambda = 0$, the conclusion holds trivially. Hence, we focus on the case that $\widehat\lambda > 0$. Note that, in this case, we must have $g(\widehat{\bx})= 0$. Taking the inner product between \eqref{eq:optcond1} and $\widehat\bx -\bx_{\text{feas}}$ gives
  \begin{eqnarray}
  \label{eq:upperboundlambdaeq1}
  \frac{1}{\widehat\lambda}\langle\widehat\bzt_f  +\hat{\rho}(\widehat\bx - \bx), \widehat\bx -\bx_{\text{feas}}\rangle
  =\langle\widehat\bzt_g+\widehat\bu/\widehat\lambda, \bx_{\text{feas}}-\widehat\bx \rangle
  \leq g(\bx_{\text{feas}})- g(\widehat\bx)=g(\bx_{\text{feas}}),
  \end{eqnarray}
  where the inequality is because of convexity of $g$ and the fact that $\widehat\bu/\widehat\lambda\in \mathcal{N}_\X(\widehat\bx)$. By Assumption~\ref{assume:allpaper} and \ref{assume:convex}C, we have $\|\widehat\bzt_f\|\leq M$, $\|\widehat\bx -\bx_{\text{feas}}\|\leq D$ and $\|\widehat\bx -\bx\|\leq D$, which imply from \eqref{eq:upperboundlambdaeq1} that 
		\small
		\begin{align*}
		\widehat\lambda \leq \frac{\langle\widehat\bzt_f  +\hat{\rho}(\widehat\bx - \bx), \bx_{\text{feas}}-\widehat\bx \rangle}{-g(\bx_{\text{feas}})} 
		\leq  \frac{MD+\hat\rho D^2 }{-g(\bx_{\text{feas}})} .
		\end{align*}
		\normalsize
\end{proof}

\textbf{The extension of Lemma~\ref{thm:boundlambda} when there exist linear constraints: }
Suppose there are linear equality constraints. We formulate problem \eqref{eq:gco} as 
\begin{eqnarray}
\label{eq:gco_le}
f^*\equiv \min_{\bx\in\mathcal{X}}f(\bx) \quad \text{s.t.} \quad h(\bx)\leq 0,\quad \bA\bx=\bb,
\end{eqnarray}
where $\bA\in\mathbb{R}^{l\times d}$ and $\bb\in\mathbb{R}^l$. Moreover, we assume $\bA$ has a full row rank. We then define 
\begin{eqnarray}
\label{eq:newg_le}
g(\bx)=\max\{h(\bx), \|\bA\bx-\bb\|_\infty\}.
\end{eqnarray}
With this $g$, Assumption~\ref{assume:convex} B cannot be satisfied. However, it is easy to show that $\widehat\lambda$ in \eqref{eq:KKTprox} can be still bounded from above if the inequality constraint in \eqref{eq:gco_le} can be satisfied strictly.

To extend Lemma \ref{thm:boundlambda} for \eqref{eq:gco_le}, we only need to change Assumption~\ref{assume:convex} B to the following
\begin{itemize}
\item[B'.] (Slater's condition) There exists $\bx_{\text{feas}}\in\text{int}(\X)$ such that $h(\bx_{\text{feas}})<0$ and  $\bA\bx_{\text{feas}}=\bb$.
\end{itemize}
Then we can still derive an upper bound of $\widehat\lambda$ similar to the one in Lemma \ref{thm:boundlambda}. 
\begin{lemma}
\label{thm:boundlambda_le}
Let $g$ defined as in \eqref{eq:newg_le}. Suppose Assumptions~\ref{assume:allpaper} and \ref{assume:convex} hold except that Assumption~\ref{assume:convex} B is replaced by B' above. Given any $\bx\in\X$, let $\widehat\bx(\bx)$ be defined as in \eqref{eq:phix} with $(\hat\rho,\tilde\rho)$ satisfying \eqref{eq:parameter1} and $\widehat\lambda$ be the associated Lagrangian multiplier satisfying \eqref{eq:KKTprox}. We have
\small
\begin{eqnarray}
\label{eq:Lambdabound_le}
\widehat\lambda \leq  \Lambda_{ec}:= \frac{MD+\hat\rho D^2 }{-h(\bx_{\textup{feas}})}
+\sqrt{l}\|(\bA\bA^\top)^{-1}\bA\|\left(\frac{M^2D+\hat\rho MD^2}{-h(\bx_{\textup{feas}})}+\frac{MD+\hat\rho D^2}{\mathrm{dist}(\bx_{\textup{feas}},\partial \mathcal{X})}+ M+\hat\rho D\right).
\end{eqnarray}
\normalsize
\end{lemma}
\begin{proof}
For simplicity of notation, we  denote $\widehat\bx(\bx)$ in \eqref{eq:phix} by $\widehat\bx$. We first show that \eqref{eq:KKTprox} holds for some Lagrangian multiplier $\widehat\lambda$ and $\widehat\bzt_g\in \partial g(\widehat\bx)$.
By Assumption B', there exist Lagrangian multipliers $\widehat\lambda_{in} \in\mathbb{R}_+$ and $\widehat\blambda_{eq}\in\mathbb{R}^{l}$ such that
 \begin{align}
\label{eq:optcond1_le}
\widehat\bzt_f  +\hat{\rho}(\widehat\bx - \bx)+\widehat\lambda_{in} \widehat\bzt_h+\bA^\top\widehat\blambda_{eq}+\widehat\bu=\mathbf{0},
\end{align}
where $\widehat\bzt_h\in\partial h(\widehat\bx)$,  and other notations are as in \eqref{eq:optcond1}. It also holds that $g(\widehat\bx)=0$ and $\widehat\lambda_{in}h(\widehat\bx)=0$. 

Since $h(\widehat\bx)\leq0$ and $\bA\widehat\bx=\bb$, we have $\mathbf{0} \in \partial g(\widehat\bx)$ by the definition of $g$. Suppose $\widehat\lambda_{in}+\|\widehat\blambda_{eq}\|_1=0$. \eqref{eq:optcond1_le} implies \eqref{eq:optcond1} with $\widehat\lambda=0$ and any $\widehat\bzt_g\in \partial g(\widehat\bx)$. Suppose $\widehat\lambda_{in}+\|\widehat\blambda_{eq}\|_1>0$. Let
$$
\widehat\bzt_g=\frac{\widehat\lambda_{in} \widehat\bzt_h+\bA^\top\widehat\blambda_{eq}}{\widehat\lambda_{in}+\|\widehat\blambda_{eq}\|_1}
\quad\text{ and }\quad\widehat\lambda=\widehat\lambda_{in}+\|\widehat\blambda_{eq}\|_1.
$$
It is easy to verify by definition that $\widehat\bzt_g\in \partial g(\widehat\bx)$ and \eqref{eq:optcond1}  also holds.

Next we derive the claimed upper bound of $\widehat\lambda$. When $\widehat\lambda_{in}+\|\widehat\blambda_{eq}\|_1=0$ and thus $\widehat\lambda = 0$, the conclusion holds trivially so we focus on the case where  $\widehat\lambda =\widehat\lambda_{in}+\|\widehat\blambda_{eq}\|_1> 0$ and $\widehat\bzt_g$ is defined as above. In this case, we follow the same analysis in the proof of Lemma 3 in \cite{lin2022complexity} to bound $\widehat\lambda_{in}$ and $\|\widehat\blambda_{eq}\|_1$.

By the convexity of $h$ and the fact that $\widehat\lambda_{in}\geq 0$, we can show that
\begin{align}
\nonumber
\widehat\lambda_{in}h(\bx_{\text{feas}}) 
\geq&~\widehat\lambda_{in}\left(h(\widehat\bx) +\langle\widehat\bzt_h,\bx_{\text{feas}}-\widehat\bx\rangle\right)\\\nonumber
=&-\langle\widehat\bzt_f  +\hat{\rho}(\widehat\bx - \bx)+\bA^\top\widehat\blambda_{eq}+\widehat\bu,\bx_{\text{feas}}-\widehat\bx\rangle\\\label{eq:boundlmabda1_le}
=&-\langle\widehat\bzt_f  +\hat{\rho}(\widehat\bx - \bx)+\widehat\bu,\bx_{\text{feas}}-\widehat\bx\rangle,
\end{align}
where the first equality is due to \eqref{eq:optcond1_le} and the second equality is because $\bA\widehat\bx=\bA\bx_{\text{feas}}=\bb$. 

Next we bound the term $-\langle\widehat\bu,\bx_{\text{feas}}-\widehat\bx\rangle$ in \eqref{eq:boundlmabda1_le} from below. Suppose $\widehat\bu\neq\mathbf{0}$. It must hold that $\widehat\bx\in\partial\X$ since, otherwise,  $\widehat\bx\in\text{int}(\X)$ so $\widehat\bu=\mathbf{0}$. Let $\mathcal{H}=\big\{\bx\in\mathbb{R}^d\,|\,\langle\widehat\bu,\bx-\widehat\bx\rangle=0\big\}$, which is a supporting hyperplane of $\X$ at $\widehat\bx$. Since $\mathrm{dist}(\bx_{\text{feas}},\mathcal{H}) = \langle\widehat\bu,\widehat\bx-\bx_{\text{feas}}\rangle/\|\widehat\bu\|$, we have 
$$
-\langle\widehat\bu,\bx_{\text{feas}}-\widehat\bx\rangle=\mathrm{dist}(\bx_{\text{feas}},\mathcal{H})\|\widehat\bu\| \ge \mathrm{dist}(\bx_{\text{feas}},\partial \mathcal{X})\|\widehat\bu\|.
$$		

Applying the inequality to \eqref{eq:boundlmabda1_le}, we have
$$
-\widehat\lambda_{in}h(\bx_{\text{feas}}) + \mathrm{dist}(\bx_{\text{feas}},\partial \mathcal{X})\|\widehat\bu\|  \leq\langle\widehat\bzt_f+\hat{\rho}(\widehat\bx - \bx),\bx_{\text{feas}}-\widehat\bx\rangle\leq MD+\hat\rho D^2,
$$
where the last inequality is by Assumption~\ref{assume:allpaper} and \ref{assume:convex}C, which imply $\|\widehat\bzt_f\|\leq M$, $\|\widehat\bx -\bx_{\text{feas}}\|\leq D$ and $\|\widehat\bx -\bx\|\leq D$.
The  inequality above then implies
\begin{eqnarray}
	\label{eq:boundlambdak_le}
\widehat\lambda_{in}\leq \frac{MD+\hat\rho D^2 }{-h(\bx_{\text{feas}})}\text{ and }\|\widehat\bu\| \leq \frac{MD+\hat\rho D^2}{\mathrm{dist}(\bx_{\text{feas}},\partial \X)}.
\end{eqnarray}
Suppose $\widehat\bu=\mathbf{0}$. We can also prove the first inequality in \eqref{eq:boundlambdak_le} by a similar argument and the second inequality in \eqref{eq:boundlambdak_le} holds trivially.   Multiplying both sides of \eqref{eq:optcond1_le} by $\bA$ and using the fact that $\bA$ has a full row rank, we have 
$$
\widehat\blambda_{eq}= -(\bA\bA^\top)^{-1}\bA\left(\widehat\bzt_f  +\hat{\rho}(\widehat\bx - \bx)+\widehat\lambda_{in} \widehat\bzt_h+\widehat\bu\right),
$$
which, together with \eqref{eq:boundlambdak_le}, further implies
\begin{align}
	\nonumber
	\|\widehat\blambda_{eq}\|_1\leq\sqrt{l}\|\widehat\blambda_{eq}\|
	\leq&~\sqrt{l}\|(\bA\bA^\top)^{-1}\bA\|\left(M(1+\widehat\lambda_{in})+\hat{\rho}D+\|\widehat\bu\|\right)\\\nonumber
	\leq&~\sqrt{l}\|(\bA\bA^\top)^{-1}\bA\|\left(\frac{M^2D+\hat\rho MD^2}{-h(\bx_{\text{feas}})}+\frac{MD+\hat\rho D^2}{\mathrm{dist}(\bx_{\text{feas}},\partial \mathcal{X})}+ M+\hat\rho D\right).
\end{align}
Applying this upper bound of $\|\widehat\blambda_{eq}\|$ and the upper bound of $\widehat\lambda_{in}$ in \eqref{eq:boundlambdak_le} to $\widehat\lambda =\widehat\lambda_{in}+\|\widehat\blambda_{eq}\|_1$, we obtain the desired upper bound of $\widehat\lambda$. 
\end{proof}

\begin{remark}
\label{remark:complexity_ec}
When there are linear equality constraints, we can just implement Algorithm~\ref{alg:dsgm} with the new definition of $g$ in \eqref{eq:newg_le}. Using the upper bound of $\widehat\lambda_t$ in \eqref{eq:Lambdabound_le} in place of the upper bound \eqref{eq:Lambdabound} in the proof of Theorem~\ref{thm:mdconverge_deterministic}, we can show that Theorem~\ref{thm:mdconverge_deterministic} still holds except that $\Lambda$ is repalced by $\Lambda_{ec}$.
\end{remark}

\subsection{Stochastic switching subgradient method}
\label{sec:stochasticassumption}
In this section, we introduce a more general setting where the subgradients and function values of $f$ and $g$ can only be accessed through stochastic oracles, and then propose a stochastic switching subgradient method using these oracles. More specifically, in addition to Assumptions~\ref{assume:allpaper} and \ref{assume:convex}, we make the following assumption. 
\begin{assumption}
\label{assume:stochasticgrad_prob}
For any $\bx\in\X$, stochastic subgradients $\bxi_f$ and $\bxi_g$ and a stochastic value $\omega$ can be generated independently such that $\mathbb{E}[\omega]=g (\bx)$, $\mathbb{E}[\bxi_f]\in\partial f (\bx)$ and $\mathbb{E}[\bxi_g]\in\partial g (\bx)$. Moreover, it holds that 
$
\mathbb{E}\left[\exp\left(\|\bxi_f\|^2/M^2\right)\right]\leq \exp(1)$, $\mathbb{E}\left[\exp\left(\|\bxi_g\|^2/M^2\right)\right]\leq \exp(1)
$, and
$
\mathbb{E}\left[\exp\left((\omega-g(\bx))^2/\sigma^2\right)\right]\leq \exp(1)
$
for a constant $\sigma$ for any $\bx\in\X$.
\end{assumption}
Assumption~\ref{assume:stochasticgrad_prob} means the distributions of $\|\bxi_f\|^2$, $\|\bxi_g\|^2$ and $(\omega-g(\bx))^2$ have light tails, which is a common assumption for proving a stochastic first-order method converges in a high probability. See (2.50) in \cite{nemirovski2009robust} for an example. We follow the convention that $\sigma=0$ when $w=g(\bx)$ deterministically. 

Under Assumption~\ref{assume:stochasticgrad_prob}, we have to use stochastic subgradients $\bxi_f$ and $\bxi_g$ in the SSG method. Moreover, since $g(\bx^{(t)})$ cannot be evaluated exactly, we sample a mini-batch $\{\omega_i^{(t)}\}_{i=1}^B$ of stochastic estimators of $g(\bx^{(t)})$ with a batch size of $B$, and then construct an estimator of $g(\bx^{(t)})$, denoted by $\bar w^{(t)}$, by averaging the samples. Then, we use the condition $\bar w^{(t)}\leq\epsilon_t$ in the SSG method to determine when the stochastic subgradient should be switched between $\bxi_f$ and $\bxi_g$. The SSG method with those modifications is called the stochastic SSG method and presented in Algorithm~\ref{alg:sgm}. Note that Algorithm~\ref{alg:sgm}  is reduced to Algorithm~\ref{alg:dsgm} when, for $t\geq0$, it holds deterministically that
\begin{eqnarray}
\label{eq:reduced_to_deterministic}
\bxi_f^{(t)}=\bzt_f^{(t)}, ~\bxi_g^{(t)}=\bzt_g^{(t)}, ~\bar\omega^{(t)}=g(\bx^{(t)}), ~\sigma=0~\text{ and }~B=1.
\end{eqnarray}
The oracle complexity of Algorithm~\ref{alg:dsgm} is given in Theorem~\ref{thm:mdconverge_deterministic} while that of  Algorithm~\ref{alg:sgm} is given in Theorems~\ref{thm:mdconverge} and \ref{thm:mdconverge_prob}.




\subsection{Technical lemmas and propositions}
\label{sec:mainproposition}



\begin{algorithm}[tb]
   \caption{Stochastic Switching Subgradient Method}
   \label{alg:sgm}
\begin{algorithmic}[1]
   \STATE {\bfseries Input:} $\bx^{(0)}\in \X$, total number of iterations $T$, tolerance of infeasibility $\epsilon_t\geq0$, stepsize $\eta_t$, mini-batch size $B$, and a starting index $S$ for generating outputs.
   \STATE {\bfseries Initialization:} $I=\emptyset$ and $J=\emptyset$.
   \FOR{iteration $t=0,1,\cdots,T-1$}
   \STATE Generate a mini-batch $\{\omega_i^{(t)}\}_{i=1}^B$ of stochastic estimators of $g(\bx^{(t)})$.
   \STATE Set $\bar\omega^{(t)}=\frac{1}{B}\sum_{i=1}^B\omega_i^{(t)}$.
   \IF{$\bar\omega^{(t)}\leq\epsilon_t$}
    \STATE Generate a stochastic subgradient $\bxi_f^{(t)}$ of $f$  at $\bx^{(t)}$.\\
    \STATE Set $\bx^{(t+1)}= \text{proj}_{\X} (\bx^{(t)}-\eta_{t}\bxi_f^{(t)} )$
    and,  if $t\geq S$, $I= I\cup\{t\}$.
    \ELSE
    \STATE Generate a stochastic subgradient $\bxi_g^{(t)}$ of $g$  at $\bx^{(t)}$.\\
    \STATE Set $\bx^{(t+1)}= \text{proj}_{\X} (\bx^{(t)}-\eta_{t}\bxi_g^{(t)} )$
     and,  if $t\geq S$, $I= I\cup\{t\}$. 
   \ENDIF
   \ENDFOR
\STATE {\bfseries Output:} $\bx^{(\tau)}$ with $\tau$ randomly sampled from  $I$ using $\text{Prob}(\tau=t)=\eta_t/\sum_{s\in I}\eta_s$.
\end{algorithmic}
\end{algorithm}

In this section, we first introduce additional notations and then present a few technical lemmas
and propositions which are needed to prove the convergence properties of Algorithms~\ref{alg:dsgm} and \ref{alg:sgm}. 

Let $\mathbb{I}(\cdot)$ be an indicator of an event, which equals one if the event occurs and zero otherwise. For each iterate $\bx^{(t)}$ in Algorithms~\ref{alg:dsgm} and ~\ref{alg:sgm}, let $\bzt_f^{(t)}:=\mathbb{E}[\bxi_f^{(t)}]\in\partial f(\bx^{(t)})$ and $\bzt_g^{(t)}:=\mathbb{E}[\bxi_g^{(t)}]\in\partial g(\bx^{(t)})$ and let $\widehat\bx^{(t)}:=\widehat\bx(\bx^{(t)})$ be defined as in \eqref{eq:phix} with $(\hat\rho,\tilde\rho)$ satisfying \eqref{eq:parameter1} and let $\widehat\lambda_t\geq 0$ be the corresponding Lagrangian multiplier satisfying \eqref{eq:KKTprox}, which exists by Assumption~\ref{assume:convex}B. Under Assumption~\ref{assume:stochasticgrad_prob}, let  $\mathbb{E}_t[\cdot]:=\mathbb{E}[\cdot~\big| ~\bar\omega^{(0)},\bxi_f^{(0)},\bxi_g^{(0)},\bar\omega^{(1)},\bxi_f^{(1)},\bxi_g^{(1)},\dots,\bar\omega^{(t-1)},\bxi_f^{(t-1)},\bxi_g^{(t-1)}]$, i.e., the conditional expectation conditioning on all stochastic events  before iteration $t$.
Let $\mathbb{E}_\tau[\cdot]$ be the expectation taken only over the random index $\tau$ when the algorithm is terminated.  
The proposition below characterizes the relationship between $\|\widehat\bx^{(t)}-\bx^{(t)}\|$ and the  parameters $T$, $S$, $\eta_t$ and $\epsilon_t$.


\begin{proposition}
\label{thm:mainprop}
Suppose Assumptions~\ref{assume:allpaper}, \ref{assume:convex} and \ref{assume:stochasticgrad_prob} hold and $g$ is $\mu$-strongly convex ($\mu$ can be zero).  Algorithm~\ref{alg:sgm} guarantees
\small
\begin{align*}
&~\sum_{t=S}^{T-1}\left[\eta_{t}\hat{\rho}\widehat\lambda_t\mathbb{I}(\bar\omega^{(t)}\leq \epsilon_t)+\eta_{t}\hat{\rho}\mathbb{I}(\bar\omega^{(t)}> \epsilon_t)\right]\cdot\frac{\mu}{2}\|\widehat\bx^{(t)}-\bx^{(t)} \|^{2}\\
&+\sum_{t=S}^{T-1}\left[\eta_{t}\hat{\rho}(\hat{\rho}-\rho)\|\widehat\bx^{(t)}-\bx^{(t)} \|^{2}-\eta_{t}\hat{\rho}\widehat\lambda_{t}g (\bx^{(t)} )\right]\mathbb{I}(\bar\omega^{(t)}\leq \epsilon_t)+\sum_{t=S}^{T-1}\eta_{t}\hat{\rho}g (\bx^{(t)} )\mathbb{I}(\bar\omega^{(t)}> \epsilon_t)\\
\leq&~\frac{\hat\rho D^2}{2}
+\dfrac{\hat{\rho}}{2}\sum_{t=S}^{T-1}\left[\eta_{t}^{2}\|\bxi_f^{(t)}\|^{2}\mathbb{I}(\bar\omega^{(t)}\leq \epsilon_t)+\eta_{t}^{2}\|\bxi_g^{(t)}\|^{2}\mathbb{I}(\bar\omega^{(t)}> \epsilon_t)\right]\\
&+\sum_{t=S}^{T-1}\left[\eta_{t}\hat{\rho}\left\langle\bxi_f^{(t)}-\bzt_f^{(t)},\widehat\bx^{(t)}-\bx^{(t)}\right\rangle\mathbb{I}(\bar\omega^{(t)}\leq \epsilon_t)+\eta_{t}\hat{\rho}\left\langle\bxi_g^{(t)}-\bzt_g^{(t)},\widehat\bx^{(t)}-\bx^{(t)}\right\rangle\mathbb{I}(\bar\omega^{(t)}> \epsilon_t)\right].
\end{align*}
\normalsize
\end{proposition}
\begin{proof}
Let $\bxi^{(t)}=\bxi_f^{(t)}$ if $t\in I$ and  $\bxi^{(t)}=\bxi_g^{(t)}$ if $t\in J$. Similarly, let $\bzt^{(t)}=\bzt_f^{(t)}\in\partial f(\bx^{(t)})$ if $t\in I$ and  $\bzt^{(t)}=\bzt_g^{(t)}\in\partial g(\bx^{(t)})$ if $t\in J$. By the updating equation of  $\bx^{(t+1)}$, we have
\begin{align*}
     &~\|\bx^{(t+1)}-\widehat\bx^{(t)} \|^{2}= \|\text{proj}_{\X} (\bx^{(t)}-\eta_{t}\bxi^{(t)} )-\widehat\bx^{(t)} \|^{2}= \|\text{proj}_{\X} (\bx^{(t)}-\eta_{t}\bxi^{(t)} )-\text{proj}_{\X} (\widehat\bx^{(t)})\|^{2} \\
    \leq&~ \|\bx^{(t)}-\eta_{t}\bxi^{(t)}-\widehat\bx^{(t)}  \|^{2}=\|\bx^{(t)}-\widehat\bx^{(t)} \|^{2}-2\eta_{t}\left\langle\bxi^{(t)},\bx^{(t)}-\widehat\bx^{(t)}\right\rangle+\eta_{t}^{2} \|\bxi^{(t)}\|^{2}.
\end{align*}
Multiplying the inequality above by $\hat{\rho}/2$ and adding $f (\widehat\bx^{(t)})$ to both sides, we obtain
\begin{align}
\nonumber
  &~f (\widehat\bx^{(t)} )+\dfrac{\hat{\rho}}{2} \|\bx^{(t+1)}-\widehat\bx^{(t)} \|^{2}\\\nonumber
  \leq&~f (\widehat\bx^{(t)} )+\dfrac{\hat{\rho}}{2} \|\bx^{(t)}-\widehat\bx^{(t)} \|^{2}-\eta_{t}\hat{\rho}\left\langle\bxi^{(t)},\bx^{(t)}-\widehat\bx^{(t)}\right\rangle+\frac{\eta_{t}^{2}\hat{\rho}}{2} \|\bxi^{(t)}\|^{2}\\\label{eq:mainthmeq1}
    =&~\varphi (\bx^{(t)})-\eta_{t}\hat{\rho}\left\langle\bzt^{(t)},\bx^{(t)}-\widehat\bx^{(t)}\right\rangle+\dfrac{\eta_{t}^{2}\hat{\rho}}{2}\|\bxi^{(t)}\|^{2}-\eta_{t}\hat{\rho}\left\langle\bxi^{(t)}-\bzt^{(t)},\bx^{(t)}-\widehat\bx^{(t)}\right\rangle,
\end{align}
where the equality is by the definition of $\varphi (\bx)$ in \eqref{eq:phi}. Since $\widehat\bx^{(t)}$ is  a feasible solution to problem \eqref{eq:phi} with $\bx=\bx^{(t+1)}$, we have 
\begin{align*}
    \varphi (\bx^{(t+1)})\leq f (\widehat\bx^{(t)} )+\dfrac{\hat{\rho}}{2} \|\bx^{(t+1)}-\widehat\bx^{(t)} \|^{2},
\end{align*}
which, together with \eqref{eq:mainthmeq1}, implies 
\begin{eqnarray}
\label{eq:mainthmeq2}
  \eta_{t}\hat{\rho}\left\langle\bzt^{(t)},\bx^{(t)}-\widehat\bx^{(t)}\right\rangle
    \leq\varphi (\bx^{(t)})-\varphi (\bx^{(t+1)})+\dfrac{\eta_{t}^{2}\hat{\rho}}{2}\|\bxi^{(t)}\|^{2}
    -\eta_{t}\hat{\rho}\left\langle\bxi^{(t)}-\bzt^{(t)},\bx^{(t)}-\widehat\bx^{(t)}\right\rangle.
\end{eqnarray}
Next, we will bound $\left\langle\bzt^{(t)},\bx^{(t)}-\widehat\bx^{(t)}\right\rangle$ from below when  $t\in I$ and $t\in J$, separately. 

Suppose $t\in I$ so $\bar\omega^{(t)}\leq \epsilon_t$, $\bzt^{(t)}=\bzt_f^{(t)}$ and $\bxi^{(t)}=\bxi_f^{(t)}$. By the $\rho$-weak convexity of $f$, we have  
\begin{align}
\nonumber
\left\langle\bzt^{(t)},\bx^{(t)}-\widehat\bx^{(t)}\right\rangle
\geq&~f (\bx^{(t)} )-f (\widehat\bx^{(t)} )-\dfrac{\rho}{2} \|\widehat\bx^{(t)}-\bx^{(t)} \|^{2}\\\label{eq:mainthmeq3}
=&~f (\bx^{(t)} )-f (\widehat\bx^{(t)} )-\dfrac{\hat\rho}{2} \|\widehat\bx^{(t)}-\bx^{(t)} \|^{2}
+\dfrac{\hat\rho-\rho}{2} \|\widehat\bx^{(t)}-\bx^{(t)} \|^{2}.
\end{align}
Consider the convex optimization problem \eqref{eq:phi} with $\bx=\bx^{(t)}$. By Assumption~\ref{assume:convex}B, there exists a Lagrangian multiplier $\widehat\lambda_t\geq 0$ such that 
 $\widehat\lambda_tg(\widehat\bx^{(t)})=0$ (complementary slackness) and 
\begin{eqnarray*}
\widehat\bx^{(t)}=\argmin_{\bx\in\X}\left\{f(\bx)+\dfrac{\hat{\rho}}{2} \|\bx-\bx^{(t)} \|^{2}+\widehat\lambda_{t} g(\bx)\right\}.
\end{eqnarray*}
Since the objective function above is $(\hat\rho-\rho+\widehat\lambda_{t}\mu)$-strongly convex, we have 
\begin{align*}
f(\bx^{(t)})+\widehat\lambda_{t} g(\bx^{(t)})
=&~ f(\bx^{(t)})+\dfrac{\hat{\rho}}{2} \|\bx^{(t)}-\bx^{(t)} \|^{2}+\widehat\lambda_{t} g(\bx^{(t)})\\
\geq&~f(\widehat\bx^{(t)})+\dfrac{\hat{\rho}}{2} \|\widehat\bx^{(t)}-\bx^{(t)} \|^{2}+\widehat\lambda_{t} g(\widehat\bx^{(t)})+\frac{\hat{\rho}-\rho+\widehat\lambda_{t}\mu}{2}\|\widehat\bx^{(t)}-\bx^{(t)} \|^{2},
\end{align*}
which, by the fact that $\widehat\lambda_tg(\widehat\bx^{(t)})=0$, implies
\begin{eqnarray*}
f(\bx^{(t)})-f(\widehat\bx^{(t)})-\dfrac{\hat{\rho}}{2} \|\widehat\bx^{(t)}-\bx^{(t)} \|^{2}
\geq -\widehat\lambda_{t}g(\bx^{(t)})+\frac{\hat{\rho}-\rho+\widehat\lambda_{t}\mu}{2}\|\widehat\bx^{(t)}-\bx^{(t)} \|^{2}.
\end{eqnarray*}
Applying this inequality and inequality \eqref{eq:mainthmeq3} to \eqref{eq:mainthmeq2} leads to 
\begin{align}
\nonumber
&~\eta_{t}\hat{\rho}(\hat{\rho}-\rho)  \|\widehat\bx^{(t)}-\bx^{(t)} \|^{2}-\eta_{t}\hat{\rho}\widehat\lambda_{t}g(\bx^{(t)})
+\dfrac{\eta_{t}\hat{\rho}\widehat\lambda_{t}\mu}{2} \|\widehat\bx^{(t)}-\bx^{(t)} \|^{2}\\\label{eq:mainthmeq4}
\leq&~ \varphi (\bx^{(t)})-\varphi (\bx^{(t+1)})+\dfrac{\eta_{t}^{2}\hat{\rho}}{2}\|\bxi_f^{(t)}\|^{2}
    -\eta_{t}\hat{\rho}\left\langle\bxi_f^{(t)}-\bzt_f^{(t)},\bx^{(t)}-\widehat\bx^{(t)}\right\rangle.
\end{align}

Suppose $t\in J$ so $\bar\omega^{(t)}> \epsilon_t$, $\bzt^{(t)}=\bzt_g^{(t)}$ and $\bxi^{(t)}=\bxi_g^{(t)}$. By the $\mu$-strong convexity of $g$ and the fact that $g (\widehat\bx^{(t)})\leq 0$, we have  
\begin{eqnarray*}
\left\langle\bzt^{(t)},\bx^{(t)}-\widehat\bx^{(t)}\right\rangle-\frac{\mu}{2}\|\widehat\bx^{(t)}-\bx^{(t)} \|^{2}
\geq g (\bx^{(t)} )-g (\widehat\bx^{(t)}) \geq g (\bx^{(t)} ).
\end{eqnarray*}
Applying this inequality to \eqref{eq:mainthmeq2} leads to 
\begin{align}
\nonumber
&~\eta_{t}\hat{\rho}g (\bx^{(t)} )+\frac{\eta_{t}\hat{\rho}\mu}{2}\|\widehat\bx^{(t)}-\bx^{(t)} \|^{2}\\\label{eq:mainthmeq5}
\leq&~\varphi (\bx^{(t)})-\varphi (\bx^{(t+1)})+\dfrac{\eta_{t}^{2}\hat{\rho}}{2}\|\bxi_g^{(t)}\|^{2} -\eta_{t}\hat{\rho}\left\langle\bxi_g^{(t)}-\bzt_g^{(t)},\bx^{(t)}-\widehat\bx^{(t)}\right\rangle.
\end{align}

Recall that  $\mathbb{I}(t\in I)=\mathbb{I}(\bar\omega^{(t)}\leq \epsilon_t)$ and $\mathbb{I}(t\in J)=\mathbb{I}(\bar\omega^{(t)}> \epsilon_t)$. Summing up \eqref{eq:mainthmeq4} and \eqref{eq:mainthmeq5} for $t=S,S+1,\dots,T-1$, we have
\small
\begin{align*}
&~\sum_{t=S}^{T-1}\left[\eta_{t}\hat{\rho}\widehat\lambda_t\mathbb{I}(\bar\omega^{(t)}\leq \epsilon_t)+\eta_{t}\hat{\rho}\mathbb{I}(\bar\omega^{(t)}> \epsilon_t)\right]\frac{\mu}{2}\|\widehat\bx^{(t)}-\bx^{(t)} \|^{2}\\
&+\sum_{t=S}^{T-1}\left[\eta_{t}\hat{\rho}(\hat{\rho}-\rho)\|\widehat\bx^{(t)}-\bx^{(t)} \|^{2}-\eta_{t}\hat{\rho}\widehat\lambda_{t}g (\bx^{(t)} )\right]\mathbb{I}(\bar\omega^{(t)}\leq \epsilon_t)+\sum_{t=S}^{T-1}\eta_{t}\hat{\rho}g (\bx^{(t)} )\mathbb{I}(\bar\omega^{(t)}> \epsilon_t)\\
\leq&~ \varphi (\bx^{(S)})-\varphi (\bx^{(T)})
+\dfrac{\hat{\rho}}{2}\sum_{t=S}^{T-1}\left[\eta_{t}^{2}\|\bxi_f^{(t)}\|^{2}\mathbb{I}(\bar\omega^{(t)}\leq \epsilon_t)+\eta_{t}^{2}\|\bxi_g^{(t)}\|^{2}\mathbb{I}(\bar\omega^{(t)}> \epsilon_t)\right]\\
&-\sum_{t=S}^{T-1}\left[\eta_{t}\hat{\rho}\left\langle\bxi_f^{(t)}-\bzt_f^{(t)},\bx^{(t)}-\widehat\bx^{(t)}\right\rangle\mathbb{I}(\bar\omega^{(t)}\leq \epsilon_t)+\eta_{t}\hat{\rho}\left\langle\bxi_g^{(t)}-\bzt_g^{(t)},\bx^{(t)}-\widehat\bx^{(t)}\right\rangle\mathbb{I}(\bar\omega^{(t)}> \epsilon_t)\right].
\end{align*}
\normalsize
The  conclusion is then implied by the facts that $\varphi (\bx^{(T)})=f(\widehat\bx^{(T)})+\frac{\hat\rho}{2}\|\widehat\bx^{(T)}-\bx^{(T)}\|^2\geq f(\widehat\bx^{(T)})$ and $\varphi (\bx^{(S)})= f(\widehat\bx^{(S)})+\frac{\hat\rho}{2}\|\widehat\bx^{(S)}-\bx^{(S)}\|^2\leq f(\widehat\bx^{(T)})+\frac{\hat\rho}{2}\|\widehat\bx^{(T)}-\bx^{(S)}\|^2\leq f(\widehat\bx^{(T)})+\frac{\hat\rho D^2}{2}.$
\end{proof}

\begin{lemma}
\label{lemma:boundomega}
Suppose Assumption~\ref{assume:stochasticgrad_prob} holds. Given any $\bx\in\X$, let $\{\omega_i\}_{i=1}^B$ be a mini-batch of stochastic estimators of $g$ at $\bx$ and $\bar\omega=\frac{1}{B}\sum_{i=1}^B\omega_i$. It holds that, for any $\delta\in(0,1)$,
\begin{align*}
\textup{Prob}\left(\bar\omega>g(\bx)+\sqrt{\frac{3}{B}}\sigma\sqrt{\ln(1/\delta)}\right)\leq \delta
~\text{ and }~
\textup{Prob}\left(\bar\omega<g(\bx)-\sqrt{\frac{3}{B}}\sigma\sqrt{\ln(1/\delta)}\right)\leq \delta.
\end{align*}
\end{lemma}
\begin{proof}
The conclusion is guaranteed by Assumption~\ref{assume:stochasticgrad_prob} and Lemma 2 (Case A) in \cite{lan2012validation} by choosing $\Omega=\sqrt{3\ln(1/\delta)}$  in their bound.
\end{proof}

\begin{lemma}
\label{lemma:boundgx}
Suppose Assumption~\ref{assume:stochasticgrad_prob} holds. For any $\delta\in (0,1)$, Algorithm~\ref{alg:sgm} guarantees with probability at least $1-\delta/(T-S)$ that
\begin{align}
\label{eq:highprobgx}
g (\bx^{(t)} )\mathbb{I}(\bar\omega^{(t)}\leq \epsilon_t)
\leq&\left(\epsilon_t+\sqrt{\frac{3}{B}}\sigma\sqrt{\ln((T-S)/\delta)}\right)\mathbb{I}(\bar\omega^{(t)}\leq \epsilon_t)
\end{align}
for $t=S,S+1,\dots,T-1$ and, consequently, with probability at least $1-\delta$ that
\begin{align*}
    \sum_{t=S}^{T-1}\eta_{t}\hat{\rho}\widehat\lambda_{t}g (\bx^{(t)} )\mathbb{I}(\bar\omega^{(t)}\leq \epsilon_t)
    \leq \sum_{t=S}^{T-1}\eta_{t}\hat{\rho}\widehat\lambda_{t} \left(\epsilon_t+\sqrt{\frac{3}{B}}\sigma\sqrt{\ln((T-S)/\delta)}\right)\mathbb{I}(\bar\omega^{(t)}\leq \epsilon_t).
\end{align*}
Similarly, Algorithm~\ref{alg:sgm} guarantees with probability at least $1-\delta$ that
\begin{align*}
    \sum_{t=S}^{T-1}\eta_{t}\hat{\rho}g (\bx^{(t)} )\mathbb{I}(\bar\omega^{(t)}> \epsilon_t)
    \geq
    \sum_{t=S}^{T-1}\eta_{t}\hat{\rho}\left(\epsilon_t-\sqrt{\frac{3}{B}}\sigma\sqrt{\ln((T-S)/\delta)}\right)
    \mathbb{I}(\bar\omega^{(t)}> \epsilon_t).
\end{align*}
\end{lemma}
\begin{proof}
For any $t$, $\bx^{(t)}$ is determined by $\bar\omega^{(0)},\bxi_f^{(0)},\bxi_g^{(0)},\bar\omega^{(1)},\bxi_f^{(1)},\bxi_g^{(1)},\dots,\bar\omega^{(t-1)},\bxi_f^{(t-1)}$ and $\bxi_g^{(t-1)}$, while $\bar\omega^{(t)}$ is generated after $\bx^{(t)}$. Hence,  according to Lemma~\ref{lemma:boundomega}, we have  with a probability of at least $1-\delta/(T-S)$ that
$g (\bx^{(t)} )\mathbb{I}(\bar\omega^{(t)}\leq \epsilon_t)
\leq\left(\bar\omega^{(t)}+\sqrt{\frac{3}{B}}\sigma\sqrt{\ln((T-S)/\delta)}\right)\mathbb{I}(\bar\omega^{(t)}\leq \epsilon_t)$, which implies \eqref{eq:highprobgx}. The first conclusion is then obtained by taking the union bound for the events above for $t=S,S+1,\dots,T-1$. The second  conclusion can be proved in a similar way.
\end{proof}

\begin{lemma}
\label{lemma:boundxinorm}
Suppose Assumption~\ref{assume:stochasticgrad_prob} holds. For any $\delta\in (0,1)$, Algorithm~\ref{alg:sgm} guarantees with probability at least $1-\delta$ that
\begin{align*}
&~\sum_{t=S}^{T-1}\left[\eta_{t}^{2}\|\bxi_f^{(t)}\|^{2}\mathbb{I}(\bar\omega^{(t)}\leq \epsilon_t)+\eta_{t}^{2}\|\bxi_g^{(t)}\|^{2}\mathbb{I}(\bar\omega^{(t)}> \epsilon_t)\right]\\
\leq&~\sum_{t=S}^{T-1}\eta_t^2M^2+\max\left\{\sqrt{12\ln(2/\delta)},\frac{4}{3}\ln(2/\delta)\right\}\sqrt{\sum_{t=S}^{T-1}\eta_t^4M^4}.
\end{align*}
\end{lemma}
\begin{proof}
For any $t$, $\bx^{(t)}$ is determined by $\bar\omega^{(0)},\bxi_f^{(0)},\bxi_g^{(0)},\bar\omega^{(1)},\bxi_f^{(1)},\bxi_g^{(1)},\dots,\bar\omega^{(t-1)},\bxi_f^{(t-1)}$ and $\bxi_g^{(t-1)}$. Also, $\bar\omega^{(t)}$, $\bxi_f^{(t)}$ and $\bxi_g^{(t)}$ are independent and generated after $\bx^{(t)}$. Hence, by Assumption~\ref{assume:stochasticgrad_prob}, we have
\begin{align*}
&\mathbb{E}_t\left[\exp\left(\frac{\eta_{t}^{2}\|\bxi_f^{(t)}\|^{2}\mathbb{I}(\bar\omega^{(t)}\leq \epsilon_t)+\eta_{t}^{2}\|\bxi_g^{(t)}\|^{2}\mathbb{I}(\bar\omega^{(t)}> \epsilon_t)}{\eta_t^2M^2}\right)\right]\\
=~&\mathbb{E}_t\left[\mathbb{I}(\bar\omega^{(t)}\leq \epsilon_t)\exp\left(\frac{\|\bxi_f^{(t)}\|^{2}}{M^2}\right)\right]+\mathbb{E}_t\left[\mathbb{I}(\bar\omega^{(t)}> \epsilon_t)\exp\left(\frac{\|\bxi_g^{(t)}\|^{2}}{M^2}\right)\right]\\
=~&\mathbb{E}_t\left[\mathbb{I}(\bar\omega^{(t)}\leq \epsilon_t)\right]\mathbb{E}_t\left[\exp\left(\frac{\|\bxi_f^{(t)}\|^{2}}{M^2}\right)\right]+\mathbb{E}_t\left[\mathbb{I}(\bar\omega^{(t)}> \epsilon_t)\right]\mathbb{E}_t\left[\exp\left(\frac{\|\bxi_g^{(t)}\|^{2}}{M^2}\right)\right]\\
\leq~&\mathbb{E}_t\left[\mathbb{I}(\bar\omega^{(t)}\leq \epsilon_t)\right]\exp\left(1\right)+\mathbb{E}_t\left[\mathbb{I}(\bar\omega^{(t)}> \epsilon_t)\right]\exp\left(1\right)=\exp\left(1\right),
\end{align*}
where the second equality is by the conditional independence between $\bar\omega^{(t)}$, $\bxi_f^{(t)}$ and $\bxi_g^{(t)}$. Then the  conclusion is guaranteed by Lemma 2 (Case B) in \cite{lan2012validation} by choosing $\Omega=\max\{\sqrt{12\ln(2/\delta)},(4/3)\cdot\ln(2/\delta)\}$ in their bound. 
\end{proof}

\begin{lemma}
\label{lemma:boundxidiff}
Suppose Assumptions~\ref{assume:allpaper} and~\ref{assume:stochasticgrad_prob} hold. For any $\delta\in (0,1)$, Algorithm~\ref{alg:sgm} guarantees with probability at least $1-\delta$ that
\begin{align*}
&~\sum_{t=S}^{T-1}\left[\eta_{t}\hat{\rho}\left\langle\bxi_f^{(t)}-\bzt_f^{(t)},\widehat\bx^{(t)}-\bx^{(t)}\right\rangle\mathbb{I}(\bar\omega^{(t)}\leq \epsilon_t)+\eta_{t}\hat{\rho}\left\langle\bxi_g^{(t)}-\bzt_g^{(t)},\widehat\bx^{(t)}-\bx^{(t)}\right\rangle\mathbb{I}(\bar\omega^{(t)}> \epsilon_t)\right]\\
\leq&~\sqrt{3\ln(1/\delta)}\sqrt{\sum_{t=S}^{T-1}4\eta_{t}^2\hat{\rho}^2M^2D^2}.
\end{align*}
\end{lemma}
\begin{proof}
By Assumption~\ref{assume:stochasticgrad_prob} and Jensen's inequality, we have
\small
\begin{align*}
    &~\exp\left(\frac{\|\bzt_f^{(t)}\|^2}{M^2}\right)\leq \mathbb{E}_t\left[\exp\left(\frac{\|\bxi_f^{(t)}\|^2}{M^2}\right)\right]\leq \exp(1)\\
    ~\text{and}~&~
    \exp\left(\frac{\|\bzt_g^{(t)}\|^2}{M^2}\right)\leq \mathbb{E}_t\left[\exp\left(\frac{\|\bxi_g^{(t)}\|^2}{M^2}\right)\right]\leq \exp(1),
\end{align*}
\normalsize
which, by Jensen's inequality again, implies
\small
\begin{align}
\nonumber
&~\mathbb{E}_t\left[\exp\left(\frac{2\|\bxi_f^{(t)}\|^{2}+2\|\bzt_f^{(t)}\|^{2}}{4M^2}\right)\right]\\\label{eq:expboundbzt1}
=&~
\mathbb{E}_t\left[\exp\left(\frac{\|\bxi_f^{(t)}\|^{2}}{2M^2}\right)\exp\left(\frac{\|\bzt_f^{(t)}\|^{2}}{2M^2}\right)\right]
\leq 
\sqrt{\mathbb{E}_t\left[\exp\left(\frac{\|\bxi_f^{(t)}\|^{2}}{M^2}\right)\right]}\exp\left(\frac{1}{2}\right)
\leq \exp(1)\\
\nonumber
~\text{and}\quad&~\mathbb{E}_t\left[\exp\left(\frac{2\|\bxi_g^{(t)}\|^{2}+2\|\bzt_g^{(t)}\|^{2}}{4M^2}\right)\right]\\\label{eq:expboundbzt2}
=&~
\mathbb{E}_t\left[\exp\left(\frac{\|\bxi_g^{(t)}\|^{2}}{2M^2}\right)\exp\left(\frac{\|\bzt_g^{(t)}\|^{2}}{2M^2}\right)\right]
\leq 
\sqrt{\mathbb{E}_t\left[\exp\left(\frac{\|\bxi_g^{(t)}\|^{2}}{M^2}\right)\right]}\exp\left(\frac{1}{2}\right)
\leq \exp(1).
\end{align}
\normalsize
For any $t$, $\bx^{(t)}$ is determined by $\bar\omega^{(0)},\bxi_f^{(0)},\bxi_g^{(0)},\bar\omega^{(1)},\bxi_f^{(1)},\bxi_g^{(1)},\dots,\bar\omega^{(t-1)},\bxi_f^{(t-1)}$ and $\bxi_g^{(t-1)}$. Also, $\bar\omega^{(t)}$, $\bxi_f^{(t)}$ and $\bxi_g^{(t)}$ are independent and generated after $\bx^{(t)}$. Hence, by Assumption~\ref{assume:stochasticgrad_prob}, we have
\small
\begin{align*}
\mathbb{E}_t\left[\eta_{t}\hat{\rho}\left\langle\bxi_f^{(t)}-\bzt_f^{(t)},\widehat\bx^{(t)}-\bx^{(t)}\right\rangle\mathbb{I}(\bar\omega^{(t)}\leq \epsilon_t)+\eta_{t}\hat{\rho}\left\langle\bxi_g^{(t)}-\bzt_g^{(t)},\widehat\bx^{(t)}-\bx^{(t)}\right\rangle\mathbb{I}(\bar\omega^{(t)}> \epsilon_t)\right]=0,
\end{align*}
\normalsize
and
\small
\begin{align*}
&~\mathbb{E}_t\left[\exp\left(\frac{\left[\eta_{t}\hat{\rho}\left\langle\bxi_f^{(t)}-\bzt_f^{(t)},\widehat\bx^{(t)}-\bx^{(t)}\right\rangle\mathbb{I}(\bar\omega^{(t)}\leq \epsilon_t)+\eta_{t}\hat{\rho}\left\langle\bxi_g^{(t)}-\bzt_g^{(t)},\widehat\bx^{(t)}-\bx^{(t)}\right\rangle\mathbb{I}(\bar\omega^{(t)}> \epsilon_t)\right]^2}{4\eta_{t}^2\hat{\rho}^2M^2D^2}\right)\right]\\
=&~\mathbb{E}_t\left[\mathbb{I}(\bar\omega^{(t)}\leq \epsilon_t)\exp\left(\frac{\left[\left\langle\bxi_f^{(t)}-\bzt_f^{(t)},\widehat\bx^{(t)}-\bx^{(t)}\right\rangle\right]^{2}}{4M^2D^2}\right)\right]+\mathbb{E}_t\left[\mathbb{I}(\bar\omega^{(t)}> \epsilon_t)\exp\left(\frac{\left[\left\langle\bxi_g^{(t)}-\bzt_g^{(t)},\widehat\bx^{(t)}-\bx^{(t)}\right\rangle\right]^{2}}{4M^2D^2}\right)\right]\\
\leq&~\mathbb{E}_t\left[\mathbb{I}(\bar\omega^{(t)}\leq \epsilon_t)\right]\mathbb{E}_t\left[\exp\left(\frac{\|\bxi_f^{(t)}-\bzt_f^{(t)}\|^{2}\|\widehat\bx^{(t)}-\bx^{(t)}\|^{2}}{4M^2D^2}\right)\right]\\
&+\mathbb{E}_t\left[\mathbb{I}(\bar\omega^{(t)}> \epsilon_t)\right]\mathbb{E}_t\left[\exp\left(\frac{\|\bxi_g^{(t)}-\bzt_g^{(t)}\|^{2}\|\widehat\bx^{(t)}-\bx^{(t)}\|^{2}}{4M^2D^2}\right)\right]\\
\leq&~\mathbb{E}_t\left[\mathbb{I}(\bar\omega^{(t)}\leq \epsilon_t)\right]\mathbb{E}_t\left[\exp\left(\frac{2\|\bxi_f^{(t)}\|^{2}+2\|\bzt_f^{(t)}\|^{2}}{4M^2}\right)\right]+\mathbb{E}_t\left[\mathbb{I}(\bar\omega^{(t)}> \epsilon_t)\right]\mathbb{E}_t\left[\exp\left(\frac{2\|\bxi_g^{(t)}\|^{2}+2\|\bzt_g^{(t)}\|^{2}}{4M^2}\right)\right]\\
\leq&~\mathbb{E}_t\left[\mathbb{I}(\bar\omega^{(t)}\leq \epsilon_t)\right]\exp\left(1\right)+\mathbb{E}_t\left[\mathbb{I}(\bar\omega^{(t)}> \epsilon_t)\right]\exp\left(1\right)=\exp\left(1\right),
\end{align*}
\normalsize
where the first inequality is by the Cauchy–Schwarz inequality and the conditional independence between $\bar\omega^{(t)}$, $\bxi_f^{(t)}$ and $\bxi_g^{(t)}$, the second inequality is by Assumption~\ref{assume:allpaper}B and the last inequality is by \eqref{eq:expboundbzt1} and \eqref{eq:expboundbzt2}. Then the conclusion is guaranteed by Lemma 2 (Case A) in \cite{lan2012validation} by choosing $\Omega=\sqrt{3\ln(1/\delta)}$ in their bound. 
\end{proof}

Suppose $g$ is convex but not strongly convex, i.e., $\mu=0$. Taking the union bound of the four events in Lemmas~\ref{lemma:boundgx},~\ref{lemma:boundxinorm} and~\ref{lemma:boundxidiff} with $\delta$ replaced by $\delta/4$ and applying the four inequalities (two from Lemma~\ref{lemma:boundgx}, one from Lemma~\ref{lemma:boundxinorm} and one from Lemma~\ref{lemma:boundxidiff}) holding in these  events to Proposition~\ref{thm:mainprop}, we have the following bounds. 
\begin{proposition}
\label{thm:mainprop_prob}
Suppose Assumptions~\ref{assume:allpaper}, \ref{assume:convex} and \ref{assume:stochasticgrad_prob} hold and $g$ is $\mu$-strongly convex ($\mu$ can be zero). Algorithm~\ref{alg:sgm} guarantees with probability at least $1-\delta/(4(T-S))$ that
\small
\begin{align}
\label{eq:highprobgx4}
g (\bx^{(t)} )\mathbb{I}(\bar\omega^{(t)}\leq \epsilon_t)
\leq&\left(\epsilon_t+\sqrt{\frac{3}{B}}\sigma\sqrt{\ln(4(T-S)/\delta)}\right)\mathbb{I}(\bar\omega^{(t)}\leq \epsilon_t)
\end{align}
\normalsize
for $t=S,S+1,\dots,T-1$ and, consequently, with probability at least $1-\delta$ that
\small
\begin{align}
\nonumber
&~\sum_{t=S}^{T-1}\left[\eta_{t}\hat{\rho}\widehat\lambda_t\mathbb{I}(\bar\omega^{(t)}\leq \epsilon_t)+\eta_{t}\hat{\rho}\mathbb{I}(\bar\omega^{(t)}> \epsilon_t)\right]\cdot\frac{\mu}{2}\|\widehat\bx^{(t)}-\bx^{(t)} \|^{2}\\\nonumber
&+\sum_{t=S}^{T-1}\left[\eta_{t}\hat{\rho}(\hat{\rho}-\rho)\|\widehat\bx^{(t)}-\bx^{(t)} \|^{2}-\eta_{t}\hat{\rho}\widehat\lambda_{t}\epsilon_t\right]\mathbb{I}(\bar\omega^{(t)}\leq \epsilon_t)+\sum_{t=S}^{T-1}\eta_{t}\hat{\rho}\epsilon_t\mathbb{I}(\bar\omega^{(t)}> \epsilon_t)\\\label{eq:propeq1}
\leq&~ \frac{\hat\rho D^2}{2}
+\frac{\hat\rho}{2}\sum_{t=S}^{T-1}\eta_t^2M^2+\frac{\hat\rho}{2}\max\left\{\sqrt{12\ln(8/\delta)},\frac{4}{3}\ln(8/\delta)\right\}\sqrt{\sum_{t=S}^{T-1}\eta_t^4M^4}\\\nonumber
&+\sqrt{3\ln(4/\delta)}\sqrt{\sum_{t=S}^{T-1}4\eta_{t}^2\hat{\rho}^2M^2D^2}+\sum_{t=S}^{T-1}\eta_{t}\hat{\rho}\sqrt{\frac{3}{B}}\sigma\sqrt{\ln(4(T-S)/\delta)}\left[\widehat\lambda_{t}\mathbb{I}(\bar\omega^{(t)}\leq \epsilon_t)+\mathbb{I}(\bar\omega^{(t)}> \epsilon_t)\right].
\end{align}
\normalsize
\end{proposition}

\subsection{Proof of Theorem~\ref{thm:mdconverge_deterministic}}
\label{sec:mdconverge_deterministic}
Although the technical results in the previous sections are derived for Algorithm~\ref{alg:sgm} in the stochastic case, they all apply to Algorithm~\ref{alg:dsgm} which is a special case of  Algorithm~\ref{alg:sgm} when all oracles are deterministic, i.e., when \eqref{eq:reduced_to_deterministic} holds. Hence, we can uses these technical results to prove Theorem~\ref{thm:mdconverge_deterministic} as follows.  

\begin{proof}[Proof of Theorem~\ref{thm:mdconverge_deterministic}]
Because Algorithm~\ref{alg:dsgm} is a deterministic special case of Algorithm~\ref{alg:sgm} when \eqref{eq:reduced_to_deterministic} holds and also  because we do not assume strong convexity in $g$ ($\mu=0$),  we can simplify the inequality in Proposition~\ref{thm:mainprop} as follows
\small
\begin{align}
\nonumber
&~\sum_{t=S}^{T-1}\left[\eta_{t}\hat{\rho}(\hat{\rho}-\rho)\|\widehat\bx^{(t)}-\bx^{(t)} \|^{2}-\eta_{t}\hat{\rho}\widehat\lambda_{t}g (\bx^{(t)} )\right]\mathbb{I}(g (\bx^{(t)} )\leq \epsilon_t)+\sum_{t=S}^{T-1}\eta_{t}\hat{\rho}g (\bx^{(t)} )\mathbb{I}(g (\bx^{(t)} )> \epsilon_t)\\\nonumber
\leq&~\frac{\hat\rho D^2}{2}
+\dfrac{\hat{\rho}}{2}\sum_{t=S}^{T-1}\left[\eta_{t}^{2}\|\bzt_f^{(t)}\|^{2}\mathbb{I}(g (\bx^{(t)} )\leq \epsilon_t)+\eta_{t}^{2}\|\bzt_g^{(t)}\|^{2}\mathbb{I}(g (\bx^{(t)} )> \epsilon_t)\right].
\end{align}
\normalsize
By the facts that $g (\bx^{(t)} )\mathbb{I}(g (\bx^{(t)} )\leq \epsilon_t)\leq \epsilon_t\mathbb{I}(g (\bx^{(t)} )\leq \epsilon_t)$ and $g (\bx^{(t)} )\mathbb{I}(g (\bx^{(t)} )> \epsilon_t)\geq \epsilon_t\mathbb{I}(g (\bx^{(t)} )> \epsilon_t)$, the inequality above implies 
\small
\begin{align}
\nonumber
&~\sum_{t=S}^{T-1}\left[\eta_{t}\hat{\rho}(\hat{\rho}-\rho)\|\widehat\bx^{(t)}-\bx^{(t)} \|^{2}-\eta_{t}\hat{\rho}\widehat\lambda_{t}\epsilon_t\right]\mathbb{I}(g (\bx^{(t)} )\leq \epsilon_t)+\sum_{t=S}^{T-1}\eta_{t}\hat{\rho}\epsilon_t\mathbb{I}(g (\bx^{(t)} )> \epsilon_t)\\\nonumber
\leq&~\frac{\hat\rho D^2}{2}
+\dfrac{\hat{\rho}}{2}\sum_{t=S}^{T-1}\left[\eta_{t}^{2}\|\bzt_f^{(t)}\|^{2}\mathbb{I}(g (\bx^{(t)} )\leq \epsilon_t)+\eta_{t}^{2}\|\bzt_g^{(t)}\|^{2}\mathbb{I}(g (\bx^{(t)} )> \epsilon_t)\right]\\\label{eq:propeq10}
\leq&~\frac{\hat\rho D^2}{2}
+\dfrac{\hat{\rho}}{2}\sum_{t=S}^{T-1}\eta_{t}^{2}M^2,
\end{align}
\normalsize
where the second inequality is because of Assumption~\ref{assume:allpaper}.

We first prove that, if $S$, $T$, $\eta_t$ and $\epsilon_t$ are chosen such that 
\begin{eqnarray}
\label{eq:propeq4}
\epsilon_t(1+\widehat\lambda_t)\leq \epsilon^2(\hat\rho-\rho)
\end{eqnarray}
and
\begin{eqnarray}
\label{eq:propeq9}
\sum_{t=S}^{T-1}\eta_{t}\hat{\rho}\epsilon_t> \frac{\hat\rho D^2}{2}
+\frac{\hat\rho}{2}\sum_{t=S}^{T-1}\eta_t^2M^2,
\end{eqnarray}
we must have $g(\bx^{(t)})\leq\epsilon_t$ for at least one $t$ in $\{S,S+1,\dots,T-1\}$ (i.e., $I\neq\emptyset$) and $\mathbb{E}_\tau[\|\widehat\bx^{(\tau)} - \bx^{(\tau)}\|^2]\leq\epsilon^2$ (so $\mathbb{E}_\tau[\|\widehat\bx^{(\tau)} - \bx^{(\tau)}\|]\leq\epsilon$). 

Suppose \eqref{eq:propeq9} holds and $g(\bx^{(t)})>\epsilon_t$ for $t=S,S+1,\dots,T-1$, i.e., $I=\emptyset$. \eqref{eq:propeq10} becomes exactly the opposite of \eqref{eq:propeq9}. This contradiction means  $g(\bx^{(t)})\leq\epsilon_t$ for at least one $t$ in $\{S,S+1,\dots,T-1\}$. Suppose \eqref{eq:propeq4} and \eqref{eq:propeq9} hold but $\mathbb{E}_\tau[\|\widehat\bx^{(\tau)} - \bx^{(\tau)}\|^2]>\epsilon^2$. Since $\tau$ is generated by Output I, we have
\begin{align}
\label{eq:stationaryopposite}
\epsilon^2<\mathbb{E}_\tau[\|\widehat\bx^{(\tau)} - \bx^{(\tau)}\|^2]
=\frac{\sum_{t=S}^{T-1}\eta_{t}\mathbb{I}(g(\bx^{(t)})\leq\epsilon_t)\|\widehat\bx^{(t)}-\bx^{(t)} \|^{2}}{\sum_{t=S}^{T-1}\eta_{t}\mathbb{I}(g(\bx^{(t)})\leq\epsilon_t)}.
\end{align}
Note that the right-hand side of \eqref{eq:stationaryopposite} is well-defined because we just proved $I\neq\emptyset$. \eqref{eq:stationaryopposite} and \eqref{eq:propeq4} imply
\small
\begin{align}
\nonumber
&~\sum_{t=S}^{T-1}\left[\eta_{t}\hat{\rho}(\hat{\rho}-\rho)\|\widehat\bx^{(t)}-\bx^{(t)} \|^{2}-\eta_{t}\hat{\rho}\widehat\lambda_{t}\epsilon_t\right]\mathbb{I}(g(\bx^{(t)})\leq \epsilon_t)+\sum_{t=S}^{T-1}\eta_{t}\hat{\rho}\epsilon_t\mathbb{I}(g(\bx^{(t)})> \epsilon_t)\\\nonumber
>&~\sum_{t=S}^{T-1}\left[\eta_{t}\hat{\rho}(\hat{\rho}-\rho)\epsilon^{2}-\eta_{t}\hat{\rho}\widehat\lambda_{t}\epsilon_t\right]\mathbb{I}(g(\bx^{(t)})\leq \epsilon_t)+\sum_{t=S}^{T-1}\eta_{t}\hat{\rho}\epsilon_t\mathbb{I}(g(\bx^{(t)})> \epsilon_t)\\\label{eq:leftright}
\geq&~\sum_{t=S}^{T-1}\eta_{t}\hat{\rho}\epsilon_t\mathbb{I}(g(\bx^{(t)})\leq \epsilon_t)+\sum_{t=S}^{T-1}\eta_{t}\hat{\rho}\epsilon_t\mathbb{I}(g(\bx^{(t)})> \epsilon_t)
=\sum_{t=S}^{T-1}\eta_{t}\hat{\rho}\epsilon_t,
\end{align}
\normalsize
where the second inequality is  because of \eqref{eq:propeq4}. Combining this inequality and \eqref{eq:propeq10} leads to the opposite of \eqref{eq:propeq9}. This contradiction means $\mathbb{E}_\tau[\|\widehat\bx^{(\tau)} - \bx^{(\tau)}\|^2]\leq\epsilon^2$.

Given the result above, we only need to show that the two choices of $S$, $T$, $\eta_t$ and $\epsilon_t$ ensure \eqref{eq:propeq4} and \eqref{eq:propeq9}.

In Case I, \eqref{eq:propeq4} holds because of Lemma~\ref{thm:boundlambda} and the choice of $\epsilon_t$. Let $\eta=\eta_t=\frac{2\epsilon^2(\hat\rho-\rho)}{5(1+\Lambda)M^2}$ for any $t$. Using Lemma~\ref{thm:boundlambda} and plugging the values of $S$, $T$, $\eta_t$ and $\epsilon_t$ in  \eqref{eq:propeq9}, we can show that \eqref{eq:propeq9} is equivalent to
\begin{align*}
\frac{T\eta\hat\rho\epsilon^2(\hat\rho-\rho)}{1+\Lambda}>&~ \frac{\hat\rho D^2}{2}
+\frac{\hat\rho}{2}T\eta^2M^2,
\end{align*}
which, after dividing both sides by $T\eta\hat\rho$, can be equivalently written as
\begin{align*}
\frac{\epsilon^2(\hat\rho-\rho)}{1+\Lambda}>&~ \frac{D^2}{2T\eta}
+\frac{\eta M^2}{2}.
\end{align*}
By the values of $\eta$ and $T$, each summand in the right-hand side of the inequality above is no more than $\frac{\epsilon^2(\hat\rho-\rho)}{5(1+\Lambda)}$ so the right-hand side of the inequality above no more than $\frac{2\epsilon^2(\hat\rho-\rho)}{5(1+\Lambda)}$ which is strictly less than the left-hand side. This means \eqref{eq:propeq9} holds with this choice of parameters and thus $\mathbb{E}_\tau[\|\widehat\bx^{(\tau)} - \bx^{(\tau)}\|]\leq\epsilon$. By the convexity of $g$ and the choices of $\eta_t$ and $\epsilon_t$, we have $\mathbb{E}_\tau [g(\bx^{(\tau)})]\leq \frac{\epsilon^2(\hat\rho-\rho)}{1+\Lambda}$.

In Case II, by the choices of $\epsilon_t$ and $T$, we have, for any $t\in\{S,S+1,\dots,T-1\}$,  
\begin{align}
\label{eq:epsilont}
\epsilon_t =\frac{5MD}{\sqrt{t+1}}\leq \frac{5MD}{\sqrt{S+1}}=\frac{5MD}{\sqrt{T/2+1}}\leq \frac{\epsilon^2(\hat\rho-\rho)}{1+\Lambda}.
\end{align}
This further implies \eqref{eq:propeq4} because of Lemma~\ref{thm:boundlambda}. Note that $\eta_t$ and $\epsilon_t$ are decreasing in $t$. Hence, the left-hand side of  \eqref{eq:propeq9} satisfies 
\begin{eqnarray}
\label{eq:propeq5}
\sum_{t=S}^{T-1}\eta_{t}\hat{\rho}\epsilon_t>\frac{T}{2}\eta_T\hat{\rho}\epsilon_T=\frac{5T}{2T+2}\hat{\rho}D^2\geq \frac{5\hat{\rho}D^2}{4}.
\end{eqnarray}
The right-hand side of \eqref{eq:propeq9} satisfies 
\begin{eqnarray}
\label{eq:propeq11}
\frac{\hat\rho D^2}{2}
+\frac{\hat\rho}{2}\sum_{t=S}^{T-1}\eta_t^2M^2
=\frac{\hat\rho D^2}{2}
+\frac{\hat\rho}{2}D^2\sum_{t=S}^{T-1}\frac{1}{t+1}
\leq\hat\rho D^2,
\end{eqnarray}
where the equality is obtained by plugging in the definition of $\eta_t$ and the inequality is because $\sum_{t=S}^{T-1}\frac{1}{t+1}\leq \int_{S}^{T}\frac{1}{t}dt=\ln(T/S)=\ln(2)\leq 1$.  The right-hand side of \eqref{eq:propeq5} is strictly greater than the right-hand side \eqref{eq:propeq11}, which means \eqref{eq:propeq9} holds and thus $\mathbb{E}_\tau[\|\widehat\bx^{(\tau)} - \bx^{(\tau)}\|]\leq\epsilon$. By the convexity of $g$ and the choices of $\eta_t$ and $\epsilon_t$, we have $\mathbb{E}_\tau [g(\bx^{(\tau)})]\leq \epsilon_{S}\leq\frac{\epsilon^2(\hat\rho-\rho)}{1+\Lambda}$ according to \eqref{eq:epsilont}.
\end{proof}

\subsection{Complexity analysis for strongly convex constraints}
\label{sec:strongconvexg}

Suppose $f$ and $g$ in \eqref{eq:gco} are deterministic, i.e., \eqref{eq:reduced_to_deterministic} holds, and $g$ is $\mu$-strongly convex and $\mu>0$. The complexity of Algorithm~\ref{alg:dsgm} is characterized by the following theorem. 

\begin{theorem}
\label{thm:gstrongconvex}
Suppose Assumptions~\ref{assume:allpaper} and \ref{assume:convex} hold.  Let $\widehat\bx(\bx^{(t)})$ be defined as in \eqref{eq:phix} with $(\hat\rho,\tilde\rho)$ satisfying \eqref{eq:parameter1}
and $\bx^{(\tau)}$ is generated by Output II. Algorithm~\ref{alg:dsgm} guarantees 
$\mathbb{E}_\tau[\|\widehat\bx^{(\tau)} - \bx^{(\tau)}\|]\leq\epsilon$ in either of the following cases. 

Case I: $S=0$, $\epsilon_t=0$, $\eta_t=\frac{\epsilon^2\min\{\hat\rho-\rho,\mu/2\}}{M^2}$ and 
$T\geq \frac{M^2D^2}{\epsilon^4\min\{(\hat\rho-\rho)^2,\mu^2/4\}}
=O(1/\epsilon^4).$

Case II: $S=T/2$, $\epsilon_t=0$, $\eta_t=\frac{D}{M\sqrt{t+1}}$ and 
$T\geq \frac{4M^2D^2}{\epsilon^4\min\{(\hat\rho-\rho)^2,\mu^2/4\}}
=O(1/\epsilon^4).$
\end{theorem}

Before presenting its proof, we would like to make a few remarks. 
\begin{remark}
According to Theorem~\ref{thm:gstrongconvex}, strong convexity in the constraint function $g$ does not reduce the $O(1/\epsilon^4)$ complexity of the SSG method for finding a nearly $\epsilon$-stationary point. However, strong convexity brings benefit on other aspects. First, we can simply set $\epsilon_t=0$ when $g$ is strongly convex, which makes stepsize $\eta_t$  the only tuning parameter. Second, the theoretical complexity no longer depends on $\Lambda$, the upper bound of the dual variables, so it can be strictly better than the one in Theorem~\ref{thm:mdconverge_deterministic} when $\Lambda$ is very large. 
\end{remark}

\begin{remark}
Theorem~\ref{thm:mdconverge_deterministic}, \ref{thm:mdconverge} and \ref{thm:mdconverge_prob} guarantee $\mathbb{E}_\tau [g(\bx^{(\tau)})]\leq O(\epsilon^2)$. This is no longer true in Theorem~\ref{thm:gstrongconvex} because $\tau$ can take value in the index set $J$ on which $\bx^{(t)}$ can be highly infeasible. However, with $\mathbb{E}_\tau[\|\widehat\bx^{(\tau)} - \bx^{(\tau)}\|]\leq\epsilon$, we can at least derive  $\mathbb{E}_\tau [g(\bx^{(\tau)})]\leq O(\epsilon)$ from Theorem~\ref{thm:gstrongconvex}, which is good enough. See Remark~\ref{remark:feasibility} for the reason. 
\end{remark}

\begin{proof}[Proof of Theorem~\ref{thm:gstrongconvex}]
Since Algorithm~\ref{alg:dsgm} is fully deterministic, we can simplify the inequality in Proposition~\ref{thm:mainprop} as follows
\small
\begin{align}
\nonumber
&~\sum_{t=S}^{T-1}\left[\eta_{t}\hat{\rho}\widehat\lambda_t\mathbb{I}(g (\bx^{(t)} )\leq 0)+\eta_{t}\hat{\rho}\mathbb{I}(g (\bx^{(t)} )> 0)\right]\frac{\mu}{2}\|\widehat\bx^{(t)}-\bx^{(t)} \|^{2}\\\nonumber
&+\sum_{t=S}^{T-1}\left[\eta_{t}\hat{\rho}(\hat{\rho}-\rho)\|\widehat\bx^{(t)}-\bx^{(t)} \|^{2}-\eta_{t}\hat{\rho}\widehat\lambda_{t}g (\bx^{(t)} )\right]\mathbb{I}(g (\bx^{(t)} )\leq 0)+\sum_{t=S}^{T-1}\eta_{t}\hat{\rho}g (\bx^{(t)} )\mathbb{I}(g (\bx^{(t)} )> 0)\\\nonumber
\leq&~ \frac{\hat\rho D^2}{2}
+\dfrac{\hat{\rho}}{2}\sum_{t=S}^{T-1}\left[\eta_{t}^{2}\|\bzt_f^{(t)}\|^{2}\mathbb{I}(g (\bx^{(t)} )\leq 0)+\eta_{t}^{2}\|\bzt_g^{(t)}\|^{2}\mathbb{I}(g (\bx^{(t)} )> 0)\right].
\end{align}
\normalsize
After dropping some non-negative terms, the left-hand side of the inequality above can be  bounded from below as
\begin{align}
\nonumber
&~\sum_{t=S}^{T-1}\left[\eta_{t}\hat{\rho}\widehat\lambda_t\mathbb{I}(g (\bx^{(t)} )\leq 0)+\eta_{t}\hat{\rho}\mathbb{I}(g (\bx^{(t)} )> 0)\right]\frac{\mu}{2}\|\widehat\bx^{(t)}-\bx^{(t)} \|^{2}\\\nonumber
&+\sum_{t=S}^{T-1}\left[\eta_{t}\hat{\rho}(\hat{\rho}-\rho)\|\widehat\bx^{(t)}-\bx^{(t)} \|^{2}-\eta_{t}\hat{\rho}\widehat\lambda_{t}g (\bx^{(t)} )\right]\mathbb{I}(g (\bx^{(t)} )\leq 0)+\sum_{t=S}^{T-1}\eta_{t}\hat{\rho}g (\bx^{(t)} )\mathbb{I}(g (\bx^{(t)} )> 0)\\\nonumber
\geq&~\sum_{t=S}^{T-1}\eta_{t}\hat{\rho}\mathbb{I}(g (\bx^{(t)} )> 0)\frac{\mu}{2}\|\widehat\bx^{(t)}-\bx^{(t)} \|^{2}+\sum_{t=S}^{T-1}\eta_{t}\hat{\rho}(\hat{\rho}-\rho)\|\widehat\bx^{(t)}-\bx^{(t)} \|^{2}\mathbb{I}(g (\bx^{(t)} )\leq 0)\\\nonumber
\geq&~\hat\rho\min\{\hat\rho-\rho,\mu/2\}\sum_{t=S}^{T-1}\eta_t\|\widehat\bx^{(t)}-\bx^{(t)} \|^{2}.
\end{align}
Combining the two inequalities above gives
\begin{align}
\nonumber
&~\hat\rho\min\{\hat\rho-\rho,\mu/2\}\sum_{t=S}^{T-1}\eta_t\|\widehat\bx^{(t)}-\bx^{(t)} \|^{2}\\\nonumber
\leq&~ \frac{\hat\rho D^2}{2}
+\dfrac{\hat{\rho}}{2}\sum_{t=S}^{T-1}\left[\eta_{t}^{2}\|\bxi_f^{(t)}\|^{2}\mathbb{I}(g (\bx^{(t)} )\leq 0)+\eta_{t}^{2}\|\bzt_g^{(t)}\|^{2}\mathbb{I}(g (\bx^{(t)} )> 0)\right]\\\label{eq:propeq10_new}
\leq&~ \frac{\hat\rho D^2}{2}
+\dfrac{\hat{\rho}}{2}\sum_{t=S}^{T-1}\eta_{t}^{2}M^2,
\end{align}
where the last inequality is because of Assumption~\ref{assume:allpaper}. Since $\tau$ is generated by Output II, after organizing terms, we have
\begin{eqnarray}
\label{eq:propeq11_new}
\mathbb{E}_\tau[\|\widehat\bx^{(t)}-\bx^{(t)} \|^{2}]\leq 
\frac{1}{\min\{\hat\rho-\rho,\mu/2\}}\left(\frac{D^2}{2}
+\dfrac{1}{2}\sum_{t=S}^{T-1}\eta_{t}^{2}M^2\right)\Bigg/\left(\sum_{t=S}^{T-1}\eta_{t}\right).
\end{eqnarray}

In Case I, let $\eta=\eta_t=\epsilon^2\min\{\hat\rho-\rho,\mu/2\}/M^2$ for any $t$. By the choices of $\eta_t$, $S$ and $T$, \eqref{eq:propeq11_new} implies
\begin{eqnarray*}
\mathbb{E}_\tau[\|\widehat\bx^{(t)}-\bx^{(t)} \|^{2}]\leq 
\frac{1}{\min\{\hat\rho-\rho,\mu/2\}}\left(\frac{D^2}{2T\eta}
+\dfrac{\eta M^2}{2}\right)\leq\frac{\epsilon^2}{2}+\frac{\epsilon^2}{2}=\epsilon^2.
\end{eqnarray*}

In Case II, by the choices of $\eta_t$, $S$ and $T$, \eqref{eq:propeq11_new} implies 
\begin{eqnarray*}
\mathbb{E}_\tau[\|\widehat\bx^{(t)}-\bx^{(t)} \|^{2}]\leq 
\frac{1}{\min\{\hat\rho-\rho,\mu/2\}}\left(\frac{DM}{\sqrt{T}}
+\frac{DM}{\sqrt{T}}\right)\leq\frac{\epsilon^2}{2}+\frac{\epsilon^2}{2}=\epsilon^2,
\end{eqnarray*}
where, in the first inequality, we use the facts that $\sum_{t=S}^{T-1}\eta_{t}\geq \frac{T}{2}\eta_{T-1}=\frac{\sqrt{T}D}{2M}$ and that $\sum_{t=S}^{T-1}\eta_t^2=\sum_{t=S}^{T-1}\frac{D^2}{M^2(t+1)}\leq \int_{S}^{T}\frac{D^2}{M^2t}dt=\frac{D^2}{M^2}\ln(T/S)=\frac{D^2}{M^2}\ln(2)\leq \frac{D^2}{M^2}$, and the second inequality is by the choice of $T$.
\end{proof}

\subsection{Complexity analysis when subgradient oracle is stochastic}
\label{sec:stochastic}



In this section, we will analyze the oracle complexity of Algorithm~\ref{alg:sgm} under Assumption~\ref{assume:stochasticgrad_prob} with an additional condition that 
\begin{eqnarray}
\label{eq:reduced_to_semi_deterministic}
\bar\omega^{(t)}=g(\bx^{(t)}), ~\sigma=0~\text{ and }~B=1,
\end{eqnarray}
which is a relaxation of \eqref{eq:reduced_to_deterministic}. In this case, the function value of $g$ is can be accessed deterministically while the subgradients of $f$ and $g$ are still stochastic. This holds in the setting of zeroth-order or derivative-free optimization where one can evaluate $f$ and $g$ exactly as black boxes and then use their values to construct stochastic approximate subgradients. We will show that the complexity of Algorithm~\ref{alg:sgm} when \eqref{eq:reduced_to_semi_deterministic} holds is the same as Algorithm~\ref{alg:dsgm} (equivalent to Algorithm~\ref{alg:sgm} when \eqref{eq:reduced_to_deterministic} holds). For simplicity of notation, we denote $\widehat\bx(\bx^{(t)})$ by $\widehat\bx^{(t)}$.


\begin{theorem}
\label{thm:mdconverge}
Suppose Assumptions~\ref{assume:allpaper}, \ref{assume:convex} and~\ref{assume:stochasticgrad_prob} hold and $\Lambda$ is as in \eqref{eq:Lambdabound}. Also, suppose \eqref{eq:reduced_to_semi_deterministic} holds. Let $\widehat\bx(\bx^{(t)})$ be defined as in \eqref{eq:phix} with $(\hat\rho,\tilde\rho)$ satisfying \eqref{eq:parameter1}
and $\bx^{(\tau)}$ is generated by Algorithm~\ref{alg:sgm}. Algorithm~\ref{alg:sgm} guarantees 
$\mathbb{E}_\tau[\|\widehat\bx^{(\tau)} - \bx^{(\tau)}\|]\leq\epsilon$  and $\mathbb{E}_\tau [g(\bx^{(\tau)})]\leq \frac{\epsilon^2(\hat\rho-\rho)}{1+\Lambda}$ with probability at least $1-\delta$ in either of the following cases.

Case I: $S$, $\epsilon_t$ and $\eta_t$ are chosen as Case I in Theorem~\ref{thm:mdconverge_deterministic} and 
\small
\begin{align*}
T\geq&~ \left\{
\frac{25M^2D^2(1+\Lambda)^2}{4\epsilon^4(\hat\rho-\rho)^2},~
\max\Big\{12\ln(8/\delta),\frac{16}{9}\ln^2(8/\delta)\Big\}, ~
\frac{300\ln(4/\delta)M^2D^2(1+\Lambda)^2}{\epsilon^4(\hat\rho-\rho)^2}
\right\}
=O(1/\epsilon^4).
\end{align*}
\normalsize

Case II: $S$ and $\eta_t$ are chosen as Case II in Theorem~\ref{thm:mdconverge_deterministic}, $\epsilon_t=\frac{EMD}{\sqrt{t+1}}$
where $E$ is any positive constant such that 
\small
\begin{align*}
E\geq &~ 4+\frac{2\pi}{\sqrt{6}}\max\left\{\sqrt{12\ln(8/\delta)},\frac{4}{3}\ln(8/\delta)\right\}+8\sqrt{3\ln(4/\delta)}
\end{align*}
\normalsize
and 
$T\geq  \frac{2E^2M^2D^2(1+\Lambda)^2}{\epsilon^4(\hat\rho-\rho)^2}   
=O(1/\epsilon^4)$.
\end{theorem}

This theorem shows that, the complexity remains $O(1/\epsilon^4)$ if the subgradient oracles are stochastic but the function value oracles remain deterministic. The only difference is that the result holds in a high probability. This complexity matches the lower-bound complexity for stochastic smooth non-convex unconstrained optimization~\cite{arjevani2022lower,drori2020complexity}, so it is optimal. 

\begin{proof}[Proof of Theorem~\ref{thm:mdconverge}]
By \eqref{eq:reduced_to_semi_deterministic} and the fact that $\mu=0$, we can drop the last term in the right-hand side of \eqref{eq:propeq1} and have
\begin{align}
\nonumber
&~\sum_{t=S}^{T-1}\left[\eta_{t}\hat{\rho}(\hat{\rho}-\rho)\|\widehat\bx^{(t)}-\bx^{(t)} \|^{2}-\eta_{t}\hat{\rho}\widehat\lambda_{t}\epsilon_t\right]\mathbb{I}(g(\bx^{(t)})\leq \epsilon_t)+\sum_{t=S}^{T-1}\eta_{t}\hat{\rho}\epsilon_t\mathbb{I}(g(\bx^{(t)})> \epsilon_t)\\\label{eq:propeq2}
\leq&~ \frac{\hat\rho D^2}{2}
+\frac{\hat\rho}{2}\sum_{t=S}^{T-1}\eta_t^2M^2+\frac{\hat\rho}{2}\max\left\{\sqrt{12\ln(8/\delta)},\frac{4}{3}\ln(8/\delta)\right\}\sqrt{\sum_{t=S}^{T-1}\eta_t^4M^4}\\\nonumber
&+\sqrt{3\ln(4/\delta)}\sqrt{\sum_{t=S}^{T-1}4\eta_{t}^2\hat{\rho}^2M^2D^2}.
\end{align}
with a probability of at least $1-\delta$. 
In the rest of the proof, we always assume \eqref{eq:propeq2} holds. We first prove that, if $S$, $T$, $\eta_t$ and $\epsilon_t$ are chosen such that \eqref{eq:propeq4} holds 
and
\begin{align}
\label{eq:propeq3}
\sum_{t=S}^{T-1}\eta_{t}\hat{\rho}\epsilon_t>&~\frac{\hat\rho D^2}{2}
+\frac{\hat\rho}{2}\sum_{t=S}^{T-1}\eta_t^2M^2+\frac{\hat\rho}{2}\max\left\{\sqrt{12\ln(8/\delta)},\frac{4}{3}\ln(8/\delta)\right\}\sqrt{\sum_{t=S}^{T-1}\eta_t^4M^4}\\\nonumber
&+\sqrt{3\ln(4/\delta)}\sqrt{\sum_{t=S}^{T-1}4\eta_{t}^2\hat{\rho}^2M^2D^2},
\end{align}
we must have $g(\bx^{(t)})\leq\epsilon_t$ for at least one $t$ in $\{S,S+1,\dots,T-1\}$ (i.e., $I\neq\emptyset$) and $\mathbb{E}_\tau[\|\widehat\bx^{(\tau)} - \bx^{(\tau)}\|^2]\leq\epsilon^2$ (so $\mathbb{E}_\tau[\|\widehat\bx^{(\tau)} - \bx^{(\tau)}\|]\leq\epsilon$). 

Suppose \eqref{eq:propeq3} holds and $g(\bx^{(t)})>\epsilon_t$ for $t=S,S+1,\dots,T-1$, i.e., $I=\emptyset$. \eqref{eq:propeq2} becomes exactly the opposite of \eqref{eq:propeq3}. This contradiction means  $g(\bx^{(t)})\leq\epsilon_t$ for at least one $t$ in $\{S,S+1,\dots,T-1\}$. Suppose \eqref{eq:propeq4} and \eqref{eq:propeq3} hold but $\mathbb{E}_\tau[\|\widehat\bx^{(\tau)} - \bx^{(\tau)}\|^2]>\epsilon^2$. Since $\tau$ is generated by Algorithm~\ref{alg:sgm}, we have \eqref{eq:stationaryopposite}. Note that the right-hand side of \eqref{eq:stationaryopposite} is well-defined because we just proved $I\neq\emptyset$. \eqref{eq:stationaryopposite} and \eqref{eq:propeq4} imply
\eqref{eq:leftright}. Combining \eqref{eq:leftright} and \eqref{eq:propeq3} leads to the opposite of \eqref{eq:propeq2}. This contradiction means $\mathbb{E}_\tau[\|\widehat\bx^{(\tau)} - \bx^{(\tau)}\|^2]\leq\epsilon^2$.

Given the result above, we only need to show that the choices of $S$, $T$, $\eta_t$ and $\epsilon_t$ ensure \eqref{eq:propeq4} and \eqref{eq:propeq3}.

In Case I, \eqref{eq:propeq4} holds because of Lemma~\ref{thm:boundlambda} and the choice of $\epsilon_t$. Let $\eta=\eta_t=\frac{2\epsilon^2(\hat\rho-\rho)}{5(1+\Lambda)M^2}$ for any $t$. Using Lemma~\ref{thm:boundlambda} and plugging the values of $S$, $T$, $\eta_t$ and $\epsilon_t$ in  \eqref{eq:propeq3}, we can show that \eqref{eq:propeq3} is equivalent to
\begin{align*}
\frac{T\eta\hat\rho\epsilon^2(\hat\rho-\rho)}{1+\Lambda}>&~ \frac{\hat\rho D^2}{2}
+\frac{\hat\rho}{2}T\eta^2M^2+\frac{\hat\rho}{2}\max\left\{\sqrt{12\ln(8/\delta)},\frac{4}{3}\ln(8/\delta)\right\}\sqrt{T}\eta^2M^2\\\nonumber
&+2\sqrt{3\ln(4/\delta)}\hat{\rho}\sqrt{T}\eta MD,
\end{align*}
which, after dividing both sides by $T\eta\hat\rho$, can be equivalently written as
\begin{align*}
\frac{\epsilon^2(\hat\rho-\rho)}{1+\Lambda}>&~ \frac{D^2}{2T\eta}
+\frac{\eta M^2}{2}+\frac{1}{2}\max\left\{\sqrt{12\ln(8/\delta)},\frac{4}{3}\ln(8/\delta)\right\}\frac{\eta M^2}{\sqrt{T}}\\\nonumber
&+\frac{2\sqrt{3\ln(4/\delta)}MD}{\sqrt{T}}.
\end{align*}
By the values of $\eta$ and $T$, each summand in the right-hand side of the inequality above is no more than $\frac{\epsilon^2(\hat\rho-\rho)}{5(1+\Lambda)}$ so the right-hand side of the inequality above no more than $\frac{4\epsilon^2(\hat\rho-\rho)}{5(1+\Lambda)}$ which is strictly less than the left-hand side. This means \eqref{eq:propeq3} holds with this choice of parameters and thus $\mathbb{E}_\tau[\|\widehat\bx^{(\tau)} - \bx^{(\tau)}\|]\leq\epsilon$. We can prove $\mathbb{E}_\tau [g(\bx^{(\tau)})]\leq \frac{\epsilon^2(\hat\rho-\rho)}{1+\lambda}$ in the same way as in  the proof of Theorem~\ref{thm:mdconverge_deterministic}.

In Case II, by the choices of $\epsilon_t$ and $T$, we have, for any $t\in\{S,S+1,\dots,T-1\}$,  
\begin{align}
\label{eq:epsilont_new}
\epsilon_t =\frac{EMD}{\sqrt{t+1}}\leq \frac{EMD}{\sqrt{S+1}}=\frac{EMD}{\sqrt{T/2+1}}\leq \frac{\epsilon^2(\hat\rho-\rho)}{1+\Lambda}.
\end{align}
This further implies \eqref{eq:propeq4} because of Lemma~\ref{thm:boundlambda}.  Note that $\eta_t$ and $\epsilon_t$ are decreasing in $t$. Hence, the left-hand side of  \eqref{eq:propeq3} satisfies
\begin{eqnarray}
\label{eq:propeq5_new}
\sum_{t=S}^{T-1}\eta_{t}\hat{\rho}\epsilon_t>\frac{T}{2}\eta_T\hat{\rho}\epsilon_T=\frac{ET}{2T+2}\hat{\rho}D^2\geq \frac{E\hat{\rho}D^2}{4}.
\end{eqnarray}

The right-hand side of \eqref{eq:propeq3} satisfies 
\small
\begin{align}
\nonumber
&~\frac{\hat\rho D^2}{2}
+\frac{\hat\rho}{2}\sum_{t=S}^{T-1}\eta_t^2M^2+\frac{\hat\rho}{2}\max\left\{\sqrt{12\ln(8/\delta)},\frac{4}{3}\ln(8/\delta)\right\}\sqrt{\sum_{t=S}^{T-1}\eta_t^4M^4}\\\nonumber
&+\sqrt{3\ln(4/\delta)}\sqrt{\sum_{t=S}^{T-1}4\eta_{t}^2\hat{\rho}^2M^2D^2}\\\nonumber
=&~\frac{\hat\rho D^2}{2}
+\frac{\hat\rho}{2}D^2\sum_{t=S}^{T-1}\frac{1}{t+1}+\frac{\hat\rho}{2}\max\left\{\sqrt{12\ln(8/\delta)},\frac{4}{3}\ln(8/\delta)\right\}D^2\sqrt{\sum_{t=S}^{T-1}\frac{1}{(t+1)^2}}\\\nonumber
&+2\sqrt{3\ln(4/\delta)}\hat{\rho}D^2\sqrt{\sum_{t=S}^{T-1}\frac{1}{t+1}}\\\label{eq:propeq6}
\leq&~\hat\rho D^2+\frac{\hat\rho\pi}{2\sqrt{6}}\max\left\{\sqrt{12\ln(8/\delta)},\frac{4}{3}\ln(8/\delta)\right\}D^2+2\sqrt{3\ln(4/\delta)}\hat{\rho}D^2,
\end{align}
\normalsize
where the equality is obtained by plugging in the definition of $\eta_t$ and the inequality is because $\sum_{t=S}^{T-1}\frac{1}{t+1}\leq \int_{S}^{T}\frac{1}{t}dt=\ln(T/S)=\ln(2)\leq 1$ and $\sum_{t=S}^{T-1}\frac{1}{(t+1)^2}\leq\pi^2/6$.  By the condition satisfied by $E$, the right-hand side of \eqref{eq:propeq5_new} is greater than or equal to the right-hand side \eqref{eq:propeq6}. This means \eqref{eq:propeq9} holds with this choice of parameters and thus $\mathbb{E}_\tau[\|\widehat\bx^{(\tau)} - \bx^{(\tau)}\|]\leq\epsilon$. We can prove $\mathbb{E}_\tau [g(\bx^{(\tau)})]\leq \frac{\epsilon^2(\hat\rho-\rho)}{1+\lambda}$ in the same way as in the proof of  Theorem~\ref{thm:mdconverge_deterministic}.
\end{proof}

\subsection{Complexity analysis when subgradient and function oracles are both stochastic}
\label{sec:mdconverge_prob}

In this section, we will analyze the oracle complexity of Algorithm~\ref{alg:sgm} under Assumption~\ref{assume:stochasticgrad_prob} without requiring \eqref{eq:reduced_to_deterministic} or \eqref{eq:reduced_to_semi_deterministic}. In this case, not only the subgradient oracles are stochastic but also the function value $g$ can only be accessed through a stochastic oracle. When the function value of $g$ is stochastic, the complexity analysis becomes fundamentally more challenging  than the case of a deterministic $g$. In fact, the challenge comes only from the stochastic function value of $g$ instead of its stochastic subgradient. Note that, in the fully deterministic case, Algorithm~\ref{alg:dsgm} essentially updates $\bx^{(t)}$ along a hybrid subgradient 
\begin{eqnarray}
\label{eq:hybridsubgradient}
\mathbb{I}(g(\bx^{(t)})\leq \epsilon_t)\bzt_f^{(t)}+\mathbb{I}(g(\bx^{(t)})> \epsilon_t)\bzt_g^{(t)}.
\end{eqnarray}
If only the subgradients are stochastic but the function value $g(\bx^{(t)})$ remains deterministic, the hybrid stochastic subgradient
$\mathbb{I}(g(\bx^{(t)})\leq \epsilon_t)\bxi_f^{(t)}+\mathbb{I}(g(\bx^{(t)})> \epsilon_t)\bxi_g^{(t)}$
provides an unbiased estimation of \eqref{eq:hybridsubgradient}. In this case, we can still obtain complexity of $O(1/\epsilon^4)$ (see Theorem~\ref{thm:mdconverge}) by a proof similar to the deterministic case. However, when $g(\bx^{(t)})$ must be queried through some unbiased estimator $w^{(t)}$, the naively constructed direction 
\begin{eqnarray*}
\mathbb{I}(w^{(t)}\leq \epsilon_t)\bxi_f^{(t)}+\mathbb{I}(w^{(t)}> \epsilon_t)\bxi_g^{(t)}
\end{eqnarray*}
is not an unbiased estimator of \eqref{eq:hybridsubgradient}. To tackle this issue, we have to query a mini-batch of $w^{(t)}$ of size $B$, i.e., $\{\omega_i^{(t)}\}_{i=1}^B$ to construct $\bar w^{(t)}$ as an estimation of $g(\bx^{(t)})$ with a high accuracy in a high probability. In this way, we can use
\begin{eqnarray*}
\mathbb{I}(\bar w^{(t)}\leq \epsilon_t)\bxi_f^{(t)}+\mathbb{I}(\bar w^{(t)}> \epsilon_t)\bxi_g^{(t)}
\end{eqnarray*}
as a nearly unbiased estimator of \eqref{eq:hybridsubgradient}, which leads to the switching condition in Algorithm~\ref{alg:sgm}. As a consequence, the oracle complexity of the function value of $g$ increases from $O(1/\epsilon^4)$ to $\tilde O(1/\epsilon^8)$.~\footnote{$\tilde O(\cdot)$ suppresses logarithmic factors in the order of complexity. } although the oracle complexity of stochastic subgradients remains $O(1/\epsilon^4)$. 

 \begin{theorem}
\label{thm:mdconverge_prob}
Suppose Assumptions~\ref{assume:allpaper}, \ref{assume:convex} and~\ref{assume:stochasticgrad_prob} hold and $\Lambda$ is as in \eqref{eq:Lambdabound}. Let $\widehat\bx(\bx^{(t)})$ be defined as in \eqref{eq:phix} with $(\hat\rho,\tilde\rho)$ satisfying \eqref{eq:parameter1} and $\bx^{(\tau)}$ is generated by Algorithm~\ref{alg:sgm}. Algorithm~\ref{alg:sgm} guarantees 
$\mathbb{E}_\tau[\|\widehat\bx^{(\tau)} - \bx^{(\tau)}\|]\leq\epsilon$  and $\mathbb{E}_\tau [g(\bx^{(\tau)})]\leq \frac{2\epsilon^2(\hat\rho-\rho)}{1+\Lambda}$ with probability at least $1-\delta$ in either of the following cases. 

Case I: If $S$, $\epsilon_t$, $\eta_t$ and $T$ are chosen as Case I in Theorem~\ref{thm:mdconverge} and 
$B=\frac{300\sigma^2\ln(4T/\delta)(1+\Lambda)^4}{\epsilon^4(\hat\rho-\rho)^2}$.

Case II: If $S$, $\epsilon_t$, $\eta_t$ and $T$ are chosen as Case II in Theorem~\ref{thm:mdconverge} except that $E$ is any positive constant such that 
\small
\begin{align*}
E\geq &~8+\frac{2\pi}{\sqrt{6}}\max\Big\{\sqrt{12\ln(8/\delta)},\frac{4}{3}\ln(8/\delta)\Big\}+8\sqrt{3\ln(4/\delta)}
\end{align*}
\normalsize
and 
$B=\frac{3T\sigma^2\ln(2T/\delta)(1+\Lambda)^2}{2M^2D^2}$.
\end{theorem}
In each iteration of Algorithm~\ref{alg:sgm}, we query one stochastic subgradient of $f$ or $g$ but $B$ stochastic function values of $g$. In both Case I and Case II, we have $T=O(1/\epsilon^4)$ and $B=\tilde O(1/\epsilon^4)$ so the subgradient oracle complexity is still $O(1/\epsilon^4)$ but the function value oracle complexity becomes $\tilde O(1/\epsilon^8)$, which is higher than the $O(1/\epsilon^6)$ complexity by the double-loop methods in \cite{boob2022stochastic, ma2020quadratically}. It is our future work to reduce the complexity when $g$ is stochastic, for example, by a single-loop primal-dual method that uses a hybrid subgradient like $\bxi_f^{(t)}+\lambda_t\bxi_g^{(t)}$ with the dual variable $\lambda_t$ updated by only one sample of $w^{(t)}$.

 

\begin{proof}[Proof of Theorem~\ref{thm:mdconverge_prob}]
By Proposition~\ref{thm:mainprop_prob}, with a probability of at least $1-\delta$, we simultaneously have 
\eqref{eq:propeq1} with $\mu=0$ and \eqref{eq:highprobgx4} for $t=S,S+1,\dots,T-1$. In the rest of the proof, we assume \eqref{eq:propeq1} with $\mu=0$ and \eqref{eq:highprobgx4} hold for $t=S,S+1,\dots,T-1$.

We first prove that, if $S$, $T$, $B$, $\eta_t$ and $\epsilon_t$ are chosen such that \eqref{eq:propeq4} holds and
\begin{align}
    \label{eq:propeq7}
    \sum_{t=S}^{T-1}\eta_{t}\hat{\rho}\epsilon_t>&~\frac{\hat\rho D^2}{2}
    +\frac{\hat\rho}{2}\sum_{t=S}^{T-1}\eta_t^2M^2+\frac{\hat\rho}{2}\max\left\{\sqrt{12\ln(8/\delta)},\frac{4}{3}\ln(8/\delta)\right\}\sqrt{\sum_{t=S}^{T-1}\eta_t^4M^4}\\\nonumber
    &+\sqrt{3\ln(4/\delta)}\sqrt{\sum_{t=S}^{T-1}4\eta_{t}^2\hat{\rho}^2M^2D^2}+\sum_{t=S}^{T-1}\eta_{t}\hat{\rho}\sqrt{\frac{3}{B}}\sigma\sqrt{\ln(4(T-S)/\delta)}(1+\Lambda),
\end{align}
we must have $\bar\omega^{(t)}\leq\epsilon_t$ for at least one $t$ in $\{S,S+1,\dots,T-1\}$ (i.e., $I\neq\emptyset$) and $\mathbb{E}_\tau[\|\widehat\bx^{(\tau)} - \bx^{(\tau)}\|^2]\leq\epsilon^2$ (so $\mathbb{E}_\tau[\|\widehat\bx^{(\tau)} - \bx^{(\tau)}\|]\leq\epsilon$). 

Suppose \eqref{eq:propeq7} holds and $\bar\omega^{(t)}>\epsilon_t$ for $t=S,S+1,\dots,T-1$, i.e., $I=\emptyset$. \eqref{eq:propeq1} contradicts with \eqref{eq:propeq7} as $\Lambda\geq \widehat\lambda_t$. This contradiction means  $\bar\omega^{(t)}\leq\epsilon_t$ for at least one $t$ in $\{S,S+1,\dots,T-1\}$ so $I\neq\emptyset$. Suppose \eqref{eq:propeq4} and \eqref{eq:propeq7} hold but $\mathbb{E}_\tau[\|\widehat\bx^{(\tau)} - \bx^{(\tau)}\|^2]>\epsilon^2$. Since $\tau$ is generated by Algorithm~\ref{alg:sgm}, we have 
\begin{align}
\label{eq:stationaryopposite_new}
\epsilon^2<\mathbb{E}_\tau[\|\widehat\bx^{(\tau)} - \bx^{(\tau)}\|^2]
=\frac{\sum_{t=S}^{T-1}\eta_{t}\mathbb{I}(\bar\omega^{(t)}\leq\epsilon_t)\|\widehat\bx^{(t)}-\bx^{(t)} \|^{2}}{\sum_{t=S}^{T-1}\eta_{t}\mathbb{I}(\bar\omega^{(t)}\leq\epsilon_t)}.
\end{align}
Note that the right-hand side of \eqref{eq:stationaryopposite_new} is well-defined because we just proved $I\neq\emptyset$.
\eqref{eq:stationaryopposite_new} and \eqref{eq:propeq4} imply
\begin{align}
    \nonumber
    &\sum_{t=S}^{T-1}\left[\eta_{t}\hat{\rho}(\hat{\rho}-\rho)\|\widehat\bx^{(t)}-\bx^{(t)} \|^{2}-\eta_{t}\hat{\rho}\widehat\lambda_{t}\epsilon_t\right]\mathbb{I}(\bar\omega^{(t)}\leq \epsilon_t)+\sum_{t=S}^{T-1}\eta_{t}\hat{\rho}\epsilon_t\mathbb{I}(\bar\omega^{(t)}> \epsilon_t)\\\nonumber
    >&~\sum_{t=S}^{T-1}\left[\eta_{t}\hat{\rho}(\hat{\rho}-\rho)\epsilon^{2}-\eta_{t}\hat{\rho}\widehat\lambda_{t}\epsilon_t\right]\mathbb{I}(\bar\omega^{(t)}\leq \epsilon_t)+\sum_{t=S}^{T-1}\eta_{t}\hat{\rho}\epsilon_t\mathbb{I}(\bar\omega^{(t)}> \epsilon_t)\\\label{eq:leftright_new}
    \geq&~\sum_{t=S}^{T-1}\eta_{t}\hat{\rho}\epsilon_t\mathbb{I}(\bar\omega^{(t)}\leq \epsilon_t)+\sum_{t=S}^{T-1}\eta_{t}\hat{\rho}\epsilon_t\mathbb{I}(\bar\omega^{(t)}> \epsilon_t)
\geq\sum_{t=S}^{T-1}\eta_{t}\hat{\rho}\epsilon_t,
\end{align}
where the second inequality is  because of \eqref{eq:propeq4}. Combining \eqref{eq:leftright_new} and \eqref{eq:propeq7} leads to the opposite of \eqref{eq:propeq1}. This contradiction means $\mathbb{E}_\tau[\|\widehat\bx^{(\tau)} - \bx^{(\tau)}\|^2]\leq\epsilon^2$.

Given the result above, we only need to show that the choices of $S$, $T$, $B$, $\eta_t$ and $\epsilon_t$ ensure \eqref{eq:propeq4} and \eqref{eq:propeq7}.

In Case I, \eqref{eq:propeq4} holds because of Lemma~\ref{thm:boundlambda} and the choice of $\epsilon_t$. Let $\eta=\eta_t=\frac{2\epsilon^2(\hat\rho-\rho)}{5(1+\Lambda)M^2}$ for any $t$. Using Lemma~\ref{thm:boundlambda} and plugging the values of $S$, $T$, $B$, $\eta_t$ and $\epsilon_t$ in  \eqref{eq:propeq7}, we can show that \eqref{eq:propeq7} holds if
\begin{align*}
\frac{T\eta\hat\rho\epsilon^2(\hat\rho-\rho)}{1+\Lambda}>&~\frac{\hat\rho D^2}{2}+\frac{\hat\rho}{2}T\eta^2M^2+\frac{\hat\rho}{2}\max\left\{\sqrt{12\ln(8/\delta)},\frac{4}{3}\ln(8/\delta)\right\}\sqrt{T}\eta^2M^2\\\nonumber
&+2\sqrt{3\ln(4/\delta)}\hat{\rho}\sqrt{T}\eta MD+T\eta\hat{\rho}\sqrt{\frac{3}{B}}\sigma\sqrt{\ln(4T/\delta)}(1+\Lambda),
\end{align*}
which, after dividing both sides by $T\eta\hat\rho$, can be equivalently written as
\begin{align*}
\frac{\epsilon^2(\hat\rho-\rho)}{1+\Lambda}>&~\frac{D^2}{2T\eta}
+\frac{\eta M^2}{2}+\frac{1}{2}\max\left\{\sqrt{12\ln(8/\delta)},\frac{4}{3}\ln(8/\delta)\right\}\frac{\eta M^2}{\sqrt{T}}\\\nonumber
&+\frac{2\sqrt{3\ln(4/\delta)}MD}{\sqrt{T}}+\sqrt{\frac{3}{B}}\sigma\sqrt{\ln(4T/\delta)}(1+\Lambda).
\end{align*}
By the values of $\eta$, $B$ and $T$, each of the first four summands on the right-hand side of the inequality above is no more than $\frac{\epsilon^2(\hat\rho-\rho)}{5(1+\Lambda)}$ while the last summand is no more than $\frac{\epsilon^2(\hat\rho-\rho)}{10(1+\Lambda)}$, so the right-hand side of the inequality above no more than $\frac{9\epsilon^2(\hat\rho-\rho)}{10(1+\Lambda)}$ which is strictly less than the left-hand side. This means \eqref{eq:propeq9} holds with this choice of parameters and thus $\mathbb{E}_\tau[\|\widehat\bx^{(\tau)} - \bx^{(\tau)}\|]\leq\epsilon$. Moreover, by \eqref{eq:highprobgx4}, we have 
\begin{align*}
\mathbb{E}_\tau [g(\bx^{(\tau)})]=&~\frac{\sum_{t=0}^{T-1}\eta_{t}g (\bx^{(t)} )\mathbb{I}(\bar\omega^{(t)}\leq \epsilon_t)}{\sum_{t=0}^{T-1}\eta_{t}\mathbb{I}(\bar\omega^{(t)}\leq \epsilon_t)}
\leq 
\frac{\sum_{t=0}^{T-1}\eta_{t}\left(\epsilon_t+\sqrt{\frac{3}{B}}\sigma\sqrt{\ln(4T/\delta)}\right)\mathbb{I}(\bar\omega^{(t)}\leq \epsilon_t)}{\sum_{t=0}^{T-1}\eta_{t}\mathbb{I}(\bar\omega^{(t)}\leq \epsilon_t)}\\
\leq&~\frac{\epsilon^2(\hat\rho-\rho)}{1+\Lambda}+\sqrt{\frac{3}{B}}\sigma\sqrt{\ln(4T/\delta)}\leq \frac{2\epsilon^2(\hat\rho-\rho)}{1+\Lambda},
\end{align*}
where the last inequality is because of the choice of $B$.

In Case II, by the choice of $\epsilon_t$, we have \eqref{eq:epsilont_new} holds. This further implies \eqref{eq:propeq4} because of Lemma~\ref{thm:boundlambda}.  Note that $\eta_t$ and $\epsilon_t$ are decreasing in $t$. Hence, we also have \eqref{eq:propeq5_new}. The right-hand side of \eqref{eq:propeq7} satisfies 
\small
\begin{align}
\nonumber
&~\frac{\hat\rho D^2}{2}
+\frac{\hat\rho}{2}\sum_{t=S}^{T-1}\eta_t^2M^2+\frac{\hat\rho}{2}\max\left\{\sqrt{12\ln(8/\delta)},\frac{4}{3}\ln(8/\delta)\right\}\sqrt{\sum_{t=S}^{T-1}\eta_t^4M^4}\\\nonumber
&+\sqrt{3\ln(4/\delta)}\sqrt{\sum_{t=S}^{T-1}4\eta_{t}^2\hat{\rho}^2M^2D^2}+\sum_{t=S}^{T-1}\eta_{t}\hat{\rho}\sqrt{\frac{3}{B}}\sigma\sqrt{\ln(4(T-S)/\delta)}(1+\Lambda)\\\nonumber
=&~\frac{\hat\rho D^2}{2}
+\frac{\hat\rho}{2}D^2\sum_{t=S}^{T-1}\frac{1}{t+1}+\frac{\hat\rho}{2}\max\left\{\sqrt{12\ln(8/\delta)},\frac{4}{3}\ln(8/\delta)\right\}D^2\sqrt{\sum_{t=S}^{T-1}\frac{1}{(t+1)^2}}\\\nonumber
&+2\sqrt{3\ln(4/\delta)}\hat{\rho}D^2\sqrt{\sum_{t=S}^{T-1}\frac{1}{t+1}}+\frac{\hat{\rho}D}{M}\sqrt{\frac{3}{B}}\sigma\sqrt{\ln(4(T-S)/\delta)}(1+\Lambda)\sum_{t=S}^{T-1}\frac{1}{\sqrt{t+1}}\\\label{eq:propeq8}
\leq&~\hat\rho D^2+\frac{\hat\rho\pi}{2\sqrt{6}}\max\left\{\sqrt{12\ln(8/\delta)},\frac{4}{3}\ln(8/\delta)\right\}D^2+2\sqrt{3\ln(4/\delta)}\hat{\rho}D^2+\hat\rho D^2,
\end{align}
\normalsize
where the equality is obtained by plugging in the definition of $\eta_t$ and the inequality is because of the definition of $B$ and the facts that $\sum_{t=S}^{T-1}\frac{1}{t+1}\leq \int_{S}^{T}\frac{1}{t}dt=\ln(T/S)=\ln(2)\leq 1$, $\sum_{t=S}^{T-1}\frac{1}{(t+1)^2}\leq\pi^2/6$ and $\sum_{t=S}^{T-1}\frac{1}{\sqrt{t+1}}\leq \int_{S}^{T}\frac{1}{\sqrt{t}}dt\leq \sqrt{T/2}$. 
By the condition  of $E$, the right-hand side of \eqref{eq:propeq5_new} is strictly greater than the right-hand side \eqref{eq:propeq8}. This means \eqref{eq:propeq7} holds with this choice of parameters and thus $\mathbb{E}_\tau[\|\widehat\bx^{(\tau)} - \bx^{(\tau)}\|]\leq\epsilon$. Moreover, since \eqref{eq:highprobgx} holds with $\delta$ replaced by $\delta/4$ for  $t=S,S+1,\dots,T-1$, we have 
\begin{align*}
\mathbb{E}_\tau [g(\bx^{(\tau)})]=&~\frac{\sum_{t=S}^{T-1}\eta_{t}g (\bx^{(t)} )\mathbb{I}(\bar\omega^{(t)}\leq \epsilon_t)}{\sum_{t=S}^{T-1}\eta_{t}\mathbb{I}(\bar\omega^{(t)}\leq \epsilon_t)}
\leq 
\frac{\sum_{t=S}^{T-1}\eta_{t}\left(\epsilon_t+\sqrt{\frac{3}{B}}\sigma\sqrt{\ln(2T/\delta)}\right)\mathbb{I}(\bar\omega^{(t)}\leq \epsilon_t)}{\sum_{t=S}^{T-1}\eta_{t}\mathbb{I}(\bar\omega^{(t)}\leq \epsilon_t)}\\
\leq&~ \epsilon_{T/2}+\sqrt{\frac{3}{B}}\sigma\sqrt{\ln(2T/\delta)}\leq \frac{2\epsilon^2(\hat\rho-\rho)}{1+\Lambda},
\end{align*}
where the last inequality is because of \eqref{eq:epsilont_new} and the choices of $B$ and $T$.
\end{proof}

\begin{remark}
When $g$ is $\mu$-strongly convex with $\mu>0$, under the same assumptions as Theorem~\ref{thm:mdconverge},  we can  establish $O(1/\epsilon^4)$ oracle complexity for  Algorithm~\ref{alg:sgm}. Similarly, under the same assumptions as Theorem~\ref{thm:mdconverge_prob},  we can also establish $\tilde O(1/\epsilon^8)$ oracle complexity for Algorithm~\ref{alg:sgm}. These two complexity results of Algorithm~\ref{alg:sgm} can be proved in a similar way as Theorem~\ref{thm:gstrongconvex}. Since those results do not provide additional insights on the complexity and analysis for the SSG method, we do not include them in the paper. Similar to Theorem~\ref{thm:gstrongconvex}, when $\mu>0$, we only need to set $\epsilon_t=0$ in Algorithm~\ref{alg:sgm} and at last generate $\tau$ by Output II in Algorithm~\ref{alg:dsgm}. This indicates that strong convexity in $g$ also reduces the number of tuning parameters in Algorithm~\ref{alg:sgm} just as it does for Algorithm~\ref{alg:dsgm}. 
\end{remark}

\subsection{Complexity analysis for convex constraints with bounded $\mathcal{S}$ (instead of $\mathcal{X}$)}
\label{sec:convexgsharpness}

In this section, we present the complexity analysis when $f$ is weakly convex and $g$ is convex after a small modification on Assumption~\ref{assume:convex}. Recall that 
\begin{align*}
    g_+(\bx)=\max\{g(\bx),0\},\quad \mathcal{L}=\{\bx\in\X~|~g(\bx)=0\}~\text{ and }~\mathcal{S}=\{\bx\in\X~|~g(\bx)\leq 0\}.  
\end{align*}
In this section, we choose parameters in \eqref{eq:phi} the same as in \eqref{eq:parameter1}. We make the following assumptions in addition to Assumption~\ref{assume:allpaper} in this section.
\begin{assumption}
\label{assume:convex_sharpness}
The following statements hold:
\begin{itemize}
    \item[A.] $f(\bx)$ is $\rho$-weakly convex on $\X$ and $g(\bx)$ is convex on $\X$. 
    \item[B.] (Slater's condition) There exists $\bx_{\text{feas}}\in\text{relint}(\X)$ such that $g(\bx_{\text{feas}})<0$.
    \item[C.]  There exists $D>0$ such that $\|\bx-\bx'\|\leq D$ for any $\bx$ and $\bx'$ in $\mathcal{S}$.
\end{itemize}
\end{assumption}

Under Assumption~\ref{assume:convex_sharpness},  $\X$ is not necessarily bounded (e.g. $\X=\mathbb{R}^d$) while $g$ must have a bounded $0$-sublevel set over $\X$. 

We next show that Assumptions~\ref{assume:convex_sharpness} B and C imply that subgradient of $g+\delta_\X$ can be bounded away from zero on $\mathcal{L}$. Here, $\delta_\X$ is the indicator function of $\X$. Given any $\bx\in\mathcal{L}$, by Assumption~\ref{assume:convex_sharpness}B and the convexity of $g+\delta_\X$, we have
\begin{align*}
g(\bx_{\text{feas}})\geq g(\bx)+\langle\bzt_g+\bu,\bx_{\text{feas}}-\bx\rangle\geq-\|\bzt_g+\bu\|\cdot\|\bx_{\text{feas}}-\bx\|
\end{align*}
for any $\bzt_g\in\partial g(\bx)$ and $\bu\in\mathcal{N}_{\X}(\bx)$. Since $\|\bx_{\text{feas}}-\bx\|\leq D$ by Assumption~\ref{assume:convex_sharpness}C, it holds that
\begin{equation}
\label{eq:sharpsubgradient_c}
\min\limits_{\bzt_g\in\partial g(\bx),\bu\in\mathcal{N}_{\X}(\bx)}\|\bzt_g+\bu\|\geq\frac{-g(\bx_{\text{feas}})}{\|\bx_{\text{feas}}-\bx\|}\geq\nu':=-g(\bx_{\text{feas}})/D>0.
\end{equation}
This result is similar to \eqref{eq:sharpsubgradient}.
 
Note that $\mathcal{S}$ is the optimal set of $\min_{\bx\in\X}g_+(\bx)$, which is a convex non-smooth optimization problem with an optimal value of zero. By Lemma 1 in \cite{yang2018rsg}, \eqref{eq:sharpsubgradient_c} implies that $g_+(\bx)$ satisfies a global linear error bound on $\mathcal{X}$, namely, for any $\bx\in\X$, 
\begin{eqnarray}
\label{eq:errorbound_c}
\nu'\textup{dist}(\bx,\mathcal{S})\leq g_+(\bx).
\end{eqnarray}
This result is stronger than \eqref{eq:errorbound_wc} because \eqref{eq:errorbound_wc} does not hold globally over $\X$. 
Moreover, this result implies $\nu'\leq M$.

Following almost the same proof as Lemma~\ref{thm:nonexpansion}, we obtain a result similar to \eqref{eq:shrinkdist_Q}. 
\begin{lemma}
\label{thm:nonexpansion_c}
Suppose Assumptions~\ref{assume:allpaper} and \ref{assume:convex_sharpness} hold. Also, suppose the sequence 
$\{\bx^{(t)}\}_{t\geq0}$ is generated by applying the projected subgradient method to $\min_{\bx\in\X}g_+(\bx)$ using a Polyak's stepsize, namely, 
\begin{equation}
\label{eq:Polyak_c}
\bx^{(t+1)}=\textup{proj}_{\X}(\bx^{(t)}-\eta_t\bzt_g^{(t)}), \hspace{0.5em} \eta_t=
\left\{
\begin{array}{ll}
v'g_+(\bx^{(t)})/(2M\|\bzt_g^{(t)}\|^2) & \text{ if }~\bzt_g^{(t)}\neq \mathbf{0}\\
0& \text{ if }~\bzt_g^{(t)}= \mathbf{0}
\end{array}
\right.,
\hspace{0.5em} \text{ for }t=0, 1, \dots,
\end{equation}
where $\bzt_g^{(t)}\in\partial g_+(\bx^{(t)})$.  We have 
\begin{eqnarray}
\label{eq:shrinkdist_Q_c}
\textup{dist}^2(\bx^{(t+1)},\mathcal{S})\leq\left(1-\frac{3\nu'^3}{4M^3}\right)\textup{dist}^2(\bx^{(t)},\mathcal{S}),\quad\forall t\geq0.
\end{eqnarray}
\end{lemma}
\begin{proof}
Suppose $\bzt_g^{(t)}=\mathbf{0}$ so $\bx^{(t+1)}=\bx^{(t)}$. This happens only when $g_+(\bx^{(t)})=0$. By \eqref{eq:errorbound_c}, we must have $\textup{dist}(\bx^{(t)},\mathcal{S})=0$ and thus $\textup{dist}(\bx^{(t+1)},\mathcal{S})=0\leq\left(1-\frac{3\nu'^3}{4M^3}\right)\text{dist}(\bx^{(t)},\mathcal{S})$ and \eqref{eq:shrinkdist_Q_c} holds. 

Suppose $\bzt_g^{(t)}\neq\mathbf{0}$.
Let $\bx^{\dagger(t)}=\text{proj}_{\mathcal{S}}(\bx^{(t)})$ for $t\geq0$. Following a proof similar  to Lemma~\ref{thm:nonexpansion}, which is originally from \cite{davis2018subgradient}, we have
\small
\begin{align*}
\nonumber
   \text{dist}^2(\bx^{(t+1)},\mathcal{S})
    \leq&~ \text{dist}^2(\bx^{(t)},\mathcal{S})-2\eta_{t}\left\langle\bzt_g^{(t)},\bx^{(t)}-\bx^{\dagger(t)}\right\rangle+\eta_{t}^{2} \|\bzt_g^{(t)}\|^{2}\\\nonumber
   \leq&~ \text{dist}^2(\bx^{(t)},\mathcal{S})+\frac{\nu'g_+(\bx^{(t)})}{M\|\bzt_g^{(t)}\|^2}\left(g_+(\bx^{\dagger(t)})-g_+(\bx^{(t)})\right)+\frac{\nu'^2g_+^2(\bx^{(t)})}{4M^2\|\bzt_g^{(t)}\|^2}\\\nonumber
    =&~ \text{dist}^2(\bx^{(t)},\mathcal{S})-\left(\frac{\nu'}{M}-\frac{\nu'^2}{4M^2}\right)\left(g_+(\bx^{(t)})\right)^2/\|\bzt_g^{(t)}\|^2\\\nonumber
    \leq&~ \text{dist}^2(\bx^{(t)},\mathcal{S})-\frac{3\nu'^3}{4M^3}\|\bx^{(t)}-\bx^{\dagger(t)}\|^2,
\end{align*}
\normalsize
where the equality is because $g_+(\bx^{\dagger(t)})=0$, the second inequality is by the convexity of $g_+$ and the last is by \eqref{eq:errorbound_c} and the fact that $\nu'\leq M$. 
\end{proof}

Using almost the same proof as Proposition~\ref{thm:feasiblesubproblem}, we can show the following proposition.
\begin{proposition}
\label{thm:feasiblesubproblem_c}
Suppose Assumptions~\ref{assume:allpaper} and \ref{assume:convex_sharpness} hold and $\epsilon>0$. Also, suppose that $\bx^{(t)}$ is generated by Algorithm~\ref{alg:dsgm} using 
$\bx^{(0)}\in\mathcal{S}$, $\epsilon_{t}=\frac{\nu'}{2}\min\left\{\epsilon^2/M,\nu'\right\}$ and 
\begin{align*}
\eta_{t}=
\left\{
\begin{array}{ll}
\frac{\nu'}{2M^2}\min\left\{\epsilon^2/M,\nu'\right\}&\text{ if }t\in I\\
\nu'g(\bx^{(t)})/(2M\|\bzt_g^{(t)}\|^2)&\text{ if }t\in J.
\end{array}
\right.
\end{align*}
Then $\textup{dist}(\bx^{(t)},\mathcal{S})\leq\min\left\{\epsilon^2/M,\nu'\right\}$ and $g(\bx^{(t)})\leq \epsilon^2$ for any $t\geq0$. As a consequence, $\bx^{(t)}$ is $\epsilon^2$-feasible to \eqref{eq:phi} where $\bx=\bx^{(t)}$ and $(\hat\rho,\tilde\rho)$ satisfies \eqref{eq:parameter1}. 
\end{proposition}
\begin{proof}
The choices of $\epsilon_t$ and $\eta_t$ imply 
\begin{equation}
\label{eq:stepetainI}
\epsilon_t/\nu'+\eta_tM^2/\nu'\leq\min\left\{\epsilon^2/M,\nu'\right\}
\end{equation}
for $t\in I$. We prove $\text{dist}(\bx^{(t)},\mathcal{S})\leq \min\left\{\epsilon^2/M,\nu'\right\}$ by induction on $t$. Since $\bx^{(0)}\in \mathcal{S}$, this conclusion holds trivially for $t=0$. Suppose it holds up to iteration $t$. We want to prove it also holds for iteration $t+1$. 

Suppose $t\in I$. Since $g(\bx^{(t)})\leq \epsilon_t$ for $t\in I$, we have
\begin{align*}
    g(\bx^{(t+1)})\leq g(\bx^{(t)})+M\|\bx^{(t+1)}-\bx^{(t)}\|
    \leq g(\bx^{(t)})+M\|\eta_t\bzt_f^{(t)}\|
    \leq\epsilon_t+\eta_tM^2.
\end{align*}
By  \eqref{eq:errorbound_c}, we have $\text{dist}(\bx^{(t+1)},\mathcal{S})\leq \epsilon_t/\nu'+\eta_tM^2/\nu'\leq \min\left\{\epsilon^2/M,\nu'\right\}$.  

Suppose $t\in J$. Let $t'$ be the largest index in $I$ that is smaller than $t$. By the same proof as in the previous case, we have $\text{dist}(\bx^{(t'+1)},\mathcal{S})\leq\min\left\{\epsilon^2/M,\nu'\right\}$.  Since indexes $t'+1$, $t'+2$, ..., $t$ are in $J$, Algorithm~\ref{alg:dsgm} essentially performs the projected subgradient method to $\min_{\bx\in\X}g_+(\bx)$ using the Polyak's stepsize \eqref{eq:Polyak_c} during iterations $t'+1$, $t'+2$, ..., and $t$.  Hence, by Lemma~\ref{thm:nonexpansion_c}, we have
$\text{dist}(\bx^{(t+1)},\mathcal{S})\leq\text{dist}(\bx^{(t'+1)},\mathcal{S}) \leq \min\left\{\epsilon^2/M,\nu'\right\}$. 

By induction, we have prove that $\text{dist}(\bx^{(t)},\mathcal{S})\leq  \min\left\{\epsilon^2/M,\nu'\right\}$ for any $t\geq0$. As a result, $g(\bx^{(t)})\leq M\cdot\text{dist}(\bx^{(t)},\mathcal{S})\leq M\cdot\min\left\{\epsilon^2/M,\nu'\right\}\leq \epsilon^2$, meaning that $\bx^{(t)}$ is $\epsilon^2$-feasible for \eqref{eq:phi}.  
\end{proof}

Similar to Lemma~\ref{thm:boundlambda}, we present an upper bound of $\widehat\lambda$ in \eqref{eq:KKTprox} that is independent of $\bx$. 
\begin{lemma}
\label{thm:boundlambda_c}
Suppose Assumptions~\ref{assume:allpaper} and \ref{assume:convex_sharpness} hold. Given any $\bx\in\X$ satisfying $\textup{dist}(\bx,\mathcal{S})\leq\nu'$, let $\widehat\bx(\bx)$ be defined as in \eqref{eq:phix} with $(\hat\rho,\tilde\rho)$ satisfying \eqref{eq:parameter1} and $\widehat\lambda$ be the associated Lagrangian multiplier satisfying \eqref{eq:KKTprox}. We have
\begin{eqnarray}
\label{eq:Lambdabound_c}
\|\widehat\bx -\bx\|\leq D+\nu'\quad\text{ and }\quad\widehat\lambda \leq  \Lambda'':=\frac{\left(M+\hat\rho(D+\nu')\right)D}{-g(\bx_{\textup{feas}})}.
\end{eqnarray}
\end{lemma}

\begin{proof}
Since $\textup{dist}(\bx,\mathcal{S})\leq\nu'$, we get by triangle inequality and Assumption~\ref{assume:convex_sharpness}C that
\begin{align*}
    \|\widehat\bx -\bx\|\leq \|\widehat\bx -\bx^\dagger\|+\|\bx^\dagger -\bx\| =\|\widehat\bx -\bx^\dagger\|+\text{dist}(\bx,\mathcal{S})\leq D+\nu'
\end{align*}
where $\bx^\dagger=\text{proj}_{\mathcal{S}}(\bx)$, which gives the first result in ~\eqref{eq:Lambdabound_c}. Following the same steps as in the proof of Lemma~\ref{thm:boundlambda}, we obtain
\begin{equation}
    \label{eq:upperboundlambdaeq2}
    \widehat\lambda \leq \left(\langle\widehat\bzt_f  +\hat{\rho}(\widehat\bx - \bx), \bx_{\text{feas}}-\widehat\bx \rangle\right)/\left(-g(\bx_{\text{feas}})\right).
\end{equation}
By Assumption~\ref{assume:allpaper} and \ref{assume:convex_sharpness}C, we have $\|\widehat\bzt_f\|\leq M$ and $\|\widehat\bx -\bx_{\text{feas}}\|\leq D$, which imply from the first result in ~\eqref{eq:Lambdabound_c} and \eqref{eq:upperboundlambdaeq2} that
\small
\begin{align*}
	\widehat\lambda \leq \frac{\langle\widehat\bzt_f  +\hat{\rho}(\widehat\bx - \bx), \bx_{\text{feas}}-\widehat\bx \rangle}{-g(\bx_{\text{feas}})} 
	\leq  \frac{\left(M+\hat\rho(D+\nu')\right)D}{-g(\bx_{\text{feas}})}.
	\end{align*}
\normalsize
\end{proof}

Let $\widehat\bx^{(t)}$ be $\widehat\bx(\bx^{(t)})$ defined in \eqref{eq:phix} with $(\hat\rho,\tilde\rho)$ satisfying \eqref{eq:parameter1}. The following theorem provides the main result in this section.

\begin{theorem}
\label{thm:mdconverge_deterministic_convex_sharpness}
Under the same assumptions as Proposition~\ref{thm:feasiblesubproblem_c}, Algorithm~\ref{alg:dsgm} guarantees 
$\mathbb{E}_\tau[\|\widehat\bx^{(\tau)} - \bx^{(\tau)}\|]\leq C'\epsilon$ and $\mathbb{E}_\tau [g(\bx^{(\tau)})]\leq\frac{\nu'}{2}\min\left\{\epsilon^2/M,\nu'\right\}$, where $C':= \sqrt{\frac{\nu'(1+\Lambda'')}{2M(\hat{\rho}-\rho)}}$, if we use Output I and set 
\small
\begin{align*}
    S=0\text{ and }T>\frac{8M^3D^2}{\nu'^3\min\left\{\epsilon^4/M^2,\nu'^2\right\}}=O(1/\epsilon^4).
\end{align*}
\normalsize
\end{theorem}
This theorem indicates that Algorithm~\ref{alg:dsgm} finds a nearly $(C'\epsilon)$-stationary point with complexity $O(1/\epsilon^4)$. To obtain a nearly $\epsilon$-stationary point, one only needs to replace $\epsilon$ in $\eta_t$, $\epsilon_t$ and $T$ in this theorem above by $\epsilon/\max\{C', 1\}$ . This will only change the constant factor in the $O(1/\epsilon^4)$ complexity. Once again, this complexity matches the start-of-the-art complexity by~\cite{boob2022stochastic, ma2020quadratically, jia2022first}. 

\begin{proof}[\textbf{Proof of Theorem~\ref{thm:mdconverge_deterministic_convex_sharpness}}]
Let $\bzt^{(t)}=\bzt_f^{(t)}\in\partial f(\bx^{(t)})$ if $t\in I$ and  $\bzt^{(t)}=\bzt_g^{(t)}\in\partial g(\bx^{(t)})$ if $t\in J$. Let $\varphi(\bx)$ and $\widehat{\bx}$ be defined in \eqref{eq:phi} and \eqref{eq:phix} with $(\hat\rho,\tilde\rho)$ satisfying \eqref{eq:parameter1}. By Assumption~\ref{assume:convex_sharpness}B, problem \eqref{eq:phi} with $\bx=\bx^{(t)}$ is strongly convex and has a strictly feasible solution, so $\varphi(\bx^{(t)})$ and $\widehat\bx^{(t)}$ are well defined for any $t\geq0$. For simplicity of notation, we denote $\widehat\bx(\bx^{(t)})$ by $\widehat\bx^{(t)}$ and let $\widehat\lambda_t\geq0$ be Lagrangian multiplier satisfying \eqref{eq:KKTprox}.

By following the proof of Proposition~\ref{thm:mainprop} using deterministic subgradients and setting $\mu=0$, we have
\small
\begin{align}
\nonumber
&~\sum_{t=S}^{T-1}\left[\eta_{t}\hat{\rho}(\hat{\rho}-\rho)\|\widehat\bx^{(t)}-\bx^{(t)} \|^{2}-\eta_{t}\hat{\rho}\widehat\lambda_{t}g (\bx^{(t)} )\right]\mathbb{I}(g (\bx^{(t)} )\leq \epsilon_t)+\sum_{t=S}^{T-1}\eta_{t}\hat{\rho}g (\bx^{(t)} )\mathbb{I}(g (\bx^{(t)} )> \epsilon_t)\\\label{eq:mainthmeq6_c_old}
\leq&~\varphi (\bx^{(S)})-\varphi (\bx^{(T)})
+\dfrac{\hat{\rho}}{2}\sum_{t=S}^{T-1}\left[\eta_{t}^{2}\|\bzt_f^{(t)}\|^{2}\mathbb{I}(g (\bx^{(t)} )\leq \epsilon_t)+\eta_{t}^{2}\|\bzt_g^{(t)}\|^{2}\mathbb{I}(g (\bx^{(t)} )> \epsilon_t)\right].
\end{align}
\normalsize
Since $S=0$ and $\bx^{(0)}\in\mathcal{S}$, we then have $\varphi (\bx^{(T)})=f(\widehat\bx^{(T)})+\frac{\hat\rho}{2}\|\widehat\bx^{(T)}-\bx^{(T)}\|^2\geq f(\widehat\bx^{(T)})$ and $\varphi (\bx^{(S)})= f(\widehat\bx^{(S)})+\frac{\hat\rho}{2}\|\widehat\bx^{(S)}-\bx^{(S)}\|^2\leq f(\widehat\bx^{(T)})+\frac{\hat\rho}{2}\|\widehat\bx^{(T)}-\bx^{(S)}\|^2\leq f(\widehat\bx^{(T)})+\frac{\hat\rho D^2}{2}.$ Note that $D$ here is defined as in Assumption~\ref{assume:convex_sharpness} instead of Assumption~\ref{assume:convex}. Applying these bounds to \eqref{eq:mainthmeq6_c_old}, we have
\small
\begin{align}
\nonumber
&~\sum_{t=S}^{T-1}\left[\eta_{t}\hat{\rho}(\hat{\rho}-\rho)\|\widehat\bx^{(t)}-\bx^{(t)} \|^{2}-\eta_{t}\hat{\rho}\widehat\lambda_{t}g (\bx^{(t)} )\right]\mathbb{I}(g (\bx^{(t)} )\leq \epsilon_t)+\sum_{t=S}^{T-1}\eta_{t}\hat{\rho}g (\bx^{(t)} )\mathbb{I}(g (\bx^{(t)} )> \epsilon_t)\\\label{eq:mainthmeq6_c}
\leq&~\frac{\hat\rho D^2}{2}
+\dfrac{\hat{\rho}}{2}\sum_{t=S}^{T-1}\left[\eta_{t}^{2}\|\bzt_f^{(t)}\|^{2}\mathbb{I}(g (\bx^{(t)} )\leq \epsilon_t)+\eta_{t}^{2}\|\bzt_g^{(t)}\|^{2}\mathbb{I}(g (\bx^{(t)} )> \epsilon_t)\right].
\end{align}
\normalsize

For $t\in I$, we have from Assumption~\ref{assume:allpaper} and the definition of $\eta_t$ that
\begin{equation}
\label{eq:etagradientnorm1_c}
\eta_{t}^{2}\|\bzt_f^{(t)}\|^{2}\leq\eta_{t}^2M^2= \eta_{t} M^2 \frac{\nu'}{2M^2}\min\left\{\epsilon^2/M,\nu'\right\}=\eta_{t}\epsilon_t.
\end{equation}
For $t\in J$, we have $\eta_t=\nu'g_+(\bx^{(t)})/(2M\|\bzt_g^{(t)}\|^2)$ so
\begin{equation}
\label{eq:etagradientnorm2_c}
\eta_{t}^{2}\|\bzt_g^{(t)}\|^{2}= \eta_{t}\nu'g_+(\bx^{(t)})/(2M)\leq \eta_{t}\nu'\text{dist}(\bx^{(t)},\mathcal{S})/2\leq\eta_{t}\epsilon_t,
\end{equation}
where the first inequality is by $M$-Lipschitz continuity of $g_+$ and the last inequality is from Proposition~\ref{thm:feasiblesubproblem_c}. 


We first prove that, if $S$, $T$, $\eta_t$ and $\epsilon_t$ are chosen such that 
\begin{eqnarray}
\label{eq:propeq4_c}
\epsilon_t(1+\widehat\lambda_t)\leq C'^2\epsilon^2(\hat\rho-\rho)
\end{eqnarray}
and 
\small
\begin{eqnarray}
\label{eq:propeq9_c}
\sum_{t=S}^{T-1}\eta_{t}\hat{\rho}\epsilon_t> ~\frac{\hat\rho D^2}{2}
+\dfrac{\hat{\rho}}{2}\sum_{t=S}^{T-1}\left[\eta_{t}^{2}\|\bzt_f^{(t)}\|^{2}\mathbb{I}(g (\bx^{(t)} )\leq \epsilon_t)+\eta_{t}^{2}\|\bzt_g^{(t)}\|^{2}\mathbb{I}(g (\bx^{(t)} )> \epsilon_t)\right],
\end{eqnarray}
\normalsize
we must have $g(\bx^{(t)})\leq\epsilon_t$ for at least one $t$ in $\{S,S+1,\dots,T-1\}$ (i.e., $I\neq\emptyset$) and $\mathbb{E}_\tau[\|\widehat\bx^{(\tau)} - \bx^{(\tau)}\|^2]\leq C'^2\epsilon^2$ (so $\mathbb{E}_\tau[\|\widehat\bx^{(\tau)} - \bx^{(\tau)}\|]\leq C'\epsilon$). 

Suppose \eqref{eq:propeq9_c} holds and $g(\bx^{(t)})>\epsilon_t$ for $t=S,S+1,\dots,T-1$, i.e., $I=\emptyset$. \eqref{eq:mainthmeq6_c} becomes exactly the opposite of \eqref{eq:propeq9_c}. This contradiction means  $g(\bx^{(t)})\leq\epsilon_t$ for at least one $t$ in $\{S,S+1,\dots,T-1\}$. Suppose \eqref{eq:propeq4_c} and \eqref{eq:propeq9_c} hold but $\mathbb{E}_\tau[\|\widehat\bx^{(\tau)} - \bx^{(\tau)}\|^2]>C'^2\epsilon^2$. Since $\tau$ is generated by Output I, we have
\begin{align}
\label{eq:stationaryopposite_c}
C'^2\epsilon^2<\mathbb{E}_\tau[\|\widehat\bx^{(\tau)} - \bx^{(\tau)}\|^2]
=\frac{\sum_{t=S}^{T-1}\eta_{t}\mathbb{I}(g(\bx^{(t)})\leq\epsilon_t)\|\widehat\bx^{(t)}-\bx^{(t)} \|^{2}}{\sum_{t=S}^{T-1}\eta_{t}\mathbb{I}(g(\bx^{(t)})\leq\epsilon_t)}.
\end{align}
\eqref{eq:stationaryopposite_c} and \eqref{eq:propeq4_c} imply
\small
\begin{align}
\nonumber
&~\sum_{t=S}^{T-1}\left[\eta_{t}\hat{\rho}(\hat{\rho}-\rho)\|\widehat\bx^{(t)}-\bx^{(t)} \|^{2}-\eta_{t}\hat{\rho}\widehat\lambda_{t}\epsilon_t\right]\mathbb{I}(g(\bx^{(t)})\leq \epsilon_t)+\sum_{t=S}^{T-1}\eta_{t}\hat{\rho}\epsilon_t\mathbb{I}(g(\bx^{(t)})> \epsilon_t)\\\nonumber
>&~\sum_{t=S}^{T-1}\left[\eta_{t}\hat{\rho}(\hat{\rho}-\rho)C'^2\epsilon^{2}-\eta_{t}\hat{\rho}\widehat\lambda_{t}\epsilon_t\right]\mathbb{I}(g(\bx^{(t)})\leq \epsilon_t)+\sum_{t=S}^{T-1}\eta_{t}\hat{\rho}\epsilon_t\mathbb{I}(g(\bx^{(t)})> \epsilon_t)\\\label{eq:leftright}
\geq&~\sum_{t=S}^{T-1}\eta_{t}\hat{\rho}\epsilon_t\mathbb{I}(g(\bx^{(t)})\leq \epsilon_t)+\sum_{t=S}^{T-1}\eta_{t}\hat{\rho}\epsilon_t\mathbb{I}(g(\bx^{(t)})> \epsilon_t)
=\sum_{t=S}^{T-1}\eta_{t}\hat{\rho}\epsilon_t,
\end{align}
\normalsize
where the second inequality is because of \eqref{eq:propeq4_c}. Combining this inequality and \eqref{eq:mainthmeq6_c} leads to the opposite of \eqref{eq:propeq9_c}. This contradiction means $\mathbb{E}_\tau[\|\widehat\bx^{(\tau)} - \bx^{(\tau)}\|^2]\leq C'^2\epsilon^2$.

Given the result above, we only need to show that the the choices of $S$, $T$, $\eta_t$ and $\epsilon_t$ ensure \eqref{eq:propeq4_c} and \eqref{eq:propeq9_c}.

Plugging in the value of $\epsilon_t$ in Proposition~\ref{thm:feasiblesubproblem_c} and using the second result in Lemma~\ref{thm:boundlambda_c}, we have
\begin{align*}
    \epsilon_t(1+\widehat\lambda_t)\leq \epsilon_t(1+\Lambda'')\leq \frac{\nu'(1+\Lambda'')\epsilon^2}{2M}
\end{align*}
which implies \eqref{eq:propeq4_c} by the definition of  $C'$. According to the definition of $\eta_t$ for $t\in I$ and $t\in J$, we have 
\begin{align*}
    \eta_t\geq\frac{\nu'}{2M}\eta_t=\frac{\nu'^2}{4M^3}\min\left\{\epsilon^2/M,\nu'\right\}
\end{align*}
for $t\in I$ by noting that $\nu'\leq M$ and
\begin{align*}
    \eta_t\geq\frac{\nu'}{2M}\cdot\frac{\epsilon_t}{M^2}=\frac{\nu'^2}{4M^3}\min\left\{\epsilon^2/M,\nu'\right\}
\end{align*}
for $t\in J$. By the choice of $T$, we have
\begin{align*}
    \frac{\hat\rho}{2}\sum_{t=S}^{T-1}\eta_t\epsilon_t\geq\frac{\hat\rho T}{2}\frac{\nu'^3}{8M^3}\min\left\{\epsilon^4/M^2,\nu'^2\right\}>\frac{\hat\rho D^2}{2}
\end{align*}
which further implies \eqref{eq:propeq9_c} together with \eqref{eq:etagradientnorm1_c} and \eqref{eq:etagradientnorm2_c}. Thus, we have $\mathbb{E}_\tau[\|\widehat\bx^{(\tau)} - \bx^{(\tau)}\|]\leq C'\epsilon$. The inequality $\mathbb{E}_\tau [g(\bx^{(\tau)})]\leq\frac{\nu'}{2}\min\left\{\epsilon^2/M,\nu'\right\}$ is a direct consequence of $g(\bx^{(t)})\leq\epsilon_t$ for any $t\in I$ and the choice of $\epsilon_t$.
\end{proof}

\section{Convergence analysis for weakly convex constraints}
\label{sec:weaklyconvexproof}

In this section, we present two non-trivial practical examples that satisfy Assumption~\ref{assume:weaklyconvex} and then analyze the complexity of Algorithm~\ref{alg:dsgm} when $f$ and $g$ are both weakly convex. 

\subsection{Examples that satisfy Assumption~\ref{assume:weaklyconvex} }
\label{sec:slaterexample}
In this section, we present two examples that satisfy Assumption~\ref{assume:weaklyconvex}, especially, Assumption~\ref{assume:weaklyconvex}B.

\subsubsection{Demographic parity constraint}
Problem  \eqref{eq:demographicparity_fairnessClassification} is an instance of \eqref{eq:gco} with 
\begin{align*}
\X=\mathbb{R}^d, ~
f(\bx)=L(\bx)+\lambda\text{SCAD}(\bx), ~
g(\bx)=R_0(\bx)-\kappa,
\end{align*}
where $L$, $R_0$ and $\text{SCAD}$ are defined in \eqref{eq:ermL}, \eqref{eq:demographicparity} and \eqref{eq:SCAD}, respectively. It is easy to verify that problem \eqref{eq:demographicparity_fairnessClassification} satisfies Assumption~\ref{assume:allpaper}. In this section, 
we will prove that problem \eqref{eq:demographicparity_fairnessClassification} also satisfies Assumption~\ref{assume:weaklyconvex}. 

Since $L(\bx)$ is convex and $\text{SCAD}(\bx)$ is $2$-weakly convex, $f(\bx)$ is $2\lambda$-weakly convex. Let 
\small
\begin{equation}
\label{eq:demographicparity_h}
h(\bx)=\frac{1}{n_p}\sum_{i=1}^{n_p} \sigma(\bx^\top\ba_i^p)-\frac{1}{n_u}\sum_{i=1}^{n_u} \sigma(\bx^\top\ba_i^u)
\end{equation}
\normalsize
so $g(\bx)=|h(\bx)|-\kappa$. Recalling that $\sigma(z)=\exp(z)/(1+\exp(z))$, we have
\small
\begin{equation}
\label{eq:demographicparity_h_grad}
\nabla h(\bx)=\frac{1}{n_p}\sum_{i=1}^{n_p} \sigma(\bx^\top\ba_i^p)(1-\sigma(\bx^\top\ba_i^p))\ba_i^p-\frac{1}{n_u}\sum_{i=1}^{n_u} \sigma(\bx^\top\ba_i^u)(1-\sigma(\bx^\top\ba_i^u))\ba_i^u.
\end{equation}
\normalsize
It is easy to prove that $h(\bx)$ is Lipschitz continuous with a constant of 
\small
\begin{equation}
\label{eq:hLip0}
\alpha=\frac{1}{4n_p}\sum_{i=1}^{n_p}\|\ba_i^p\|+\frac{1}{4n_u}\sum_{i=1}^{n_u}\|\ba_i^u\|
\end{equation}
\normalsize
and that $\nabla h(\bx)$ is Lipschitz continuous with a constant of 
\small
\begin{equation}
\label{eq:hLip}
\beta=\frac{1}{4n_p}\sum_{i=1}^{n_p}\|\ba_i^p\|^2+\frac{1}{4n_u}\sum_{i=1}^{n_u}\|\ba_i^u\|^2,
\end{equation}
\normalsize
so $g$ is $\beta$-weakly convex according to Lemma 4.2 in \cite{drusvyatskiy2019efficiency}. Hence, \eqref{eq:demographicparity_fairnessClassification} satisfies Assumption~\ref{assume:weaklyconvex}A with $\rho=\max\{2\lambda,\beta\}$.

Recall $\mathcal{D}_p=\{\ba_i^p\}_{i=1}^{n_p}$ and $\mathcal{D}_u=\{\ba_i^u\}_{i=1}^{n_u}$ are  the feature vectors of protected and unprotected groups, respectively. We assume without generality that the first feature of  $\ba_i^p$, denoted by  $[\ba_i^p]_1$, equals one and the first feature of  $\ba_i^u$, denoted by  $[\ba_i^u]_1$, equals negative one. In fact, this feature can be the one used to split data into the protected group or the unprotected group (see the group variables in Table~\ref{tbl:data} for example). We can also simply add this feature to the data if it does not exist originally. With this feature, we can show that \eqref{eq:demographicparity_fairnessClassification} satisfies Assumption~\ref{assume:weaklyconvex}B under mild conditions. 

\begin{theorem}
Suppose $[\ba_i^p]_1=1$ for $i=1,\dots,n_p$ and $[\ba_i^u]_1=-1$ for $i=1,\dots,n_u$ and
\begin{align*}
\kappa\notin\mathcal{K}=\left\{\left|\frac{i}{n_p}-\frac{j}{n_u}\right| ~\Big|~ i=0,1,\dots,n_p,~ j=0,1,\dots,n_u\right\}.
\end{align*}
The following statements hold.
\begin{enumerate}
\item[1.] There exist $\overline\kappa$ and $\underline\kappa$ in $\mathcal{K}$ such that $(\underline\kappa, \overline\kappa)$ contains $\kappa$ but no numbers in $\mathcal{K}$. 
\item[2.] It holds that
\begin{align}
\label{eq:qgeqzero}
\underline{q}:=\inf\left\{q(\bx)~\big|~|h(\bx)|\in [(\kappa+\underline\kappa)/2,(\kappa+\overline\kappa)/2 ]\right\}>0,
\end{align} 
where $h$ is defined in \eqref{eq:demographicparity_h} and 
\small
\begin{align}
\label{eq:def}
q(\bx)=\max\left\{ \max\limits_{i=1,\dots,n_p}\frac{ \sigma(\bx^\top\ba_i^p)(1-\sigma(\bx^\top\ba_i^p)) }{n_p}, 
\max\limits_{i=1,\dots,n_u}\frac{ \sigma(\bx^\top\ba_i^u)(1-\sigma(\bx^\top\ba_i^u)) }{n_u}
\right\}.
\end{align}
\normalsize
\item[3.] Problem \eqref{eq:demographicparity_fairnessClassification} satisfies Assumption~\ref{assume:weaklyconvex}B with 
\small
\begin{align*}
\theta:=\min\left\{
\frac{\kappa-\underline\kappa}{2}, \frac{\underline{q}^2}{4(\beta+\bar\rho)},\frac{\underline{q}^2\kappa}{16\alpha^2}
\right\}
\end{align*}
\normalsize
for any $\bar\rho>\rho$ and any $\bar\epsilon>0$ that satisfy 
\small
\begin{align}
\label{eq:condbarepsilon}
\bar\epsilon^2\leq \min\left\{\frac{\overline\kappa-\kappa}{2},  \frac{\underline{q}^2}{4(\beta+\bar\rho)},\frac{\underline{q}^2\kappa}{16\alpha^2}\right\}.
\end{align}
\normalsize
\end{enumerate}
\end{theorem}
\begin{proof}
Since $\kappa\notin\mathcal{K}$ and $\mathcal{K}$ is discrete, Statement 1 holds trivially. In fact, one can sort the numbers in $\mathcal{K}$ and let $\overline\kappa$ and $\underline\kappa$ be the two consecutive numbers with $\kappa\in (\underline\kappa, \overline\kappa)$. 

Suppose $\underline{q}=0$. There must exist a sequence $\{\bx^{(t)}\}_{t\geq0}$ such that $|h(\bx^{(t)})|\in [(\kappa+\underline\kappa)/2,(\kappa+\overline\kappa)/2 ]$ for any $t\geq0$ and $\lim\limits_{t\rightarrow \infty} q(\bx^{(t)})=0$. Since $\sigma(z)\in (0,1)$ for any $z\in\mathbb{R}$, any sequence like $\sigma(z^{(t)})(1-\sigma(z^{(t)}))$ will converge to zero only when all limiting points of $\sigma(z^{(t)})$ are in $\{0,1\}$. By this observation and the definition of $q$ in \eqref{eq:def}, after passing to a subsequence if necessary, we have
\begin{align}
\lim\limits_{t\rightarrow \infty} \sigma(\bx^{(t)\top}\ba_i^p)\in \{0, 1\}, i=1,\dots,n_p~\text{ and }~
\lim\limits_{t\rightarrow \infty} \sigma(\bx^{(t)\top}\ba_i^u)\in \{0, 1\}, i=1,\dots,n_u.
\end{align}
This means $\lim\limits_{t\rightarrow \infty} |h(\bx^{(t)})|\in\mathcal{K}$, contradicting with Statement 1 and 
the fact that $|h(\bx^{(t)})|\in [(\kappa+\underline\kappa)/2,(\kappa+\overline\kappa)/2 ]\subset (\underline\kappa,\overline\kappa)$. This contradiction indicates that  $\underline{q}>0$, which proves Statement 2.


Consider any $\bar\rho>\rho$ and any $\bar\epsilon>0$ that satisfy \eqref{eq:condbarepsilon}. Because $[\ba_i^p]_1=1$ for $i=1,\dots,n_p$
and $[\ba_i^u]_1=-1$ for $i=1,\dots,n_u$, it is easy to show that 
\small
\begin{align*}
\|\nabla h(\bx)\|\geq \left|\frac{1}{n_p}\sum_{i=1}^{n_p} \sigma(\bx^\top\ba_i^p)(1-\sigma(\bx^\top\ba_i^p))[\ba_i^p]_1-\frac{1}{n_u}\sum_{i=1}^{n_u} \sigma(\bx^\top\ba_i^u)(1-\sigma(\bx^\top\ba_i^u))[\ba_i^u]_1\right|
\geq q(\bx),
\end{align*}
\normalsize
where $\nabla h$ is given in \eqref{eq:demographicparity_h_grad}. As a result, it holds that
\small
\begin{align}
\nonumber
&~\inf\left\{\|\nabla h(\bx)\|~\big|~g(\bx)\in [(\underline\kappa-\kappa)/2, \bar\epsilon^2 ]\right\}
\geq \inf\left\{q(\bx)~\big|~g(\bx)\in [(\underline\kappa-\kappa)/2, \bar\epsilon^2 ]\right\}\\ \nonumber
=&~\inf\left\{q(\bx)~\big|~|h(\bx)|\in [(\kappa+\underline\kappa)/2, \kappa+\bar\epsilon^2 ]\right\}
\geq \inf\left\{q(\bx)~\big|~|h(\bx)|\in [(\kappa+\underline\kappa)/2,(\kappa+\overline\kappa)/2 ]\right\}\\\label{eq:defl}
=&~\underline{q}>0,
\end{align}
\normalsize
where the second inequality is because $\bar\epsilon^2\leq(\overline\kappa-\kappa)/2$ so $ \kappa+\bar\epsilon^2 \leq(\kappa+\overline\kappa)/2$. 

Consider a solution $\bx$ with $g(\bx)\leq \bar\epsilon^2$. Suppose $g(\bx)\leq -(\kappa-\underline\kappa)/2$. We can set $\by=\bx$ and have 
\small
\begin{align}
\label{eq:uniformslater1}
g(\by)+\frac{\bar\rho}{2}\|\by-\bx\|^2=g(\bx)\leq -\frac{\kappa-\underline\kappa}{2}\leq -\theta.
\end{align}
\normalsize
Suppose $g(\bx)\geq -(\kappa-\underline\kappa)/2$ so $|h(\bx)|\geq  (\underline\kappa+\kappa)/2>0$. We assume $h(\bx)>0$ and the proof when $h(\bx)<0$ is the same. Let $\by$ be generated from $\bx$ by running a gradient descent step on $h(\bx)$ with a stepsize of $\eta>0$, namely, $\by=\bx-\eta \nabla h(\bx)$. We choose $\eta= \min\{1/(\beta+\bar\rho),\kappa/(4\alpha^2)\}$. Since $h(\by)+\frac{\bar\rho}{2}\|\by-\bx\|^2$ is smooth and its gradient with respect to $\by$ is $(\beta+\bar\rho)$-Lipschitz continuous, we have
\small
\begin{align}
\nonumber
h(\by)+\frac{\bar\rho}{2}\|\by-\bx\|^2
\leq&~ 
h(\bx)+\frac{\bar\rho}{2}\|\bx-\bx\|^2+\left\langle \nabla h(\bx)+\bar\rho(\bx-\bx), \by-\bx\right\rangle+\frac{\beta+\bar\rho}{2}\|\by-\bx\|^2\\\nonumber
\leq&~
\kappa+\bar\epsilon^2 -\eta\left(1-\frac{\beta+\bar\rho}{2}\eta\right)\| \nabla h(\bx)\|^2\\\nonumber
\leq&~\kappa+\bar\epsilon^2 -\frac{1}{2}\min\left\{\frac{1}{\beta+\bar\rho},\frac{\kappa}{4\alpha^2}\right\}\| \nabla h(\bx)\|^2\\
\label{eq:hgradientdescent}
\leq&~\kappa+\bar\epsilon^2 -\frac{\underline{q}^2}{2}\min\left\{\frac{1}{\beta+\bar\rho},\frac{\kappa}{4\alpha^2}\right\},
\end{align}
\normalsize
where the second inequality is by the definition of $\by$ and the fact that $h(\bx) =|h(\bx) |=g(\bx)+\kappa\leq \kappa+\bar\epsilon^2$, the third inequality is because $\eta= \min\{1/(\beta+\bar\rho),\kappa/(4\alpha^2)\}$, and the last is by \eqref{eq:defl} and the fact that $g(\bx)\in [(\underline\kappa-\kappa)/2, \bar\epsilon^2 ]$. Moreover, since $h(\bx)=|h(\bx)|=g(\bx)+\kappa\geq (\underline\kappa+\kappa)/2$, $\eta\leq \kappa/(4\alpha^2)$ and $h$ is $\alpha$-Lipschitz continuous, we have 
\small
\begin{align*}
h(\by)\geq h(\bx)- \alpha\eta\|\nabla h(\bx)\|\geq \frac{\underline\kappa+\kappa}{2}-\eta\alpha^2\geq \frac{\kappa}{4}>0,
\end{align*}
\normalsize
which indicates $g(\by)=|h(\by)|-\kappa=h(\by)-\kappa$. Hence, \eqref{eq:hgradientdescent} and \eqref{eq:condbarepsilon} imply
\small
\begin{align}
\label{eq:uniformslater2}
g(\by)+\frac{\bar\rho}{2}\|\by-\bx\|^2\leq \bar\epsilon^2-\frac{\underline{q}^2}{2}\min\left\{\frac{1}{\beta+\bar\rho},\frac{\kappa}{4\alpha^2}\right\}\leq -\frac{\underline{q}^2}{4}\min\left\{\frac{1}{\beta+\bar\rho},\frac{\kappa}{4\alpha^2}\right\}\leq -\theta.
\end{align}
\normalsize
Inequalities \eqref{eq:uniformslater1} and \eqref{eq:uniformslater2} mean  \eqref{eq:demographicparity_fairnessClassification} satisfies Assumption~\ref{assume:weaklyconvex}B, which proves Statement 3. 
\end{proof}

\subsubsection{Smoothly clipped absolute deviation constraint}
Suppose $f(\bx)$  in \eqref{eq:gco} is any training loss that satisfies Assumptions~\ref{assume:allpaper}, \ref{assume:weaklyconvex}A and \ref{assume:weaklyconvex}C. We consider a SCAD constraint for promoting sparsity of model $\bx$. It is formulated as
\begin{align}
\label{eq:scadconstraint}
g(\bx)=\text{SCAD}(\bx)-\kappa=\sum_{i=1}^ds(x_i)-\kappa\leq 0,
\end{align}
where $s(x_i)$ is defined in \eqref{eq:SCAD}. Jia and Grimmer~\cite{jia2022first} show that, if $\kappa$ is not divisible by three,  constraint \eqref{eq:scadconstraint} satisfies the strong Mangasarian-Fromovitz constraint qualification (MFCQ), which is used in~\cite{jia2022first} to bound the Lagrangian multiplier in a way similar to \eqref{eq:Lambdabounda}. Motivated by the finding in~\cite{jia2022first},  we show that 
\eqref{eq:scadconstraint} also satisfies the uniform Slater's condition in Assumption~\ref{assume:weaklyconvex}B. 
\begin{proposition}
Suppose $\kappa$ is not divisible by three. The following statements hold.
\begin{enumerate}
\item[1.]There exists a non-negative integer $k$ such that $\kappa\in (3k, 3(k+1))$. 
\item[2.] It holds that
\begin{align}
\label{eq:gnorm_lowerbound_SCAD}
\inf\left\{\max_{i}\textup{dist}(x_i,\mathcal{Z})~\big|~ \textup{SCAD}(\bx)\in [(\kappa+3k)/2, (\kappa+3(k+1))/2] \right\}=:\underline q>0,
\end{align}
where 
$\mathcal{Z}=(-\infty, -2] \cup \{0\} \cup [2, \infty)$.
\item[3.] Constraint \eqref{eq:scadconstraint} satisfies Assumption~\ref{assume:weaklyconvex}B with 
\small
\begin{align*}
\theta:=\min\left\{
\frac{\kappa-3k}{2}, \frac{\underline{q}^2}{2+\bar\rho}, \frac{\underline{q}^3}{2}
\right\}
\end{align*}
\normalsize
for any $\bar\rho>\rho$ and any $\bar\epsilon>0$ that satisfy 
\small
\begin{align}
\label{eq:condbarepsilon_SCAD}
\bar\epsilon^2\leq \min\left\{
\frac{3(k+1)-\kappa}{2}, \frac{\underline{q}^2}{2+\bar\rho}, \frac{\underline{q}^3}{2}
\right\}.
\end{align}
\normalsize
\end{enumerate}
\end{proposition}
\begin{proof}
Since $\kappa$ is not divisible by three, Statement 1 holds trivially. In addition, by \eqref{eq:condbarepsilon_SCAD}, we also have 
\begin{align}
\label{eq:interval_scad}
3k< \frac{\kappa+3k}{2} <\kappa <\kappa+\bar\epsilon^2\leq \frac{\kappa+3(k+1)}{2} < 3(k+1).
\end{align}

As observed by~\cite{jia2022first}, $0\in \partial s(x_i)$ only when $x_i\in \mathcal{Z}$, or equivalently, $s(x_i)\in\{0,3\}$ for any $i$. Suppose \eqref{eq:gnorm_lowerbound_SCAD} is not true.  There must exist a sequence $\{\bx^{(t)}\}_{t\geq0}$ such that $\text{SCAD}(\bx^{(t)})\in [(\kappa+3k)/2, (\kappa+3(k+1))/2] $ for any $t\geq0$ and $\lim\limits_{t\rightarrow \infty}\text{dist}(x_i^{(t)},\mathcal{Z})=0$ for any $i$. Therefore, by passing to a subsequence if necessary, we have $\lim\limits_{t\rightarrow \infty}s(x_i^{(t)})\in\{0,3\}$ for $i=1,\dots,d$, so $\lim\limits_{t\rightarrow \infty}\text{SCAD}(\bx^{(t)})$ is divisible by three, contradicting with \eqref{eq:interval_scad} and the fact that $\text{SCAD}(\bx^{(t)})\in [(\kappa+3k)/2, (\kappa+3(k+1))/2]$. This contradiction proves Statement 2. 

Consider a solution $\bx$ with $g(\bx)\leq \bar\epsilon^2$. Suppose $g(\bx)\leq -(\kappa-3k)/2$. We can set $\by=\bx$ and have 
\small
\begin{align}
\label{eq:uniformslater1_SCAD}
g(\by)+\frac{\bar\rho}{2}\|\by-\bx\|^2=g(\bx)\leq -\frac{\kappa-3k}{2}\leq -\theta.
\end{align}
\normalsize
Suppose $g(\bx)\geq -(\kappa-3k)/2$ so $ \text{SCAD}(\bx)\in [(\kappa+3k)/2, (\kappa+3(k+1))/2] $. By \eqref{eq:gnorm_lowerbound_SCAD}, there exists an index $i$ such that $\text{dist}(x_i,\mathcal{Z})\geq \underline q$, which implies $x_i\in [\underline q, 2-\underline q]$ or $x_i\in [-2+\underline q, -\underline q]$. Without loss of generality, we assume $i=1$ and $x_i\in [\underline q, 2-\underline q]$ since the proof when $x_1\in [-2+\underline q, -\underline q]$ is the same. By the definition of $s$ in \eqref{eq:SCAD}, we have $s'(x_1)\in [2 \underline q, 1]$. Let $y_i=x_i$ for $i\neq 1$ and $y_1=x_1-\eta s'(x_1)$ with $\eta= \min\{\underline q/2, 1/(2+\bar\rho)\}$. We then have $y_1\geq x_1-\underline q/2\geq \underline q/2$. Since the function $s(y_1)+\frac{\bar\rho}{2}(y_1-x_1)^2$ is smooth on $[\underline q/2,x_1]$ and its gradient with respect to $y_1$ is $(2+\bar\rho)$-Lipschitz continuous, we have   
\small
\begin{align}
\nonumber
g(\by)+\frac{\bar\rho}{2}\|\by-\bx\|^2
\leq&~ 
g(\bx)+s'(x_1)\cdot (y_1-x_1)+\frac{2+\bar\rho}{2}(y_1-x_1)^2\\\nonumber
\leq&~
\bar\epsilon^2 -\eta\left(1-\frac{2+\bar\rho}{2}\eta\right)[s'(x_1)]^2\\\nonumber
\leq&~\bar\epsilon^2  -\frac{1}{2}\min\left\{\frac{\underline q}{2}, \frac{1}{2+\bar\rho}\right\}[s'(x_1)]^2\\
\label{eq:hgradientdescent_SCAD}
\leq&~\bar\epsilon^2  -2\min\left\{\frac{\underline q}{2}, \frac{1}{2+\bar\rho}\right\}\underline q^2
\leq -\min\left\{\frac{\underline q}{2}, \frac{1}{2+\bar\rho}\right\}\underline q^2\leq -\theta,
\end{align}
\normalsize
where the second inequality is by the definition of $y_1$ and the fact $g(\bx)\leq\bar\epsilon^2$, the third inequality is because $\eta= \min\{\underline q/2, 1/(2+\bar\rho)\}$, the fourth is because $s'(x_1)\in [2 \underline q, 1]$, the fifth is because $\bar\epsilon^2 \leq \min\{\underline q/2, 1/(2+\bar\rho)\}\underline q^2$, and the last by the definition of $\theta$. Inequalities \eqref{eq:uniformslater1_SCAD} and \eqref{eq:hgradientdescent_SCAD} mean  \eqref{eq:scadconstraint} satisfies Assumption~\ref{assume:weaklyconvex}B, which proves Statement 3. 
\end{proof}

\subsection{Proofs of Lemmas~\ref{thm:boundlambda_wc},~\ref{thm:sharpsubgradient} and~\ref{thm:errorbound_wc} }
\label{sec:boundlambda_wc}
\begin{proof}[\textbf{Proof of Lemma~\ref{thm:boundlambda_wc}}]
For simplicity of notation, we denote $\widehat\bx(\bx)$ in \eqref{eq:phix} by $\widehat\bx$. Suppose $\bx$ is $\epsilon^2$-feasible. According to Assumption~\ref{assume:weaklyconvex}B, there exists $\by\in\text{relint}(\X)$ such that 
\begin{eqnarray}
\label{eq:slator1}
 g(\by) + \frac{\hat\rho}{2}\| \by - \bx\|^2\leq 
 g(\by) + \frac{\bar\rho}{2}\| \by - \bx\|^2\leq
 -\theta.
\end{eqnarray}
This means \eqref{eq:phi} satisfies the Slater's condition, so there exists $\widehat\lambda\geq 0$ that satisfies the KKT condition together with $\widehat\bx$. In particular, we have $\widehat\lambda (g(\widehat{\bx})+\frac{\hat\rho}{2}\|\widehat{\bx}-\bx\|^2) = 0$ and
\begin{align}
		\label{eq:optcond1_wc}
  \widehat\bzt_f  +\hat{\rho}(\widehat\bx - \bx)+\widehat\lambda (\widehat\bzt_g + \hat{\rho} (\widehat\bx -\bx ))+\widehat\bu=\mathbf{0},
\end{align}
where $\widehat\bzt_f\in\partial f(\widehat\bx)$, $\widehat\bzt_g\in\partial g(\widehat\bx)$, $\widehat\bu\in \mathcal{N}_\X(\widehat\bx)$ and $\mathcal{N}_\X(\widehat\bx)$ is the normal cone of $\X$ at $\widehat\bx$. Taking inner product between \eqref{eq:optcond1_wc} and $\widehat\bx - \bx$ gives 
\begin{align*}
 0\geq& -\left\langle\widehat\bu,\widehat\bx - \bx\right\rangle = \left\langle\widehat\bzt_f+\widehat\lambda \widehat\bzt_g,\widehat\bx - \bx\right\rangle +\hat{\rho}(1+\widehat\lambda)\|\widehat\bx - \bx\|^2\\
 \geq& -\frac{\|\widehat\bzt_f+\widehat\lambda \widehat\bzt_g\|^2}{2\hat{\rho}(1+\widehat\lambda)} -\frac{\hat{\rho}(1+\widehat\lambda)}{2}\|\widehat\bx - \bx\|^2+\hat{\rho}(1+\widehat\lambda)\|\widehat\bx - \bx\|^2,
\end{align*}
where the first inequality is because  $\widehat\bu\in \mathcal{N}_\X(\widehat\bx)$, the second inequality is by Young's inequality. Reorganizing the terms in this inequality and using the facts that  $\|\widehat\bzt_f\|\leq M$ and $\|\widehat\bzt_g\|\leq M$, we obtain  
\begin{align*}
\hat{\rho}(1+\widehat\lambda)\|\widehat\bx - \bx\|^2\leq \frac{\|\widehat\bzt_f+\widehat\lambda \widehat\bzt_g\|^2}{\hat{\rho}(1+\widehat\lambda)}\leq \frac{(1+\widehat\lambda)M^2}{\hat\rho},
\end{align*}
which further implies the first inequality in \eqref{eq:Lambdabounda}.

If $\widehat\lambda = 0$, the conclusion holds trivially. Hence, we focus on the case that $\widehat\lambda > 0$. Note that, in this case, we must have $g(\widehat{\bx})+\frac{\hat\rho}{2}\|\widehat{\bx}-\bx\|^2=0$.
  

		

Since $g(\bz) + \frac{\hat{\rho}}{2}\| \bz - \bx \|^2+\delta_{\X}(\bz) $ is $(\hat{\rho} - \rho)$-strongly convex in $\bz$ and $\widehat\bu/\widehat\lambda\in \mathcal{N}_\X(\widehat\bx)=\partial\delta_{\X}(\widehat\bx)$,  we have
\small
\begin{align*}
    &~g(\by) + \frac{\hat{\rho}}{2}\| \by  - \bx  \|^2\\
    \geq&~ g(\widehat\bx ) + \frac{\hat{\rho}}{2}\| \widehat\bx  -\bx  \|^2 + \left\langle \widehat\bzt_g + \hat{\rho} (\widehat\bx  -\bx  )+\widehat\bu/\widehat\lambda, \by - \widehat\bx  \right\rangle + \frac{\hat{\rho} - \rho}{2} \| \by - \widehat\bx  \|^2 \\
	=&~\left\langle \widehat\bzt_g + \hat{\rho} (\widehat\bx  -\bx  )+\widehat\bu/\widehat\lambda, \by - \widehat\bx  \right\rangle + \frac{\hat{\rho} - \rho}{2} \| \by - \widehat\bx  \|^2.
\end{align*}
\normalsize
Applying \eqref{eq:slator1} to the inequality above and arranging terms give
\small
\begin{align*}
    -\theta \geq \left\langle \widehat\bzt_g + \hat{\rho} (\widehat\bx  -\bx  )+\widehat\bu/\widehat\lambda, \by - \widehat\bx  \right\rangle + \frac{\hat{\rho} - \rho}{2} \| \by - \widehat\bx  \|^2  
	\geq - \frac{ \| \widehat\bzt_g + \hat{\rho} (\widehat\bx  -\bx  )+\widehat\bu/\widehat\lambda \|^2 }{2(\hat{\rho} - \rho)},
\end{align*}
\normalsize
which, together with \eqref{eq:optcond1_wc} and  the first inequality in \eqref{eq:Lambdabounda}, implies 
\small
\begin{align*}
	\widehat\lambda = \frac{\| \widehat\bzt_f + \hat{\rho} (\widehat\bx  -\bx  ) \| }{\| \widehat\bzt_g + \hat{\rho} (\widehat\bx  -\bx  ) +\widehat\bu/\widehat\lambda\| } 
	\leq  \frac{2M}{\sqrt{2 \theta (\hat{\rho} - \rho)}}.
\end{align*}
\normalsize
\end{proof}

\begin{proof}[\textbf{Proof of Lemma~\ref{thm:sharpsubgradient}}]
Consider $\bx\in\mathcal{L}$. Since $\bx$ is $\epsilon^2$-feasible, by Assumption~\ref{assume:weaklyconvex}B, there exists $\by\in\text{relint}(\X)$ such that $g(\by)+\frac{\bar\rho}{2}\|\by-\bx\|^2\leq-\theta$. Note that function $g(\bz)+\frac{\hat\rho}{2}\|\bz-\bx\|^2+\delta_{\X}(\bz)$ is $(\hat\rho-\rho)$-strongly convex with respect to $\bz$ and its subdifferential with respect to $\bz$ at location $\bz=\bx$ is $\partial g(\bx)+\mathcal{N}_{\X}(\bx)$. We then have 
\begin{align*}
g(\bx)+\frac{\hat\rho}{2}\|\bx-\bx\|^2+ \langle \bzt_g+\bu, \by-\bx \rangle+\frac{\hat\rho-\rho}{2}\|\by-\bx\|^2
\leq g(\by)+\frac{\hat\rho}{2}\|\by-\bx\|^2\leq -\theta
\end{align*}
for any $\bzt_g\in\partial g(\bx)$ and any $\bu\in \mathcal{N}_{\X}(\bx)$. Since $g(\bx)=0$ when $\bx\in\mathcal{L}$, applying Young's inequality to the inequality above yields 
\begin{align*}
-\frac{\|\bzt_g+\bu\|^2}{2(\hat\rho-\rho)}\leq \langle \bzt_g+\bu, \by-\bx \rangle+\frac{\hat\rho-\rho}{2}\|\by-\bx\|^2\leq -\theta,
\end{align*}
which implies the conclusion. 
\end{proof}

\begin{proof}[\textbf{Proof of Lemma~\ref{thm:errorbound_wc}}]
For  any $\bx\in\X$ satisfying $\text{dist}(\bx,\mathcal{S})\leq \frac{\nu}{\rho}$, we define
\begin{equation}
\label{eq:projtoS}
\bx^\dagger=\textup{proj}_{\mathcal{S}}(\bx)\in\argmin_{\by\in\X,g(\by)\leq0}\frac{1}{2}\|\by-\bx\|^2.
\end{equation} 
Since $g$ can be non-convex, $\bx^\dagger$ is not necessarily unique, but this proof works for any $\bx^\dagger$ satisfying \eqref{eq:projtoS}. Since the first conclusion of this lemma holds trivially if $\bx\in \mathcal{S}$, we assume $\bx\in\X\backslash \mathcal{S}$, which implies $g(\bx^\dagger)=0$ and $\bx^\dagger\in\mathcal{L}$. It is easy to show that $\bx^\dagger$ is also the optimal solution of 
\begin{equation}
\label{eq:projtoS_new}
 \min_{\by\in\X} \Big\{ \frac{1}{2}\|\by - \bx\|^2, 
	~{s.t.}~ g(\by)+\frac{\bar{\rho}}{2}\|\by - \bx^\dagger\|^2\leq 0 \Big\},
\end{equation} 
which is a convex optimization problem. By Assumption~\ref{assume:weaklyconvex}B, there exists $\by\in\text{relint}(\X)$ such that 
$g(\by)+\frac{\bar\rho}{2}\|\by-\bx^\dagger\|^2\leq-\theta$. Hence, the Slater's condition holds for \eqref{eq:projtoS_new}, so there exist a scalar $\lambda^\dagger\geq0$, a subgradient $\bzt_g\in\partial g(\bx^\dagger)$ and
a vector $\bu\in\mathcal{N}_{\X}(\bx^\dagger)$ such that 
\begin{equation}
\label{eq:projtoS_KKT}
\bx^\dagger-\bx+\lambda^\dagger \bzt_g+\bu=\mathbf{0}.
\end{equation} 
Since $\langle\bu,\bx-\bx^\dagger\rangle\leq 0$, we assert that $\lambda^\dagger>0$ because, otherwise,  
\eqref{eq:projtoS_KKT} implies $0=\langle\bx^\dagger-\bx+\bu,\bx-\bx^\dagger\rangle\leq -\|\bx-\bx^\dagger\|^2\leq 0$ and thus $\bx=\bx^\dagger$, contradicting with the fact that $\bx\in\X\backslash \mathcal{S} $.

By the $\rho$-weak convexity of $g$ and \eqref{eq:projtoS_KKT}, we have
\begin{align*}
&~\lambda^\dagger\left(g(\bx)-g(\bx^\dagger)\right)\geq \lambda^\dagger\left(\langle\bzt_g,\bx-\bx^\dagger\rangle
-\frac{\rho}{2}\|\bx-\bx^\dagger\|^2\right)\\
\geq&~\langle\lambda^\dagger\bzt_g+\bu,\bx-\bx^\dagger\rangle-\frac{\lambda^\dagger\rho}{2}\|\bx-\bx^\dagger\|^2
=\|\bx-\bx^\dagger\|^2-\frac{\lambda^\dagger\rho}{2}\|\bx-\bx^\dagger\|^2.
\end{align*}
Since $\lambda^\dagger>0$, dividing both sides of the inequalities above leads to 
\begin{align*}
&~g(\bx)-g(\bx^\dagger)\geq \frac{\|\bx-\bx^\dagger\|^2}{\lambda^\dagger}-\frac{\rho}{2}\|\bx-\bx^\dagger\|^2\\
=&~\|\bx-\bx^\dagger\|\cdot\|\bzt_g+\bu/\lambda^\dagger\|-\frac{\rho}{2}\|\bx-\bx^\dagger\|^2
\geq \nu \|\bx-\bx^\dagger\|-\frac{\rho}{2}\|\bx-\bx^\dagger\|^2\geq\frac{\nu}{2}\|\bx-\bx^\dagger\|,
\end{align*}
where the second inequality is by Lemma~\ref{thm:sharpsubgradient} and the last inequality is because $\|\bx-\bx^\dagger\|=\text{dist}(\bx,\mathcal{S})\leq \frac{\nu}{\rho}$. The first conclusion is thus proved by the facts that $g(\bx)=g_+(\bx)$ for $\bx\in\X\backslash\mathcal{S}$ and that $g(\bx^\dagger)=0$. The second conclusion is directly from Lemma 3.1 and 3.2 in \cite{davis2018subgradient}. 
\end{proof}

\subsection{Proof of Proposition~\ref{thm:feasiblesubproblem} and Theorem~\ref{thm:mdconverge_deterministic_weaklyconvex}}
\label{sec:proof_mainresult_wc}


To prove  Proposition~\ref{thm:feasiblesubproblem} and Theorem~\ref{thm:mdconverge_deterministic_weaklyconvex}, we need the following lemma from Theorem 4.1 in \cite{davis2018subgradient} whose proof is provided only for completeness. It shows that, when the Polyak's stepsize is used, the subgradient method can utilize the local error bound condition in Lemma~\ref{thm:errorbound_wc} to ensure $\text{dist}(\bx^{(t)},\mathcal{S})$ Q-linearly converges to zero, which makes sure  $\textup{dist}(\bx^{(t)},\mathcal{S})$ is small even when $t\in J$ and prevents $\bx^{(t)}$ from being trapped in an infeasible stationary point given the second conclusion of Lemma~\ref{thm:errorbound_wc}.
\begin{lemma}[\cite{davis2018subgradient}]
\label{thm:nonexpansion}
Suppose Assumptions~\ref{assume:allpaper} and \ref{assume:weaklyconvex} hold. Also, suppose the sequence 
$\{\bx^{(t)}\}_{t\geq0}$ is generated by applying the projected subgradient method to $\min_{\bx\in\X}g_+(\bx)$ using a Polyak's stepsize, namely, 
\begin{equation}
\label{eq:Polyak}
\bx^{(t+1)}=\textup{proj}_{\X}(\bx^{(t)}-\eta_t\bzt_g^{(t)}), \hspace{0.5em} \eta_t=
\left\{
\begin{array}{ll}
g_+(\bx^{(t)})/\|\bzt_g^{(t)}\|^2 & \text{ if }~\bzt_g^{(t)}\neq \mathbf{0}\\
0& \text{ if }~\bzt_g^{(t)}= \mathbf{0}
\end{array}
\right.,
\hspace{0.5em} \text{ for }t=0, 1, \dots,
\end{equation}
where $\bzt_g^{(t)}\in\partial g_+(\bx^{(t)})$.  If $\textup{dist}(\bx^{(0)},\mathcal{S})\leq \nu/(4\rho)$,  we have 
\begin{eqnarray}
\label{eq:shrinkdist_Q}
\textup{dist}^2(\bx^{(t+1)},\mathcal{S})\leq\left(1-\frac{\nu^2}{8M^2}\right)\textup{dist}^2(\bx^{(t)},\mathcal{S})
\end{eqnarray}
and
\begin{eqnarray}
\label{eq:shrinkdist}
\textup{dist}(\bx^{(t)},\mathcal{S})\leq \textup{dist}(\bx^{(0)},\mathcal{S}),\quad\forall t\geq0.
\end{eqnarray}
\end{lemma}
\begin{proof}

We prove \eqref{eq:shrinkdist} by induction. \eqref{eq:shrinkdist} holds trivially for $t=0$. Suppose that \eqref{eq:shrinkdist} holds up to iteration $t$ of \eqref{eq:Polyak}. We want to prove that it also holds for iteration $t+1$. By the induction hypothesis, we have $\text{dist}(\bx^{(t)},\mathcal{S})\leq \nu/(4\rho)$. 

Suppose $\bzt_g^{(t)}=\mathbf{0}$ so $\bx^{(t+1)}=\bx^{(t)}$. By Lemma~\ref{thm:errorbound_wc}, we must have $\textup{dist}(\bx^{(t)},\mathcal{S})=0$ and thus $\textup{dist}(\bx^{(t+1)},\mathcal{S})=0\leq\text{dist}(\bx^{(0)},\mathcal{S})$ and \eqref{eq:shrinkdist_Q} holds.  

Suppose $\bzt_g^{(t)}\neq\mathbf{0}$.
Let $\bx^{\dagger(t)}=\text{proj}_{\mathcal{S}}(\bx^{(t)})$ for $t\geq0$. By the nonexpansiveness  of $\text{prox}_{\X}(\cdot)$, we have
\small
\begin{align*}
\nonumber
   \text{dist}^2(\bx^{(t+1)},\mathcal{S})
   \leq&~ \|\bx^{(t+1)}-\bx^{\dagger(t)}\|^{2}= \|\text{proj}_{\X} (\bx^{(t)}-\eta_{t}\bzt_g^{(t)} )-\text{proj}_{\X} (\bx^{\dagger(t)})\|^{2} \\ \nonumber
   \leq&~ \|\bx^{(t)}-\eta_{t}\bzt_g^{(t)}-\bx^{\dagger(t)}  \|^{2}\\\nonumber
    =&~ 
    \text{dist}^2(\bx^{(t)},\mathcal{S})-2\eta_{t}\left\langle\bzt_g^{(t)},\bx^{(t)}-\bx^{\dagger(t)}\right\rangle+\eta_{t}^{2} \|\bzt_g^{(t)}\|^{2}\\\nonumber
   \leq&~ \text{dist}^2(\bx^{(t)},\mathcal{S})+\frac{2 g_+(\bx^{(t)})}{\|\bzt_g^{(t)}\|^2}\left(g_+(\bx^{\dagger(t)})-g_+(\bx^{(t)})+\frac{\rho}{2}\|\bx^{(t)}-\bx^{\dagger(t)}\|^2\right)+\frac{g_+^2(\bx^{(t)})}{\|\bzt_g^{(t)}\|^2}\\\label{eq:renegar1}
    =&~ \text{dist}^2(\bx^{(t)},\mathcal{S})+\frac{g_+(\bx^{(t)})}{\|\bzt_g^{(t)}\|^2}\left( -g_+(\bx^{(t)})+\rho\|\bx^{(t)}-\bx^{\dagger(t)}\|^2\right)\\\nonumber
    \leq&~ \text{dist}^2(\bx^{(t)},\mathcal{S})+\frac{g_+(\bx^{(t)})}{\|\bzt_g^{(t)}\|^2}\left( -\frac{\nu}{2}+\rho\|\bx^{(t)}-\bx^{\dagger(t)}\|\right)\|\bx^{(t)}-\bx^{\dagger(t)}\|\\\nonumber
    \leq&~ \text{dist}^2(\bx^{(t)},\mathcal{S})-\frac{ g_+(\bx^{(t)})}{\|\bzt_g^{(t)}\|^2}\frac{\nu}{4}\|\bx^{(t)}-\bx^{\dagger(t)}\|\\\nonumber
    \leq&~ \text{dist}^2(\bx^{(t)},\mathcal{S})-\frac{\nu^2}{8M^2}\|\bx^{(t)}-\bx^{\dagger(t)}\|^2,
\end{align*}
\normalsize
where the second equality is because $g_+(\bx^{\dagger(t)})=0$, the third inequality is by the $\rho$-weak convexity of $g_+$, the fourth is by Lemma~\ref{thm:errorbound_wc}, the fifth is by the hypothesis that $\text{dist}(\bx^{(t)},\mathcal{S})\leq\text{dist}(\bx^{(0)},\mathcal{S})\leq \nu/(4\rho)$, and the last is by  Lemma~\ref{thm:errorbound_wc} and Assumption~\ref{assume:allpaper}A. 
This inequality further implies 
\begin{eqnarray*}
\textup{dist}^2(\bx^{(t+1)},\mathcal{S})\leq\left(1-\frac{\nu^2}{8M^2}\right)\textup{dist}^2(\bx^{(t)},\mathcal{S})\leq\textup{dist}^2(\bx^{(0)},\mathcal{S}),
\end{eqnarray*}
which proves \eqref{eq:shrinkdist_Q} and also proves \eqref{eq:shrinkdist} by induction.
\end{proof}

\begin{proof}[\textbf{Proof of Proposition~\ref{thm:feasiblesubproblem}}]
The choices of $\epsilon_t$ and $\eta_t$ imply 
\begin{equation}
\label{eq:stepetainI}
2\epsilon_t/\nu+2\eta_tM^2/\nu\leq\min\left\{\epsilon^2/M,\nu/(4\rho)\right\}
\end{equation}
for $t\in I$. We first prove $\text{dist}(\bx^{(t)},\mathcal{S})\leq \min\left\{\epsilon^2/M,\nu/(4\rho)\right\}$  by induction on $t$. Since $\bx^{(0)}\in \mathcal{S}$, this conclusion holds trivially for $t=0$. Suppose it holds up to iteration $t$. We want to prove it also holds for iteration $t+1$. 


Suppose $t\in I$. Recall that $\nu\leq 2M$. We have $\eta_tM\leq 2\eta_tM^2/\nu\leq\nu/(4\rho)$ and
$
\|\bx^{(t+1)}-\bx^{(t)}\|\leq 
\|\eta_t\bzt_f^{(t)}\|\leq \eta_tM, 
$
which, by triangle inequality and the induction hypothesis, implies that 
\begin{align*}
    \text{dist}(\bx^{(t+1)},\mathcal{S})\leq \text{dist}(\bx^{(t)},\mathcal{S})+\eta_tM\leq\frac{\nu}{4\rho}+ \eta_tM\leq\frac{\nu}{2\rho}.
\end{align*}
Since $g(\bx^{(t)})\leq \epsilon_t$ for $t\in I$, we have $g(\bx^{(t+1)})\leq g(\bx^{(t)})+M\|\bx^{(t+1)}-\bx^{(t)}\|\leq \epsilon_t+\eta_tM^2$. By Lemma~\ref{thm:errorbound_wc}, we have   $\text{dist}(\bx^{(t+1)},\mathcal{S})\leq 2\epsilon_t/\nu+2\eta_tM^2/\nu\leq \min\left\{\epsilon^2/M,\nu/(4\rho)\right\}$.  

Suppose $t\in J$. Let $t'$ be the largest index in $I$ that is smaller than $t$. By the same proof as in the previous case, we have $\text{dist}(\bx^{(t'+1)},\mathcal{S})\leq  \min\left\{\epsilon^2/M,\nu/(4\rho)\right\}$.  Since indexes $t'+1$, $t'+2$, ..., $t$ are in $J$, Algorithm~\ref{alg:dsgm} essentially performs the projected subgradient method to $\min_{\bx\in\X}g_+(\bx)$ using a Polyak's stepsize during iterations $t'+1$, $t'+2$, ..., and $t$.  Hence, by Lemma~\ref{thm:nonexpansion}, we have
$\text{dist}(\bx^{(t+1)},\mathcal{S})\leq\text{dist}(\bx^{(t'+1)},\mathcal{S}) \leq \min\left\{\epsilon^2/M,\nu/(4\rho)\right\}$.    

By induction, we have prove that $\text{dist}(\bx^{(t)},\mathcal{S})\leq \min\left\{\epsilon^2/M,\nu/(4\rho)\right\}$ for any $t\geq0$. As a result, $g(\bx^{(t)})+\frac{\tilde\rho}{2}\|\bx^{(t)}-\bx^{(t)}\|^2=g(\bx^{(t)})\leq M\cdot\text{dist}(\bx^{(t)},\mathcal{S})\leq M\cdot\min\left\{\epsilon^2/M,\nu/(4\rho)\right\}\leq \epsilon^2$, meaning that $\bx^{(t)}$ is $\epsilon^2$-feasible for \eqref{eq:phi}.  
\end{proof}

The following proposition provides the main inequality needed for proving Theorem~\ref{thm:mdconverge_deterministic_weaklyconvex}. 

\begin{proposition}
\label{thm:mainprop_weaklyconvex}
Under the same assumptions as  Proposition~\ref{thm:feasiblesubproblem},  Algorithm~\ref{alg:dsgm} guarantees
\small
\begin{align*}
\nonumber
\sum_{t=S}^{T-1}\dfrac{\eta_{t}\hat{\rho}(\hat{\rho}-\rho)}{2} \|\widehat\bx^{(t)}-\bx^{(t)} \|^{2}
\leq&~f(\bx^{(S)})-\underline f+\frac{3M^2}{2\hat\rho}
+\sum_{t\in I}(1+\widehat\lambda_{t+1})\eta_{t}\hat{\rho}\widehat\lambda_{t}\epsilon_t\\\nonumber
&+\dfrac{\hat{\rho}}{2}\sum_{t=S}^{T-1}(1+\widehat\lambda_{t+1})\left[\eta_{t}^{2}\|\bzt_f^{(t)}\|^{2}\mathbb{I}(g(\bx^{t})\leq \epsilon_t)+\eta_{t}^{2}\|\bzt_g^{(t)}\|^{2}\mathbb{I}(g(\bx^{t})> \epsilon_t)\right].
\end{align*}
\normalsize
\end{proposition}
\begin{proof}
Let $\bzt^{(t)}=\bzt_f^{(t)}\in\partial f(\bx^{(t)})$ if $t\in I$ and  $\bzt^{(t)}=\bzt_g^{(t)}\in\partial g(\bx^{(t)})$ if $t\in J$. Let $\varphi(\bx)$ and $\widehat{\bx}$ be defined in \eqref{eq:phi} and \eqref{eq:phix} with $(\hat\rho,\tilde\rho)$ satisfying \eqref{eq:parameter2}. For simplicity of notation, we denote $\widehat\bx(\bx^{(t)})$ by $\widehat\bx^{(t)}$. Since $\bx^{(t)}$ is $\epsilon^2$-feasible by Proposition~\ref{thm:feasiblesubproblem}, we get from Assumption~\ref{assume:weaklyconvex}B that $\varphi(\bx^{(t)})$ and $\widehat\bx^{(t)}$ are well defined for any $t\geq0$. 


By Assumption~\ref{assume:weaklyconvex}B, problem \eqref{eq:phi} with $\bx=\bx^{(t)}$ is strongly convex and has a strictly feasible solution, so there exists a Lagrangian multiplier $\widehat\lambda_t\geq 0$ satisfying \eqref{eq:KKTprox} and
\begin{eqnarray}
\label{eq:optimalityhatx_wc}
\widehat\bx^{(t)}=\argmin_{\bx\in\X}\left\{f(\bx)+\dfrac{\hat{\rho}}{2} \|\bx-\bx^{(t)} \|^{2}+\widehat\lambda_t\left(g(\bx)+\frac{\hat\rho}{2}\| \bx -\bx^{(t)}\|^2\right)\right\}.
\end{eqnarray}

By the updating equation of  $\bx^{(t+1)}$, we have
\small
\begin{align*}
     \|\bx^{(t+1)}-\widehat\bx^{(t)} \|^{2}=&~ \|\text{proj}_{\X} (\bx^{(t)}-\eta_{t}\bzt^{(t)} )-\widehat\bx^{(t)} \|^{2}= \|\text{proj}_{\X} (\bx^{(t)}-\eta_{t}\bzt^{(t)} )-\text{proj}_{\X} (\widehat\bx^{(t)})\|^{2} \\
    \leq&~ \|\bx^{(t)}-\eta_{t}\bzt^{(t)}-\widehat\bx^{(t)}  \|^{2}= 
    \|\bx^{(t)}-\widehat\bx^{(t)} \|^{2}-2\eta_{t}\left\langle\bzt^{(t)},\bx^{(t)}-\widehat\bx^{(t)}\right\rangle+\eta_{t}^{2} \|\bzt^{(t)}\|^{2}.
\end{align*}
\normalsize
Multiplying the inequality above by $\hat{\rho}(1+\widehat\lambda_{t+1})/2$ and adding $f (\widehat\bx^{(t)})+\widehat\lambda_{t+1}g(\widehat\bx^{(t)})$ to both sides, we obtain
\begin{align}
\nonumber
  &~f (\widehat\bx^{(t)} )+\dfrac{\hat{\rho}}{2} \|\bx^{(t+1)}-\widehat\bx^{(t)} \|^{2}
  +\widehat\lambda_{t+1}\left(g(\widehat\bx^{(t)} )+\frac{\hat\rho}{2}\|\bx^{(t+1)}-\widehat\bx^{(t)}\|^2\right)\\\nonumber
  \leq&~f (\widehat\bx^{(t)} )+\dfrac{\hat{\rho}}{2} \|\bx^{(t)}-\widehat\bx^{(t)} \|^{2}+\widehat\lambda_{t+1}\left(g(\widehat\bx^{(t)} )+\frac{\hat\rho}{2}\|\bx^{(t)}-\widehat\bx^{(t)} \|^2\right)\\\nonumber
  &-(1+\widehat\lambda_{t+1})\eta_{t}\hat{\rho}\left\langle\bzt^{(t)},\bx^{(t)}-\widehat\bx^{(t)}\right\rangle+\frac{(1+\widehat\lambda_{t+1})\eta_{t}^{2}\hat{\rho}}{2} \|\bzt^{(t)}\|^{2}\\\label{eq:mainthmeq1_wc}
  \leq&~ \varphi (\bx^{(t)})-(1+\widehat\lambda_{t+1})\eta_{t}\hat{\rho}\left\langle\bzt^{(t)},\bx^{(t)}-\widehat\bx^{(t)}\right\rangle+\frac{(1+\widehat\lambda_{t+1})\eta_{t}^{2}\hat{\rho}}{2} \|\bzt^{(t)}\|^{2},
\end{align}
where the second inequality is by the definition of $\varphi (\bx)$ in \eqref{eq:phi} and the fact that $\widehat\lambda_{t+1}\geq0$ and $g(\widehat\bx^{(t)} )+\frac{\hat\rho}{2}\|\bx^{(t)}-\widehat\bx^{(t)} \|^2\leq 0$. 

By the definitions of $\varphi (\bx^{(t+1)})$ and $\widehat\lambda_{t+1}$, we have the following equation similar to \eqref{eq:optimalityhatx_wc}
\small
\begin{align*}
 \varphi (\bx^{(t+1)})=&~  f (\widehat\bx^{(t+1)} )+\dfrac{\hat{\rho}}{2} \|\bx^{(t+1)}-\widehat\bx^{(t+1)} \|^{2}
  +\widehat\lambda_{t+1}\left(g(\widehat\bx^{(t+1)} )+\frac{\hat\rho}{2}\|\bx^{(t+1)}-\widehat\bx^{(t+1)}\|^2\right)\\
 \leq&~f (\widehat\bx^{(t)} )+\dfrac{\hat{\rho}}{2} \|\bx^{(t+1)}-\widehat\bx^{(t)} \|^{2}
  +\widehat\lambda_{t+1}\left(g(\widehat\bx^{(t)} )+\frac{\hat\rho}{2}\|\bx^{(t+1)}-\widehat\bx^{(t)}\|^2\right),
\end{align*}
\normalsize
which, together with \eqref{eq:mainthmeq1_wc}, implies 
\begin{align}
\label{eq:mainthmeq2_wc}
  (1+\widehat\lambda_{t+1})\eta_{t}\hat{\rho}\left\langle\bzt^{(t)},\bx^{(t)}-\widehat\bx^{(t)}\right\rangle
    \leq&~\varphi (\bx^{(t)})-\varphi (\bx^{(t+1)})+\dfrac{(1+\widehat\lambda_{t+1})\eta_{t}^{2}\hat{\rho}}{2}\|\bzt^{(t)}\|^{2}.
\end{align}
Next, we will bound $\left\langle\bzt^{(t)},\bx^{(t)}-\widehat\bx^{(t)}\right\rangle$ from below when  $t\in I$ and $t\in J$, separately. 

Suppose $t\in I$ so $g(\bx^{(t)})\leq \epsilon_t$ and $\bzt^{(t)}=\bzt_f^{(t)}$. By the $\rho$-weak convexity of $f$, we have  
\begin{align}
\nonumber
\left\langle\bzt^{(t)},\bx^{(t)}-\widehat\bx^{(t)}\right\rangle
\geq&~ f (\bx^{(t)} )-f (\widehat\bx^{(t)} )-\dfrac{\rho}{2} \|\widehat\bx^{(t)}-\bx^{(t)} \|^{2}\\\label{eq:mainthmeq3_wc}
=&~ f (\bx^{(t)} )-f (\widehat\bx^{(t)} )-\dfrac{\hat\rho}{2} \|\widehat\bx^{(t)}-\bx^{(t)} \|^{2}
+\dfrac{\hat\rho-\rho}{2} \|\widehat\bx^{(t)}-\bx^{(t)} \|^{2}.
\end{align}

Since the objective function in \eqref{eq:optimalityhatx_wc} is $(1+\widehat\lambda_{t})(\hat\rho-\rho)$-strongly convex, we have 
\begin{align*}
f(\bx^{(t)})+\widehat\lambda_{t} g(\bx^{(t)})
=&~ f(\bx^{(t)})+\dfrac{\hat{\rho}}{2} \|\bx^{(t)}-\bx^{(t)} \|^{2}+\widehat\lambda_{t} \left(g(\bx^{(t)})+\dfrac{\hat{\rho}}{2} \|\bx^{(t)}-\bx^{(t)} \|^{2}\right)\\
\geq&~f(\widehat\bx^{(t)})+\dfrac{\hat{\rho}}{2} \|\widehat\bx^{(t)}-\bx^{(t)} \|^{2}+\widehat\lambda_{t} \left(g(\widehat\bx^{(t)})+\dfrac{\hat{\rho}}{2} \|\widehat\bx^{(t)}-\bx^{(t)} \|^{2}\right)\\
&+\frac{(1+\widehat\lambda_{t})(\hat{\rho}-\rho)}{2}\|\widehat\bx^{(t)}-\bx^{(t)} \|^{2},
\end{align*}
which, by the facts that $g(\bx^{(t)})\leq \epsilon_t$, $\widehat\lambda_t\geq 0$ and $\widehat\lambda_{t} \big(g(\widehat\bx^{(t)})+\frac{\hat{\rho}}{2} \|\widehat\bx^{(t)}-\bx^{(t)} \|^{2}\big)=0$, implies
\begin{eqnarray*}
f(\bx^{(t)})-f(\widehat\bx^{(t)})-\dfrac{\hat{\rho}}{2} \|\widehat\bx^{(t)}-\bx^{(t)} \|^{2}
\geq -\widehat\lambda_{t}\epsilon_t+\frac{\hat{\rho}-\rho}{2}\|\widehat\bx^{(t)}-\bx^{(t)} \|^{2}.
\end{eqnarray*}
Applying this inequality and inequality \eqref{eq:mainthmeq3_wc} to \eqref{eq:mainthmeq2_wc} leads to 
\begin{align}
\nonumber
&~(1+\widehat\lambda_{t+1})\eta_{t}\hat{\rho}(\hat{\rho}-\rho)  \|\widehat\bx^{(t)}-\bx^{(t)} \|^{2}\\\label{eq:mainthmeq4_wc}
\leq&~ \varphi (\bx^{(t)})-\varphi (\bx^{(t+1)})+\dfrac{(1+\widehat\lambda_{t+1})\eta_{t}^{2}\hat{\rho}}{2}\|\bzt_f^{(t)}\|^{2}+ (1+\widehat\lambda_{t+1})\eta_{t}\hat{\rho}\widehat\lambda_{t}\epsilon_t.
\end{align}

Suppose $t\in J$ so $g(\bx^{(t)})> \epsilon_t$ and $\bzt^{(t)}=\bzt_g^{(t)}$. By the $\rho$-weak convexity of $g$ and the fact that $g (\widehat\bx^{(t)})+\frac{\hat\rho}{2}\|\widehat\bx^{(t)}-\bx^{(t)}\|^2\leq 0$, we have 
\begin{align*}
\nonumber
\left\langle\bzt^{(t)},\bx^{(t)}-\widehat\bx^{(t)}\right\rangle
\geq&~ g (\bx^{(t)} )-g (\widehat\bx^{(t)} )-\dfrac{\rho}{2} \|\widehat\bx^{(t)}-\bx^{(t)} \|^{2}\\
=&~ g (\bx^{(t)} )-g (\widehat\bx^{(t)} )-\dfrac{\hat\rho}{2} \|\widehat\bx^{(t)}-\bx^{(t)} \|^{2}
+\dfrac{\hat\rho-\rho}{2} \|\widehat\bx^{(t)}-\bx^{(t)} \|^{2}\\
>&~ \epsilon_t+\dfrac{\hat\rho-\rho}{2} \|\widehat\bx^{(t)}-\bx^{(t)} \|^{2}.
\end{align*}
Applying this inequality to \eqref{eq:mainthmeq2_wc} leads to 
\begin{align}
\nonumber
&~\dfrac{(1+\widehat\lambda_{t+1})\eta_{t}\hat{\rho}(\hat{\rho}-\rho)}{2} \|\widehat\bx^{(t)}-\bx^{(t)} \|^{2}\\\label{eq:mainthmeq5_wc}
\leq &~\varphi (\bx^{(t)})-\varphi (\bx^{(t+1)})+\dfrac{(1+\widehat\lambda_{t+1})\eta_{t}^{2}\hat{\rho}}{2}\|\bzt_g^{(t)}\|^{2}-(1+\widehat\lambda_{t+1})\eta_{t}\hat{\rho}\epsilon_t.
\end{align}

Summing up \eqref{eq:mainthmeq4_wc} and \eqref{eq:mainthmeq5_wc} for $t=S,S+1,\dots,T-1$, we have
\small
\begin{align}
\nonumber
&~\sum_{t=S}^{T-1}\dfrac{\eta_{t}\hat{\rho}(\hat{\rho}-\rho)}{2} \|\widehat\bx^{(t)}-\bx^{(t)} \|^{2}\leq \sum_{t=S}^{T-1}\dfrac{(1+\widehat\lambda_{t+1})\eta_{t}\hat{\rho}(\hat{\rho}-\rho)}{2} \|\widehat\bx^{(t)}-\bx^{(t)} \|^{2}\\\nonumber
\leq&~ \varphi (\bx^{(S)})-\varphi (\bx^{(T)})+\dfrac{\hat{\rho}}{2}\sum_{t=S}^{T-1}(1+\widehat\lambda_{t+1})\left[\eta_{t}^{2}\|\bzt_f^{(t)}\|^{2}\mathbb{I}(g(\bx^{t})\leq \epsilon_t)+\eta_{t}^{2}\|\bzt_g^{(t)}\|^{2}\mathbb{I}(g(\bx^{t})> \epsilon_t)\right]\\\nonumber
&+\sum_{t\in I}\left[(1+\widehat\lambda_{t+1})\eta_{t}\hat{\rho}\widehat\lambda_{t}\epsilon_t\right]+\sum_{t\in J}\left[-(1+\widehat\lambda_{t+1})\eta_{t}\hat{\rho}\epsilon_t\right]\\\nonumber
\leq&~ \varphi (\bx^{(S)})-\varphi (\bx^{(T)})+\sum_{t\in I}(1+\widehat\lambda_{t+1})\eta_{t}\hat{\rho}\widehat\lambda_{t}\epsilon_t\\\label{eq:mainthmeq6_wc}
&+\dfrac{\hat{\rho}}{2}\sum_{t=S}^{T-1}(1+\widehat\lambda_{t+1})\left[\eta_{t}^{2}\|\bzt_f^{(t)}\|^{2}\mathbb{I}(g(\bx^{t})\leq \epsilon_t)+\eta_{t}^{2}\|\bzt_g^{(t)}\|^{2}\mathbb{I}(g(\bx^{t})> \epsilon_t)\right].
\end{align}
\normalsize

Finally, from Assumption~\ref{assume:weaklyconvex}C and the first inequality in \eqref{eq:Lambdabounda}, it can be easily shown that
\begin{eqnarray*}
\varphi (\bx^{(T)})=f(\widehat\bx^{(T)})+\dfrac{\hat{\rho}}{2} \|\widehat\bx^{(T)}-\bx^{(T)} \|^{2}
\geq \underline f
\end{eqnarray*}
and
\begin{align*}
&~\varphi (\bx^{(S)})=f(\widehat\bx^{(S)})+\dfrac{\hat{\rho}}{2} \|\widehat\bx^{(S)}-\bx^{(S)} \|^{2}\leq f(\bx^{(S)})+\frac{M^2}{\hat\rho}+\dfrac{\hat{\rho}}{2} \frac{M^2}{\hat\rho^2}\leq f(\bx^{(S)})+\frac{3M^2}{2\hat\rho}.
\end{align*}
Then, the conclusion is derived from by applying these two inequalities to \eqref{eq:mainthmeq6_wc}. 
\end{proof}

We are now ready to prove the main theorem for the weakly convex case.
\begin{proof}[\textbf{Proof of Theorem~\ref{thm:mdconverge_deterministic_weaklyconvex}}]
The inequality $\mathbb{E}_\tau [g(\bx^{(\tau)})]\leq \epsilon^2$ is a direct consequence of Proposition~\ref{thm:feasiblesubproblem}.

Applying the upper bound $\Lambda'$ of $\widehat\lambda_t$ from Lemma~\ref{thm:boundlambda_wc} to the inequality in Proposition~\ref{thm:mainprop_weaklyconvex}, we obtain
\small
\begin{align}
\nonumber
\sum_{t=S}^{T-1}\dfrac{\eta_{t}\hat{\rho}(\hat{\rho}-\rho)}{2} \|\widehat\bx^{(t)}-\bx^{(t)} \|^{2}
\leq&~f(\bx^{(S)})-\underline f+\frac{3M^2}{2\hat\rho}
+\sum_{t\in I}(1+\Lambda')\eta_{t}\hat{\rho}\Lambda'\epsilon_t\\\label{eq:tempmainpropo_wc1}
&+\dfrac{\hat{\rho}}{2}\sum_{t=S}^{T-1}(1+\Lambda')\left[\eta_{t}^{2}\|\bzt_f^{(t)}\|^{2}\mathbb{I}(g(\bx^{t})\leq \epsilon_t)+\eta_{t}^{2}\|\bzt_g^{(t)}\|^{2}\mathbb{I}(g(\bx^{t})> \epsilon_t)\right].
\end{align}
\normalsize
For $t\in I$, we have from Assumption~\ref{assume:allpaper} and the definition of $\eta_t$ that
\begin{equation}
\label{eq:etagradientnorm1}
\eta_{t}^{2}\|\bzt_f^{(t)}\|^{2}\leq\eta_{t}^2M^2\leq \eta_{t} M \cdot \frac{2\eta_tM^2}{\nu}\leq \eta_{t}\epsilon^2.
\end{equation}
For $t\in J$, we have $\eta_t=g_+(\bx^{(t)})/\|\bzt_g^{(t)}\|^2$ so
\begin{equation}
\label{eq:etagradientnorm2}
\eta_{t}^{2}\|\bzt_g^{(t)}\|^{2}\leq \eta_{t}g_+(\bx^{(t)})\leq\eta_{t}\epsilon^2,
\end{equation}
where the second inequality is from Proposition~\ref{thm:feasiblesubproblem} that implies $g_+(\bx^{(t)})=g(\bx^{(t)})\leq \epsilon^2$
.

Applying \eqref{eq:etagradientnorm1} and \eqref{eq:etagradientnorm2} to \eqref{eq:tempmainpropo_wc1} and organizing term lead to 
\small
\begin{align}
\nonumber
&~\dfrac{\hat{\rho}(\hat{\rho}-\rho)}{2}\mathbb{E}_\tau[\|\widehat\bx^{(\tau)} - \bx^{(\tau)}\|^2]=\dfrac{\hat{\rho}(\hat{\rho}-\rho)}{2}\frac{\sum_{t=S}^{T-1}\eta_{t}\|\widehat\bx^{(t)}-\bx^{(t)} \|^{2}}{\sum_{t=S}^{T-1}\eta_{t}}\\\nonumber
\leq&~ \frac{1}{\sum_{t=S}^{T-1}\eta_{t}}\left(f(\bx^{(S)})-\underline f+\frac{3M^2}{2\hat\rho}\right)+(1+\Lambda')\hat{\rho}\frac{\sum_{t\in I}\eta_{t}\Lambda'\epsilon_t}{\sum_{t=S}^{T-1}\eta_{t}}+\dfrac{\hat{\rho}(1+\Lambda')}{2}\epsilon^2\\\label{eq:tempmainpropo_wc2}
\leq&~ \frac{1}{\sum_{t=S}^{T-1}\eta_{t}}\left(f(\bx^{(S)})-\underline f+\frac{3M^2}{2\hat\rho}\right)
+\dfrac{3\hat{\rho}(1+\Lambda')}{2}\epsilon^2,
\end{align}
\normalsize
where the second inequality is by the definition of $\epsilon_t$ which implies $\epsilon_t\leq\epsilon^2\nu/(2M)=\epsilon^2/\Lambda'$.
According to the definition of $\eta_t$ for $t\in I$ and $t\in J$, we have 
\small
\begin{align}
\nonumber
\sum_{t=S}^{T-1}\eta_{t}\geq &~ \frac{\nu|I|}{4M^2}\min\left\{\epsilon^2/M,\nu/(4\rho)\right\}+\sum_{t\in J} \frac{ \epsilon_t}{\|\bzt_g^{(t)}\|^2}\\\nonumber
\geq&~  \frac{\nu|I|}{4M^2}\min\left\{\epsilon^2/M,\nu/(4\rho)\right\}+\frac{ \nu|J|}{4M^2}\min\left\{\epsilon^2/M,\nu/(4\rho)\right\}\\\label{eq:lowerbound_sumeta}
=&~\frac{\nu(T-S)}{4M^2}\min\left\{\epsilon^2/M,\nu/(4\rho)\right\}.
\end{align}
\normalsize
Applying \eqref{eq:lowerbound_sumeta} as well as the definitions of $S$, $\eta_t$ and $\epsilon_t$
to the right hand side of \eqref{eq:tempmainpropo_wc2} gives 
\small
\begin{align*}
\mathbb{E}_\tau[\|\widehat\bx^{(\tau)} - \bx^{(\tau)}\|^2]
\leq\frac{8M^2\left(f(\bx^{(0)})-\underline f+3M^2/(2\hat\rho)\right)}{\hat\rho(\hat{\rho}-\rho)\nu\min\left\{\epsilon^2/M,\nu/(4\rho)\right\}T}+\frac{3(1+\Lambda')}{\hat{\rho}-\rho}\epsilon^2,
\end{align*}
\normalsize
which implies the conclusion of this theorem after plugging in the value of $T$. 
\end{proof}

\section{Additional numerical results}
\label{sec:add_exp}

\begin{figure*}
     \begin{tabular}[h]{@{}c|ccc@{}}
      & a9a & bank & COMPAS \\
		\hline \vspace*{-0.1in}\\
		\raisebox{12ex}{\small{\rotatebox[origin=c]{90}{Objective}}}
		& \hspace*{-0.06in}\includegraphics[width=0.30\textwidth]{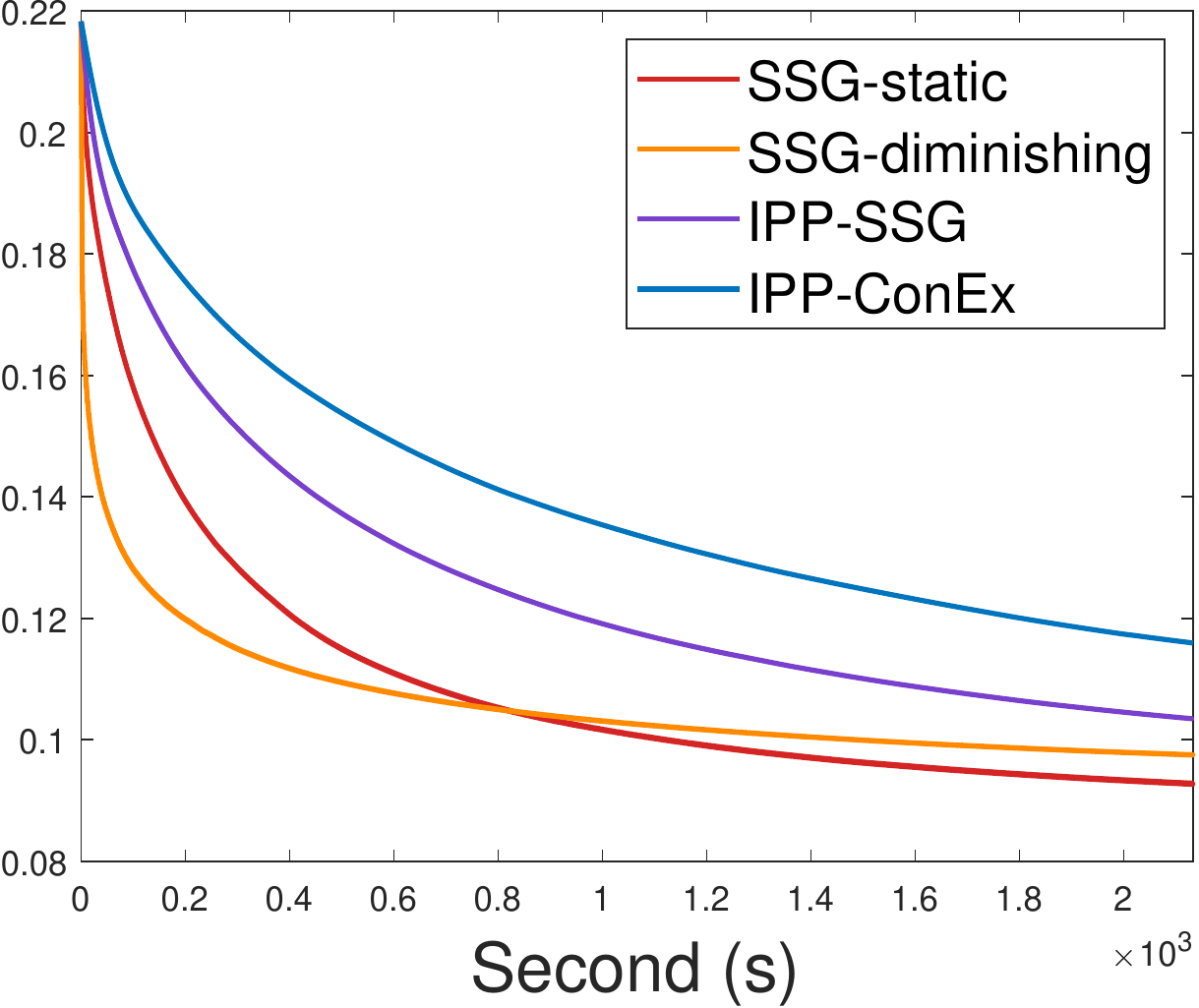}
		& \hspace*{-0.06in}\includegraphics[width=0.30\textwidth]{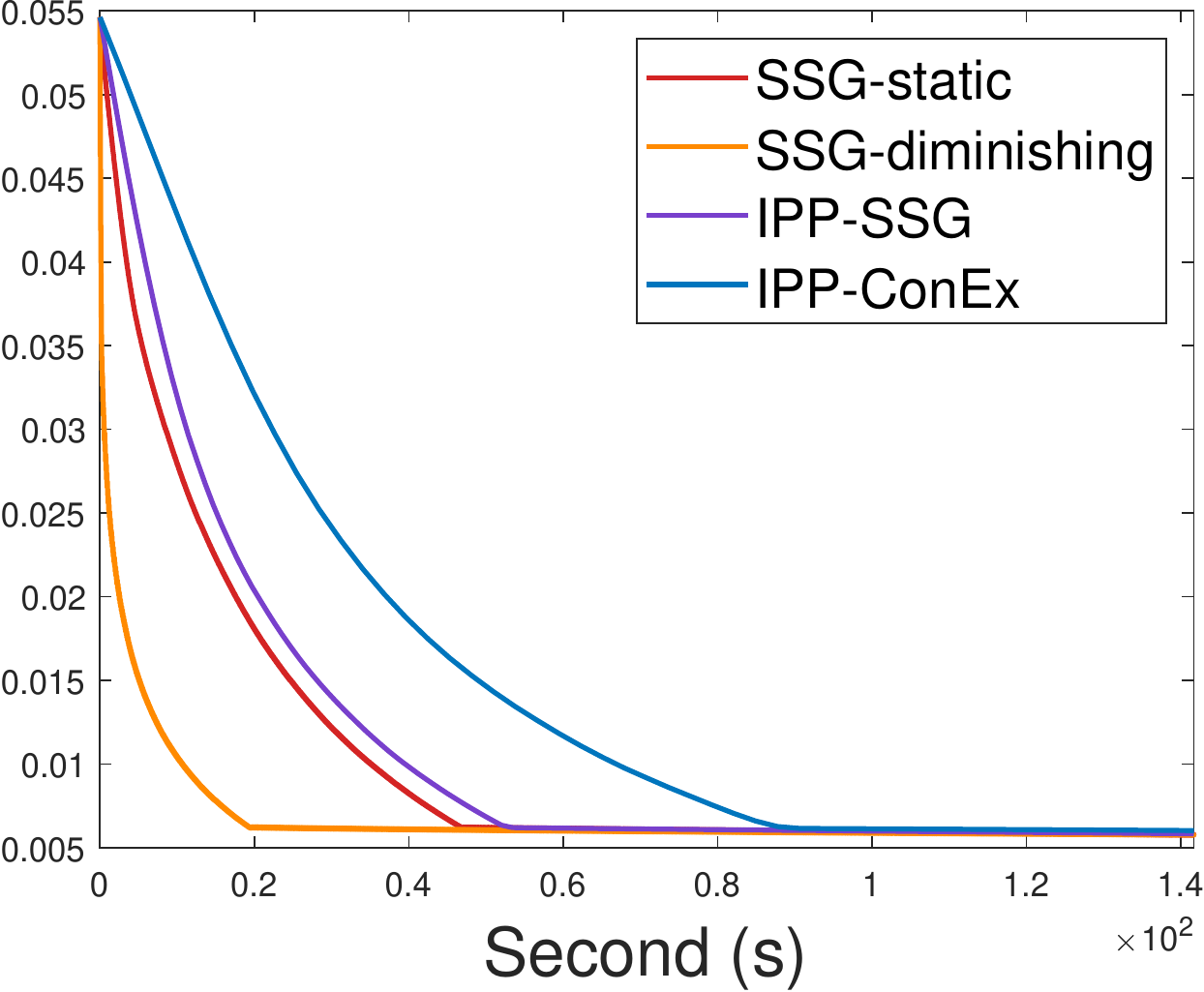}
		& \hspace*{-0.06in}\includegraphics[width=0.30\textwidth]{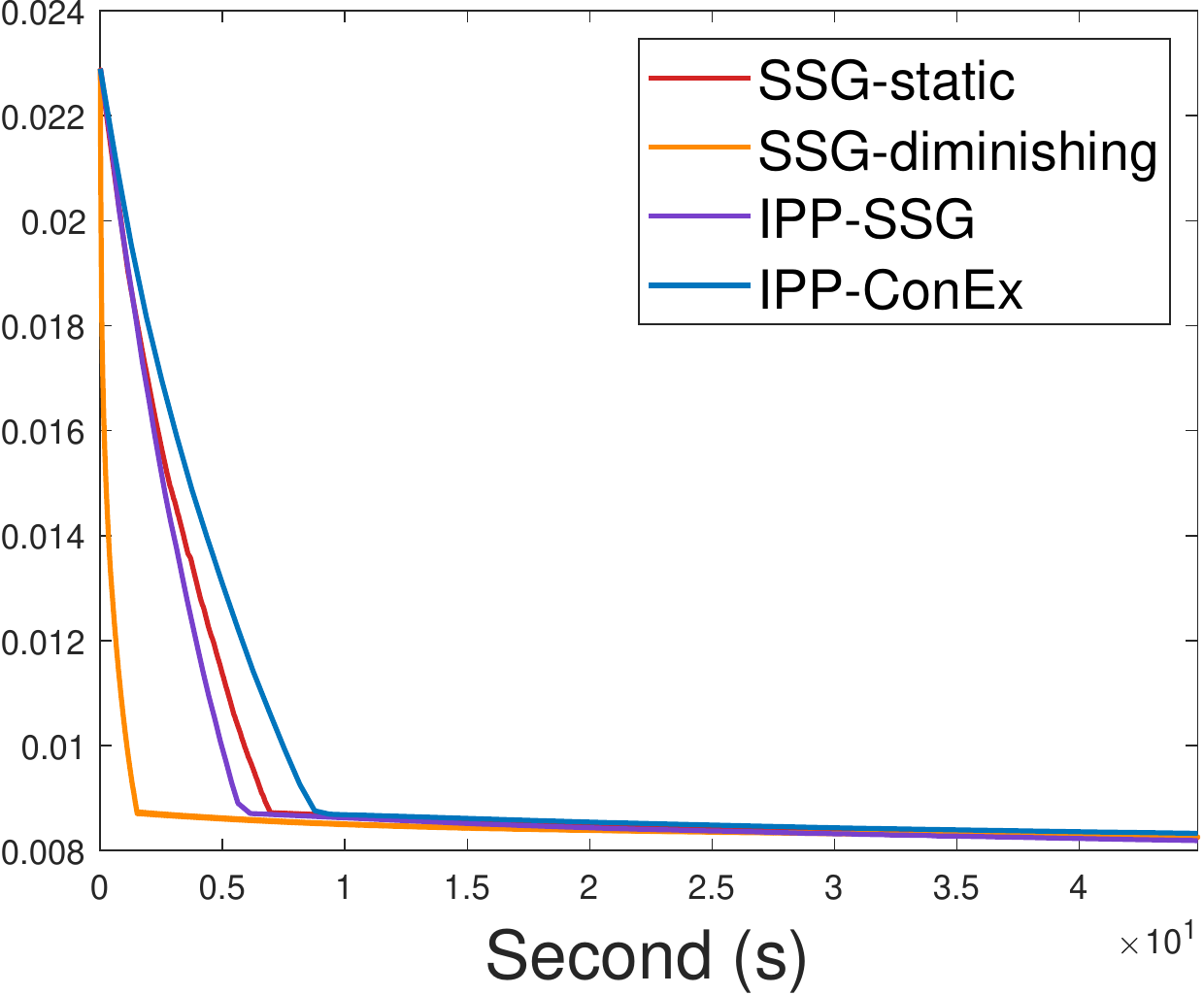}
          \\
		\raisebox{12ex}{\small{\rotatebox[origin=c]{90}{Infeasibility}}}
		& \hspace*{-0.06in}\includegraphics[width=0.30\textwidth]{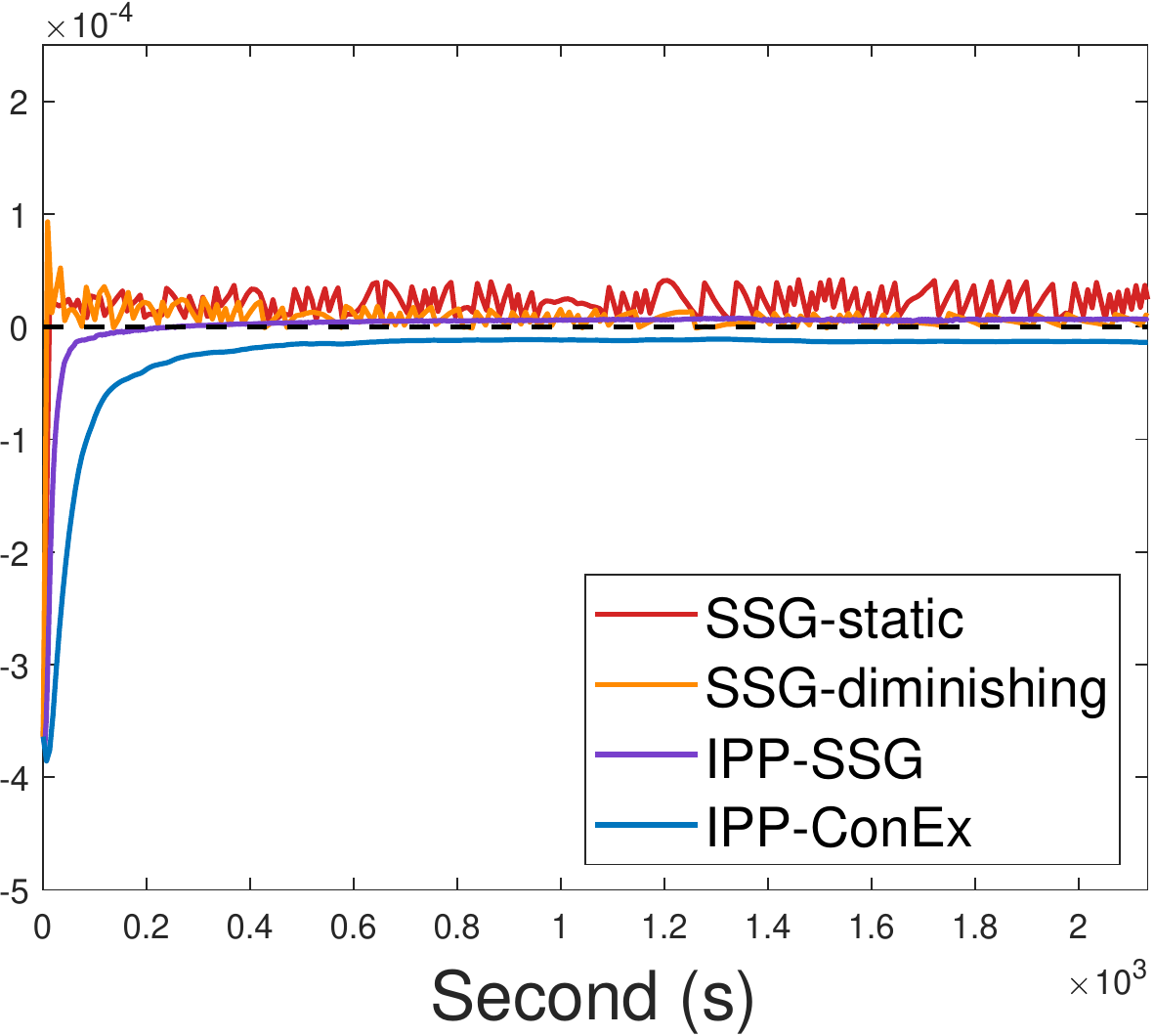}
		& \hspace*{-0.06in}\includegraphics[width=0.30\textwidth]{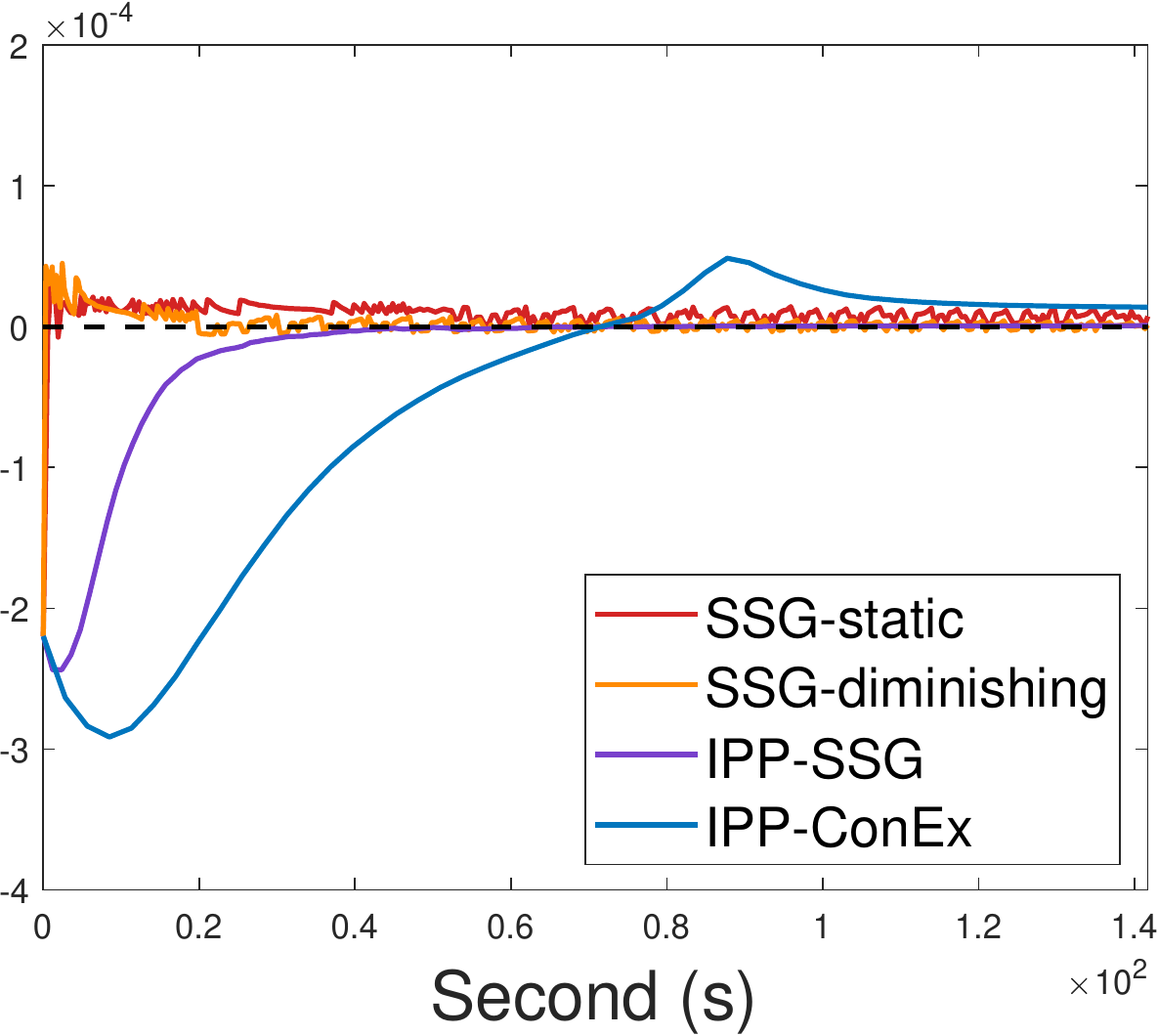}
		& \hspace*{-0.06in}\includegraphics[width=0.30\textwidth]{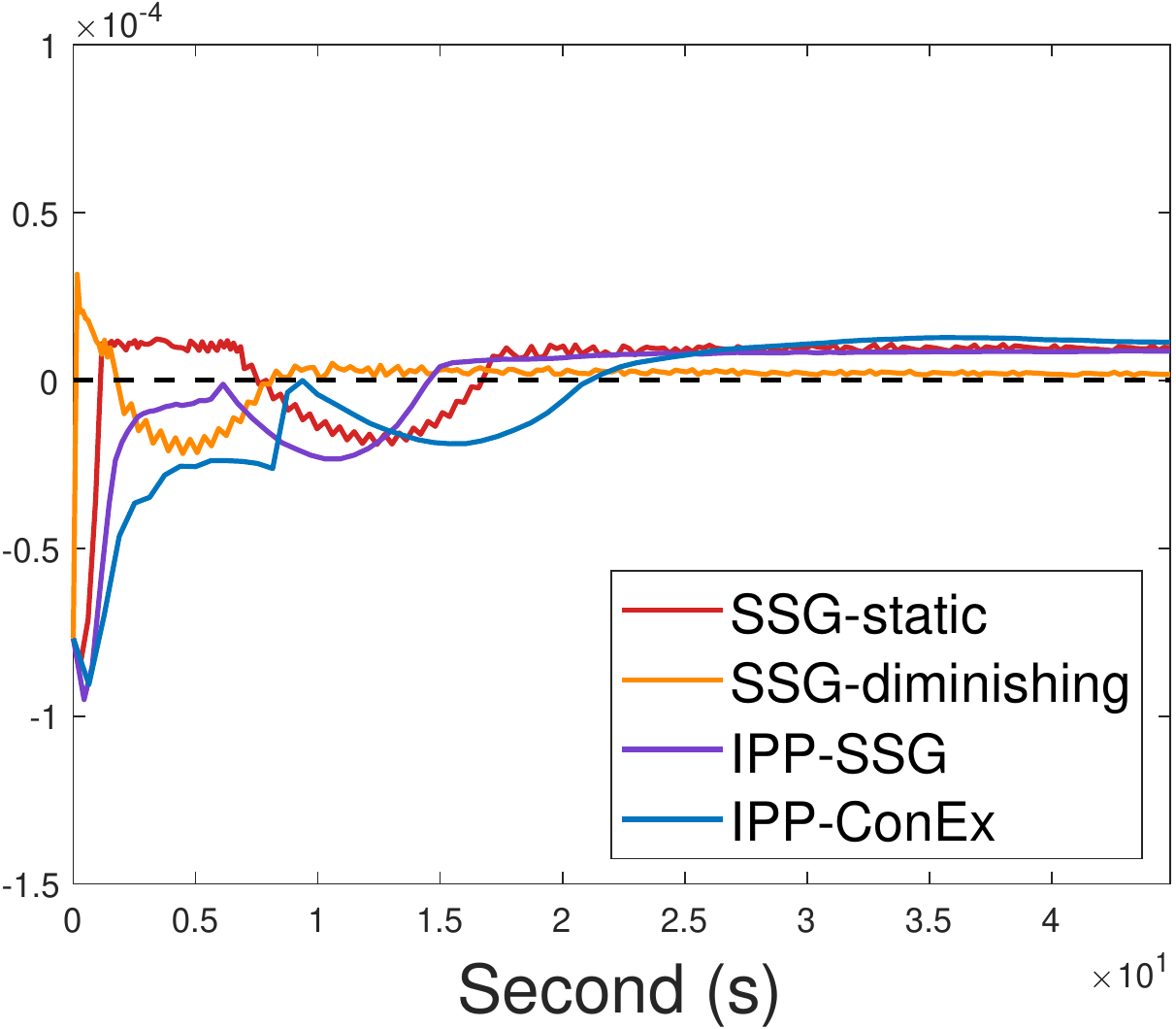}
           \\
        \raisebox{12ex}{\small{\rotatebox[origin=c]{90}{Near Stationarity}}}
		& \hspace*{-0.06in}\includegraphics[width=0.30\textwidth]{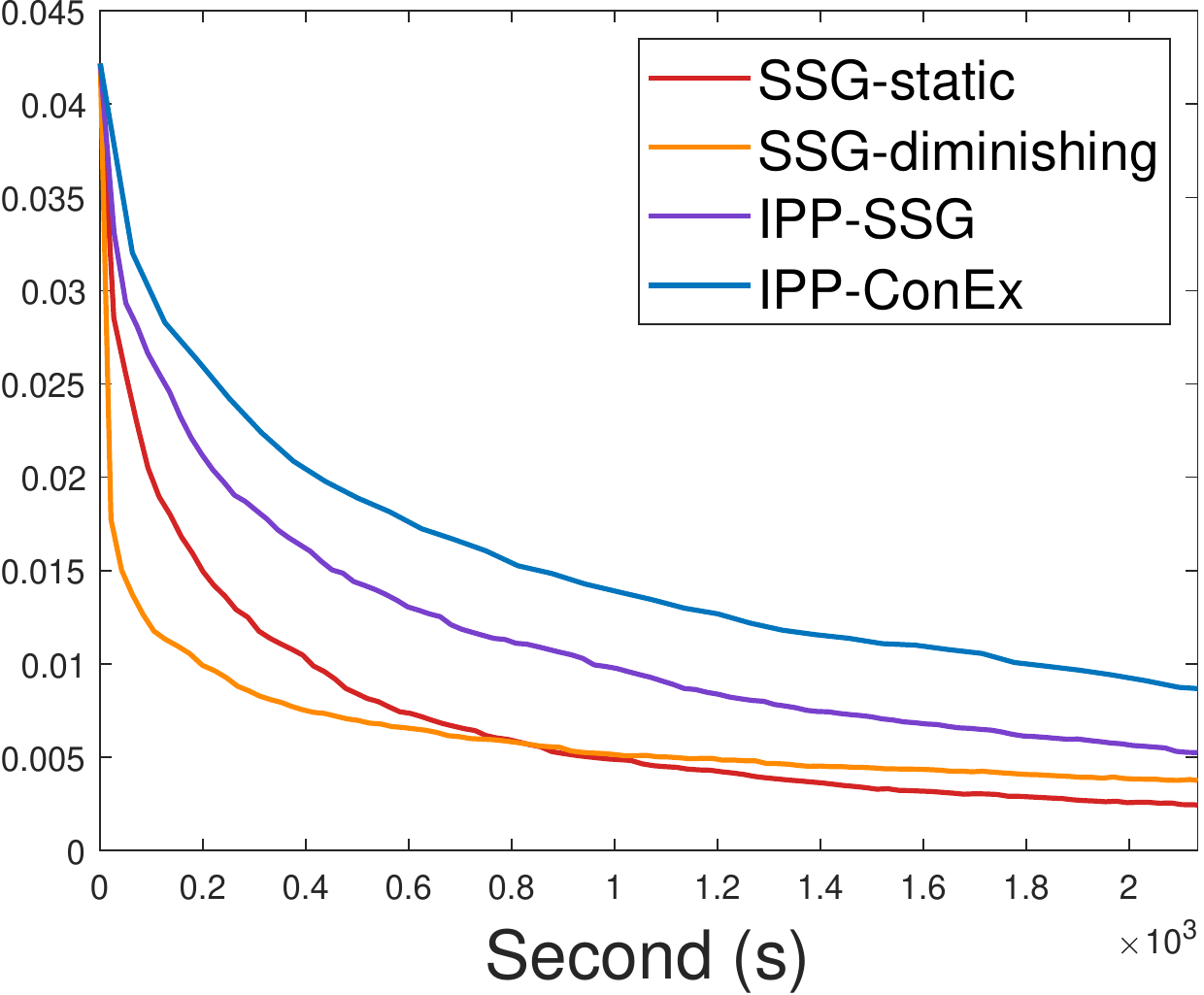}
		& \hspace*{-0.06in}\includegraphics[width=0.30\textwidth]{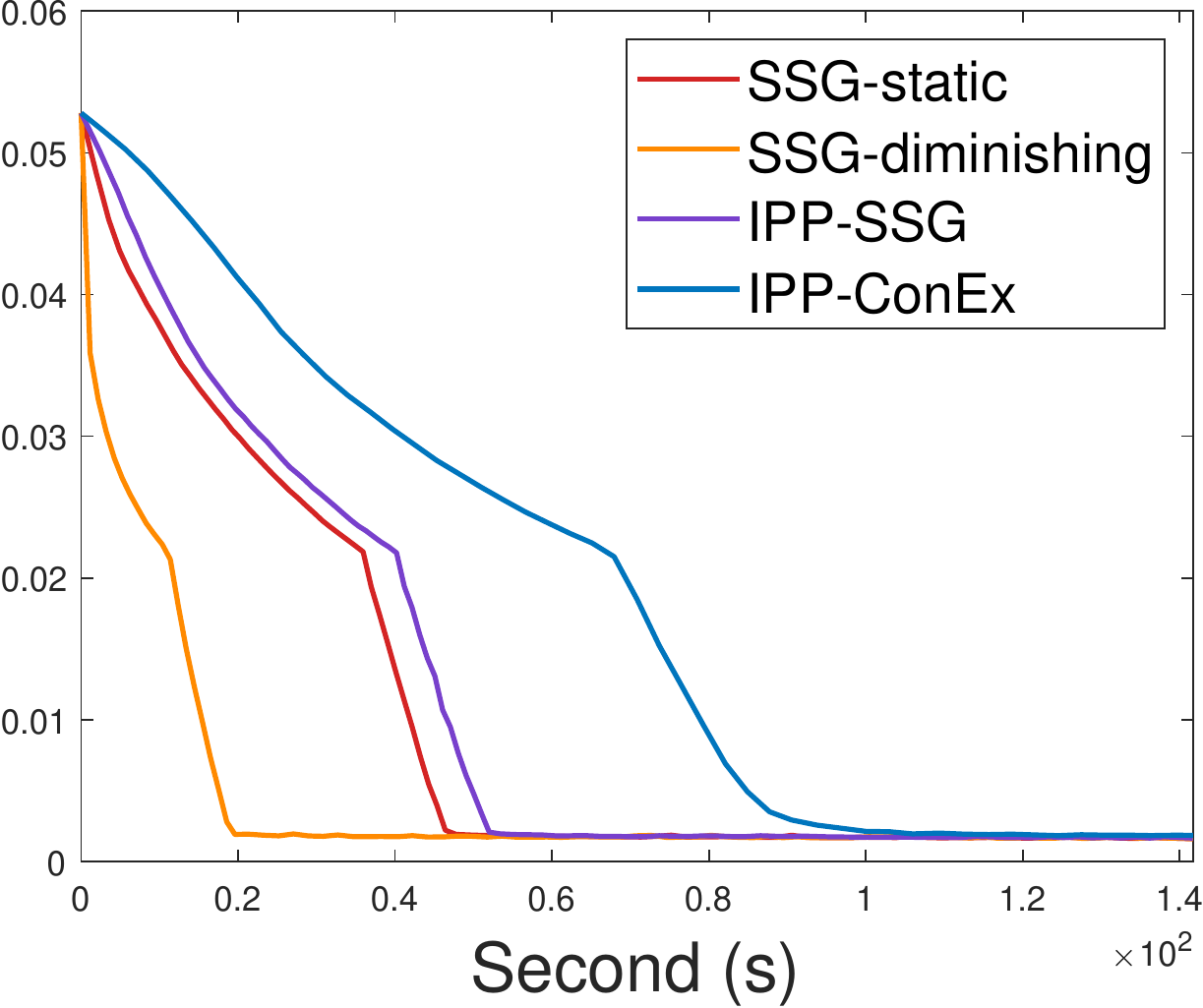}
		& \hspace*{-0.06in}\includegraphics[width=0.30\textwidth]{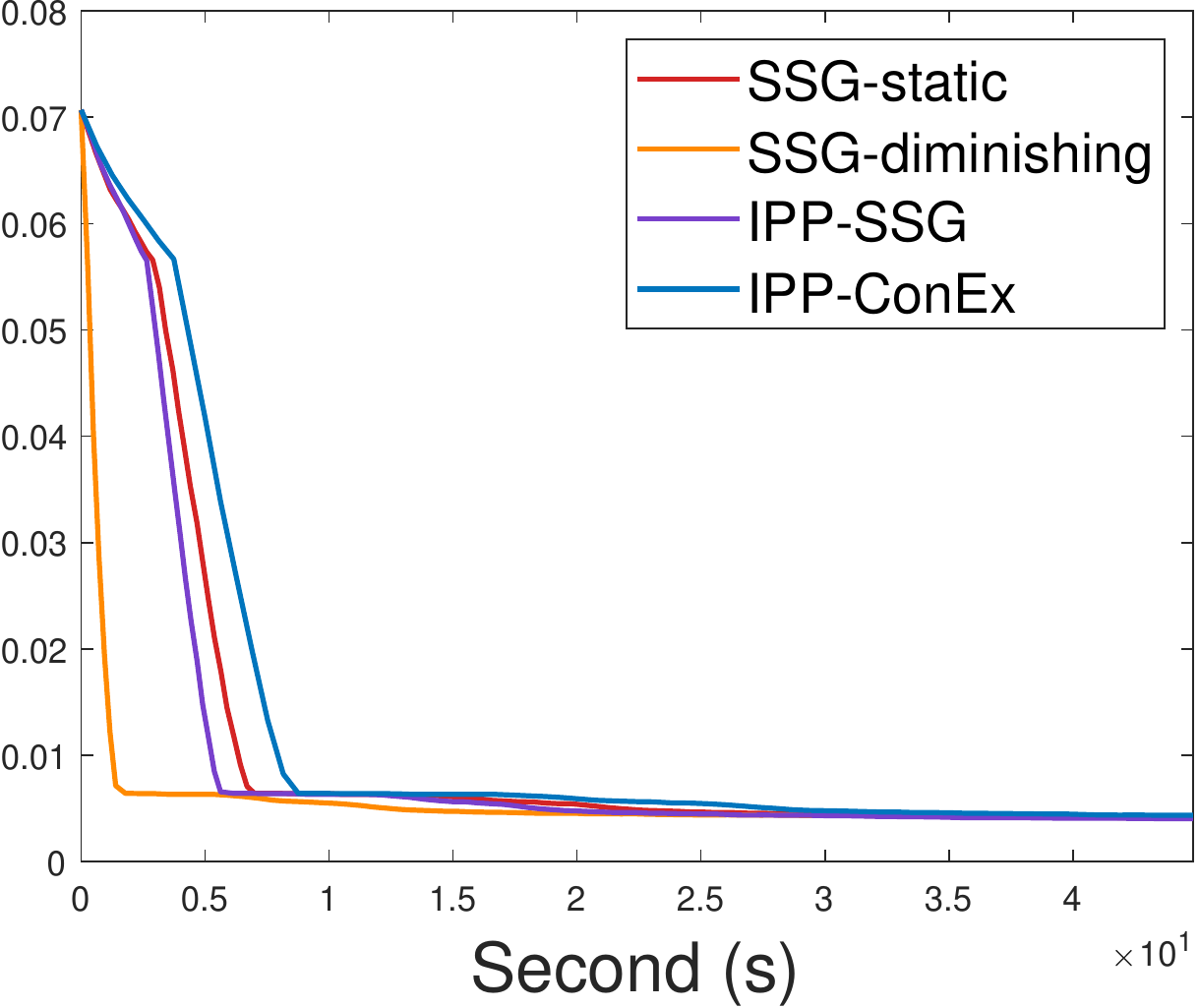}
        \end{tabular}
	\caption{Performances vs CPU time on classification problems with ROC-based fairness.}  \label{fig:figure_convex_experiment_cputime}
	\vspace{-0.1in}
\end{figure*}

In this section, we first present the performances versus the used CPU time of the methods in comparisons on problem \eqref{eq:ROCfairnessClassification} in Figure~\ref{fig:figure_convex_experiment_cputime}. The meanings of $y$-axes in Figure~\ref{fig:figure_convex_experiment_cputime} remain the same as in Figure~\ref{fig:figure_convex_experiment_iteration}, but the $x$-axis in each plot represents the CPU time in seconds the methods use. According to Figure~\ref{fig:figure_convex_experiment_cputime}, the relative speed among both SSG methods and IPP-SSG remain nearly the same as in Figure~\ref{fig:figure_convex_experiment_iteration} while IPP-ConEx is about twice time-consuming compared to others. This can be explained by the algorithmic designs of these methods. In fact, the two SSG methods and IPP-SSG only need to compute either $\bzt_f^{(t)}$ or $\bzt_g^{(t)}$ in each (inner) iteration while IPP-SSG computes both $\bzt_f^{(t)}$ and $\bzt_g^{(t)}$ per inner iteration and thus induces  additional computational cost.

\begin{figure*}
     \begin{tabular}[h]{@{}c|ccc@{}}
      & a9a & bank & COMPAS \\
		\hline \vspace*{-0.1in}\\
		\raisebox{12ex}{\small{\rotatebox[origin=c]{90}{Objective}}}
		& \hspace*{-0.06in}\includegraphics[width=0.30\textwidth]{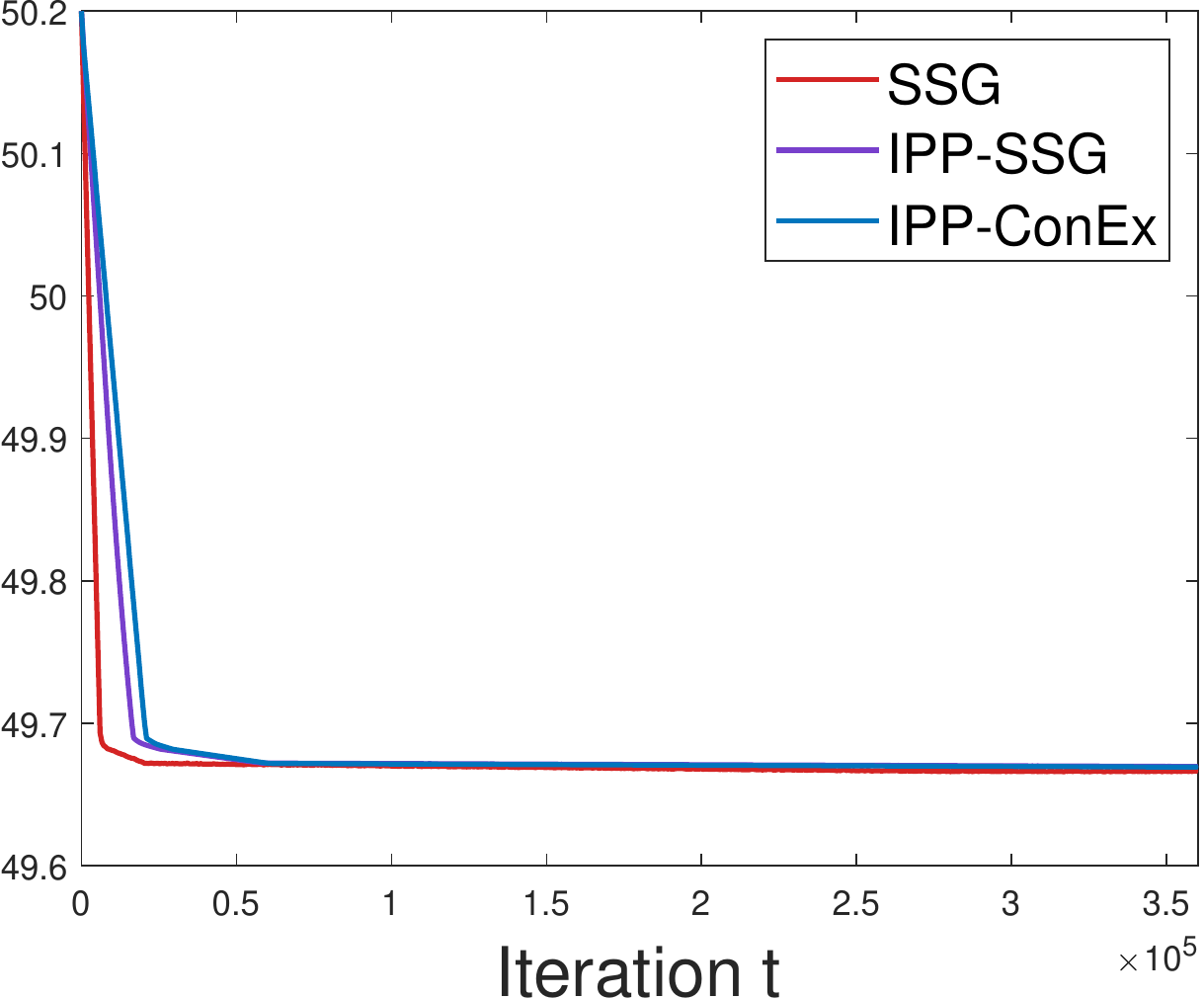}
		& \hspace*{-0.06in}\includegraphics[width=0.30\textwidth]{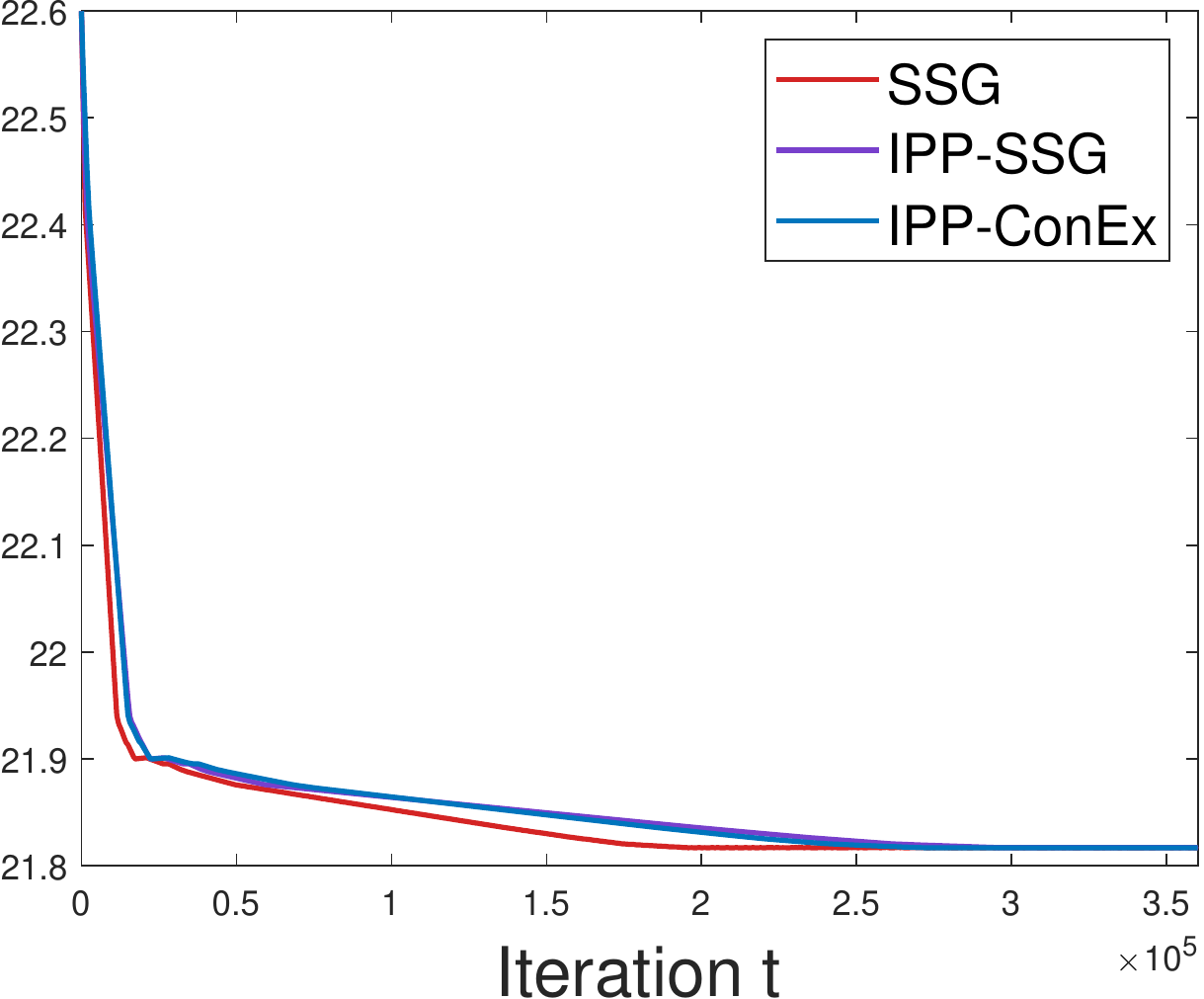}
		& \hspace*{-0.06in}\includegraphics[width=0.30\textwidth]{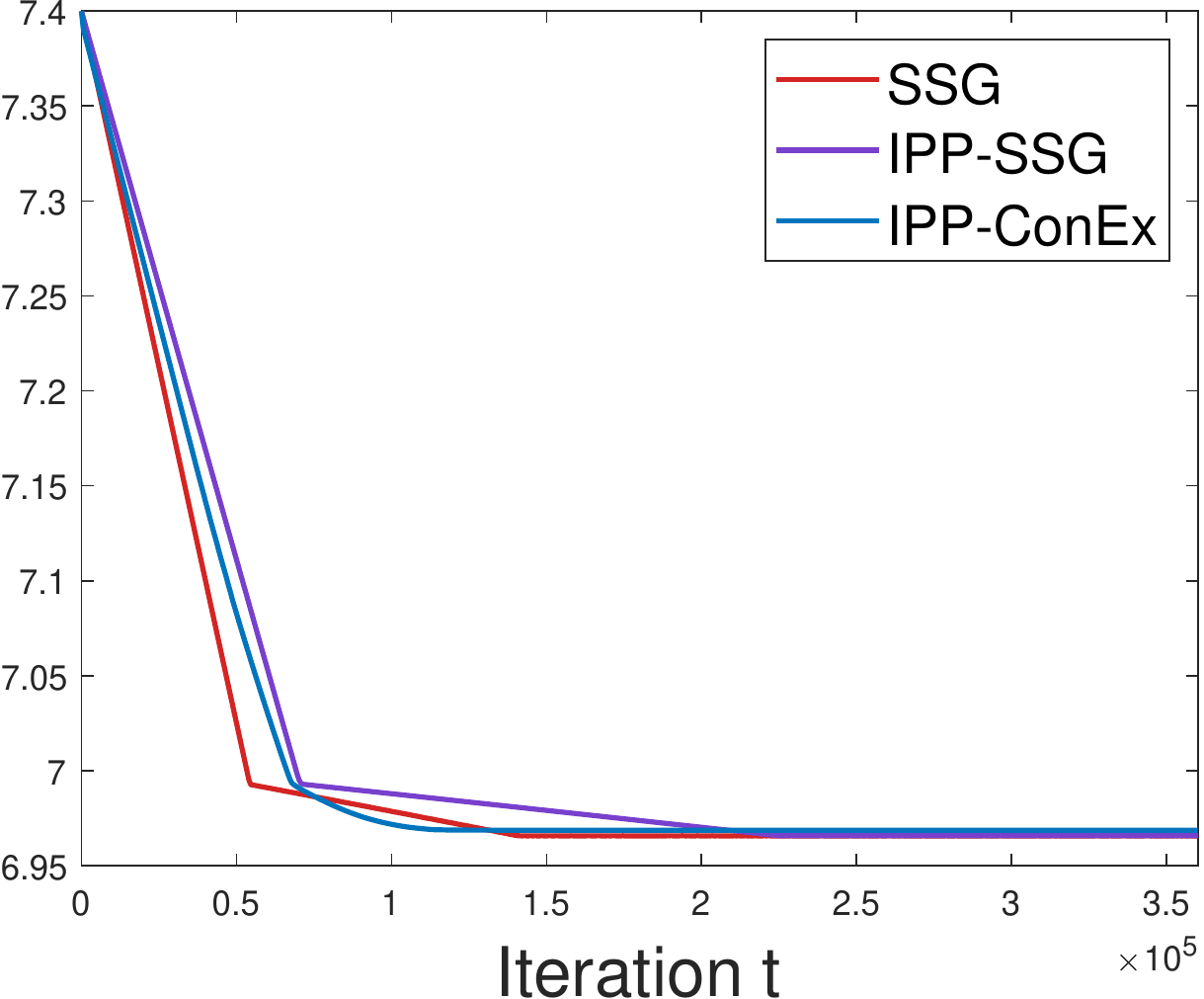}
          \\
		\raisebox{12ex}{\small{\rotatebox[origin=c]{90}{Infeasibility}}}
		& \hspace*{-0.06in}\includegraphics[width=0.30\textwidth]{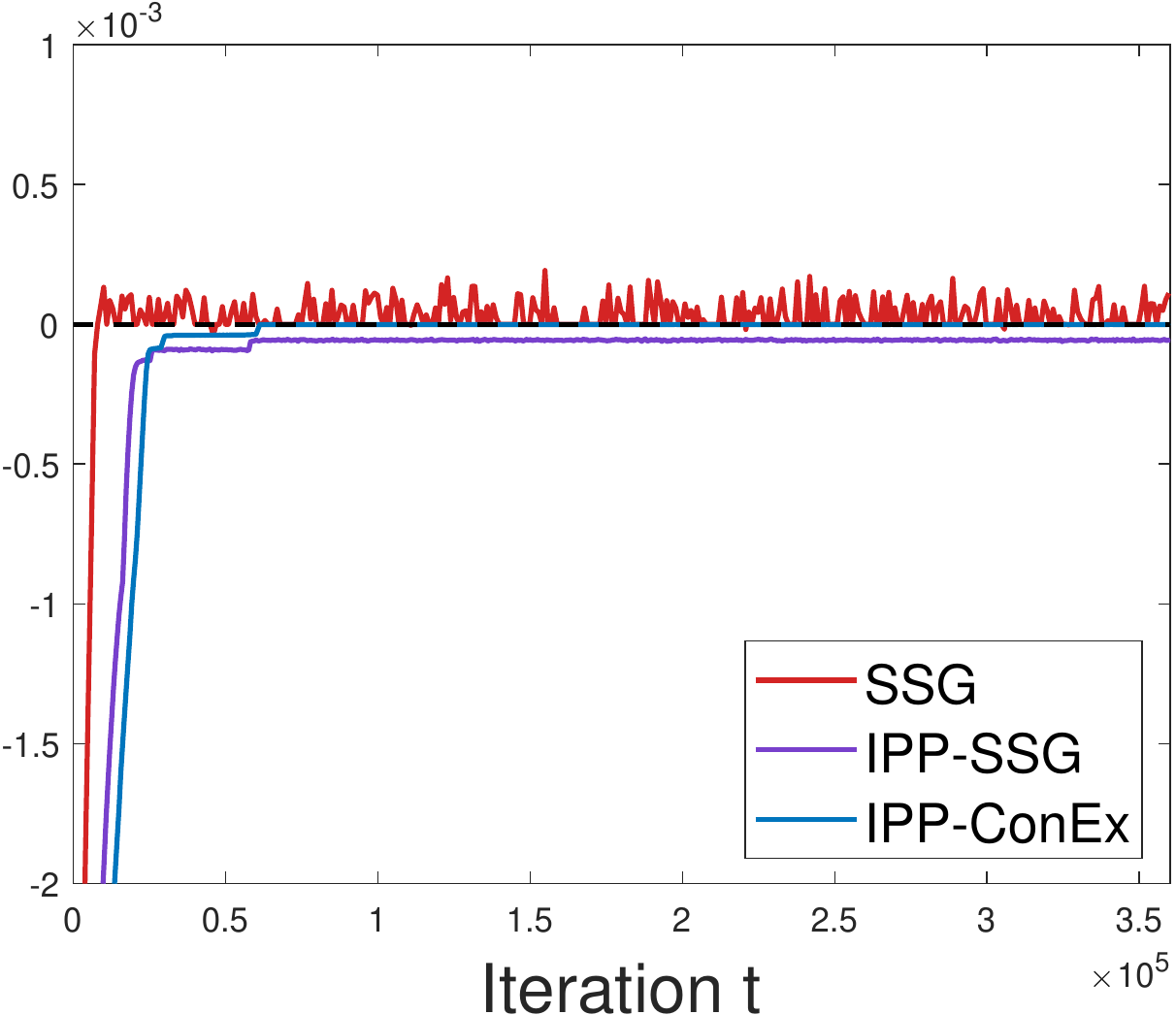}
		& \hspace*{-0.06in}\includegraphics[width=0.30\textwidth]{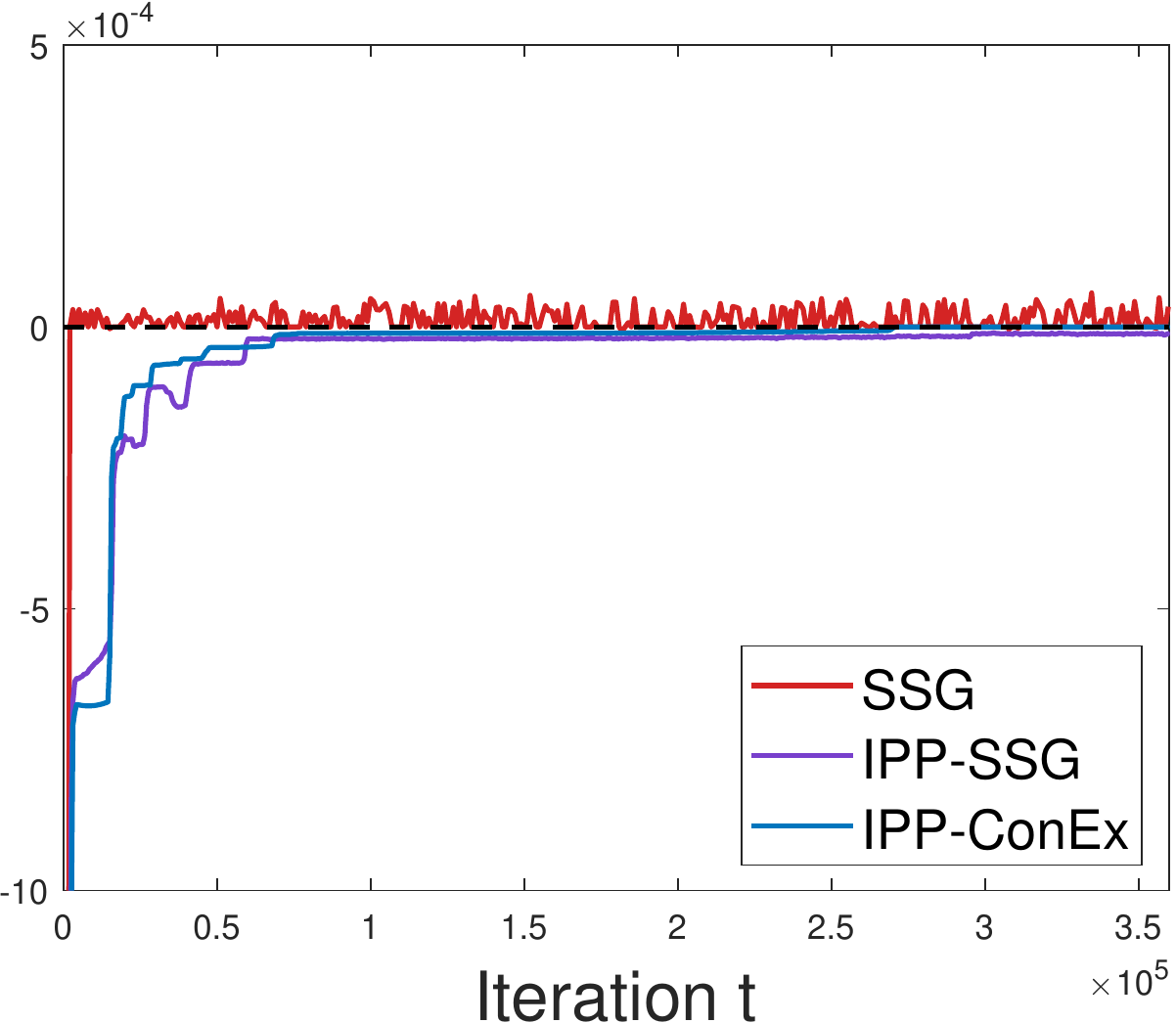}
		& \hspace*{-0.06in}\includegraphics[width=0.30\textwidth]{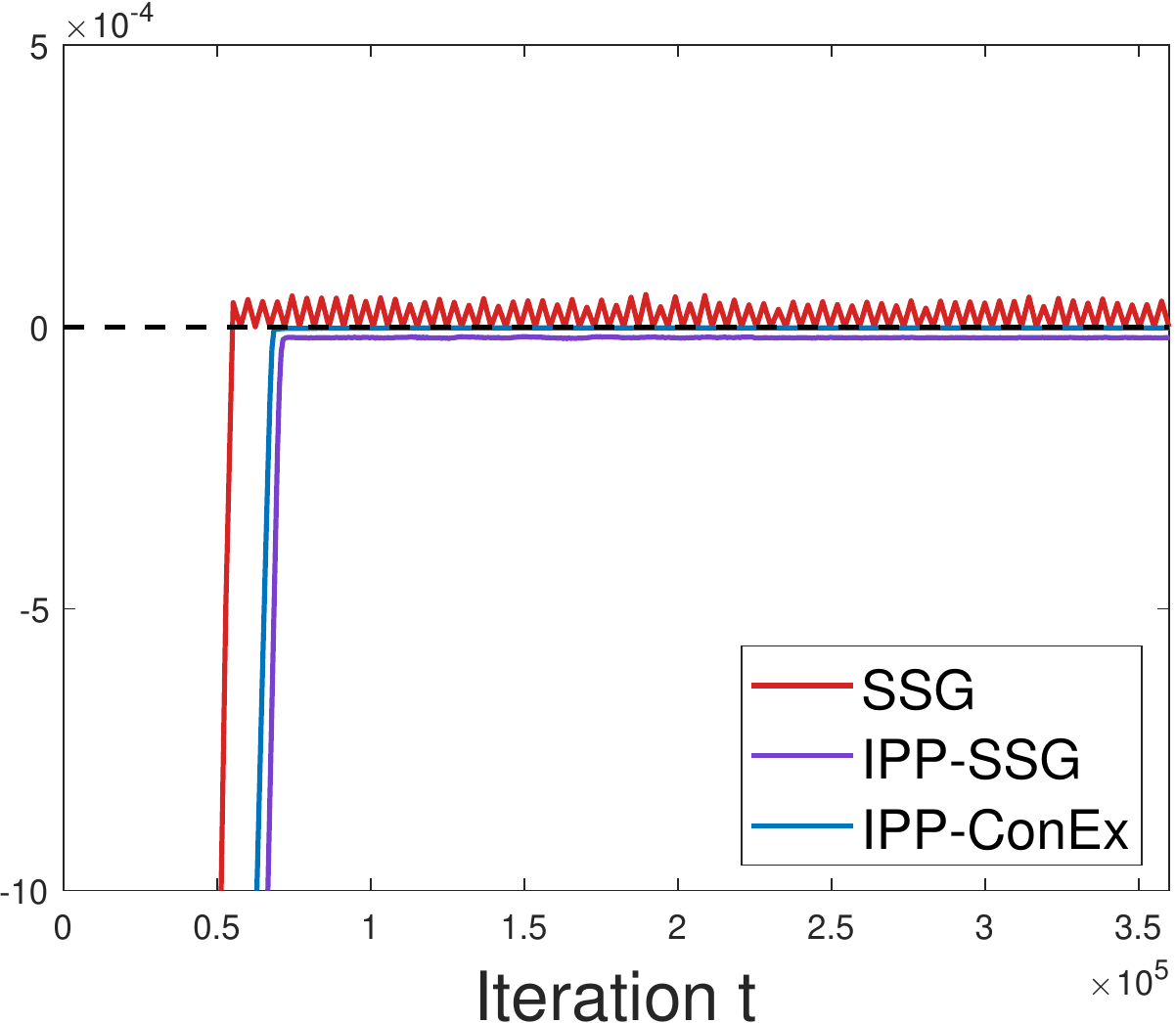}
           \\
        \raisebox{12ex}{\small{\rotatebox[origin=c]{90}{Near Stationarity}}}
		& \hspace*{-0.06in}\includegraphics[width=0.30\textwidth]{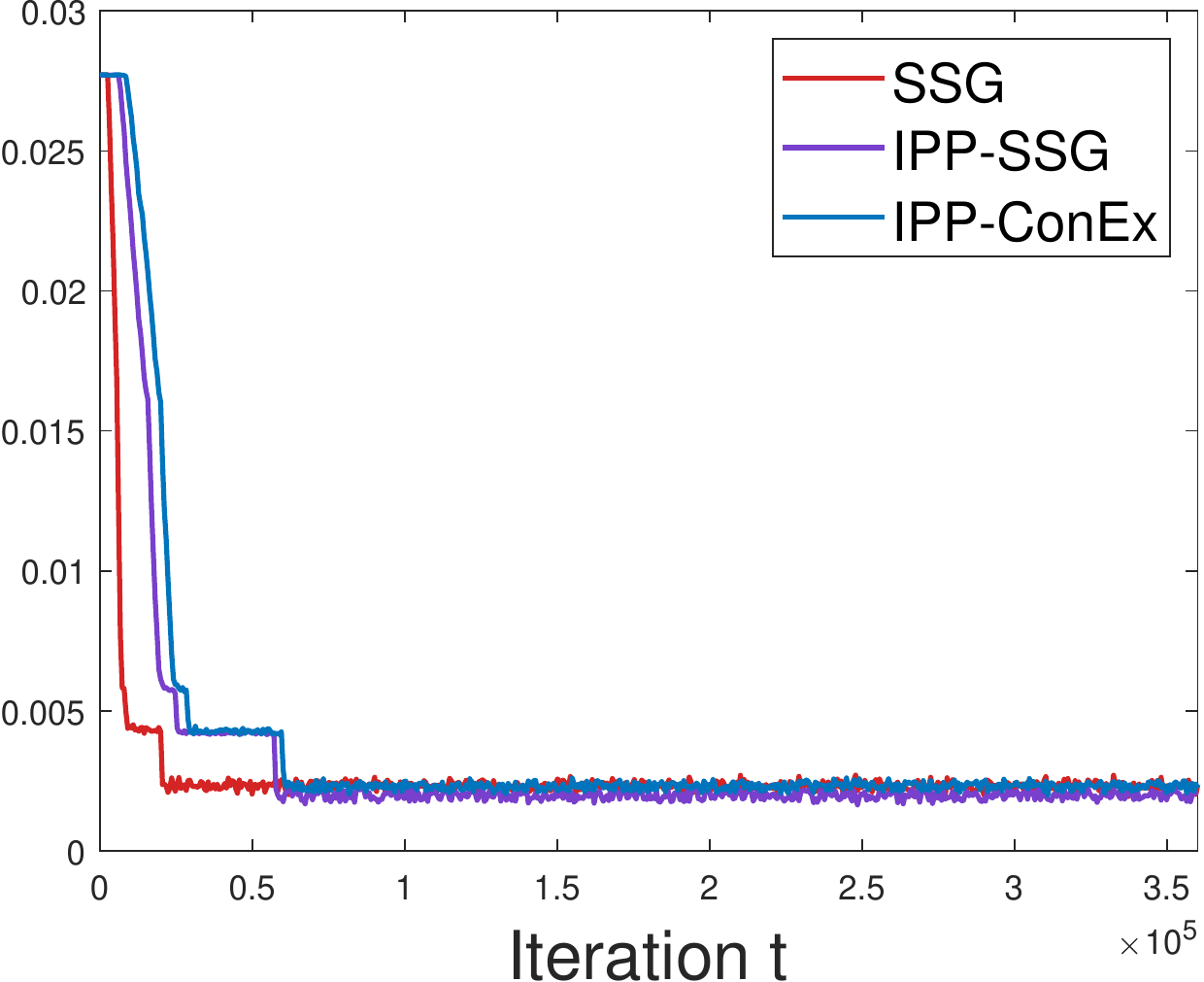}
		& \hspace*{-0.06in}\includegraphics[width=0.30\textwidth]{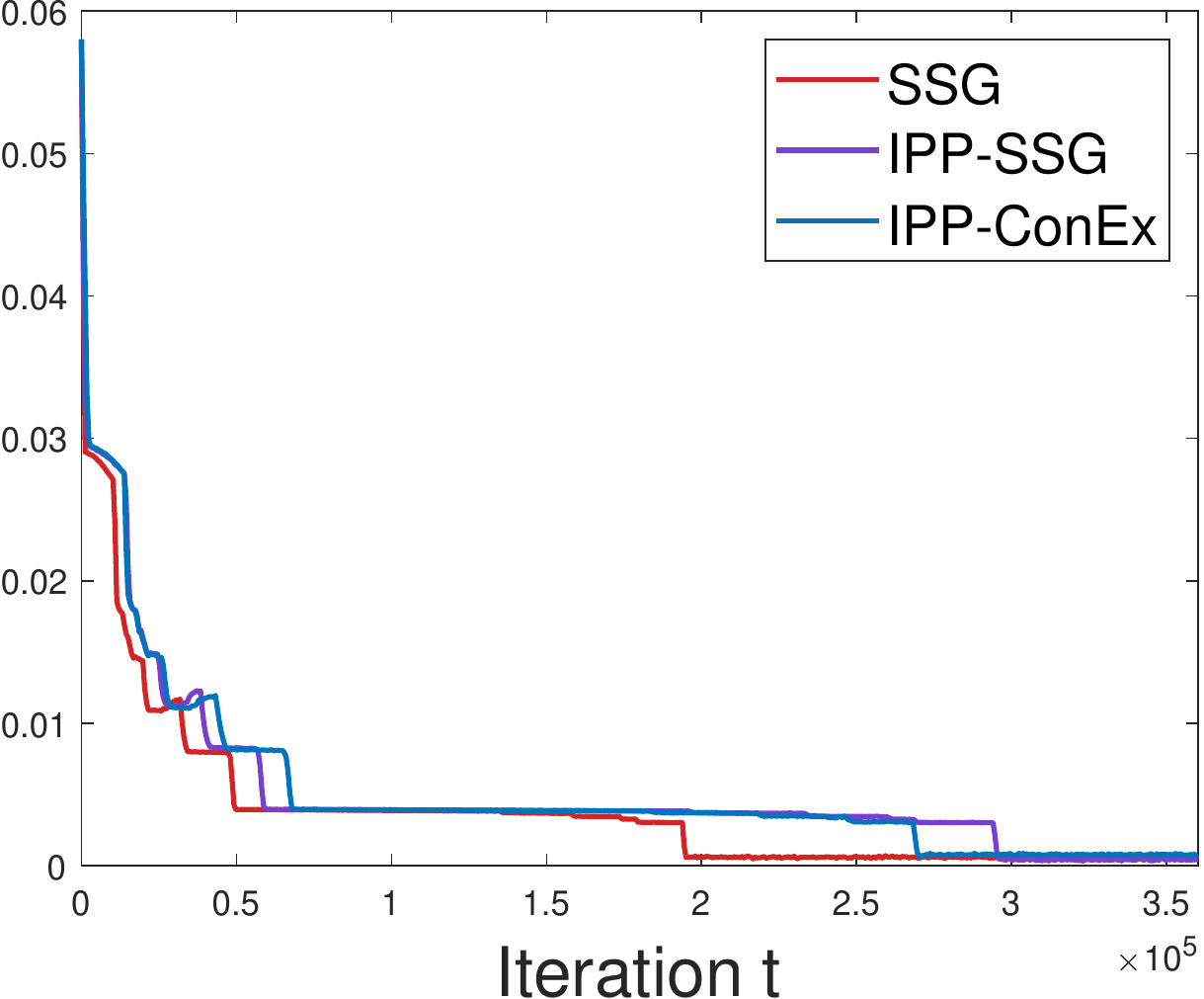}
		& \hspace*{-0.06in}\includegraphics[width=0.30\textwidth]{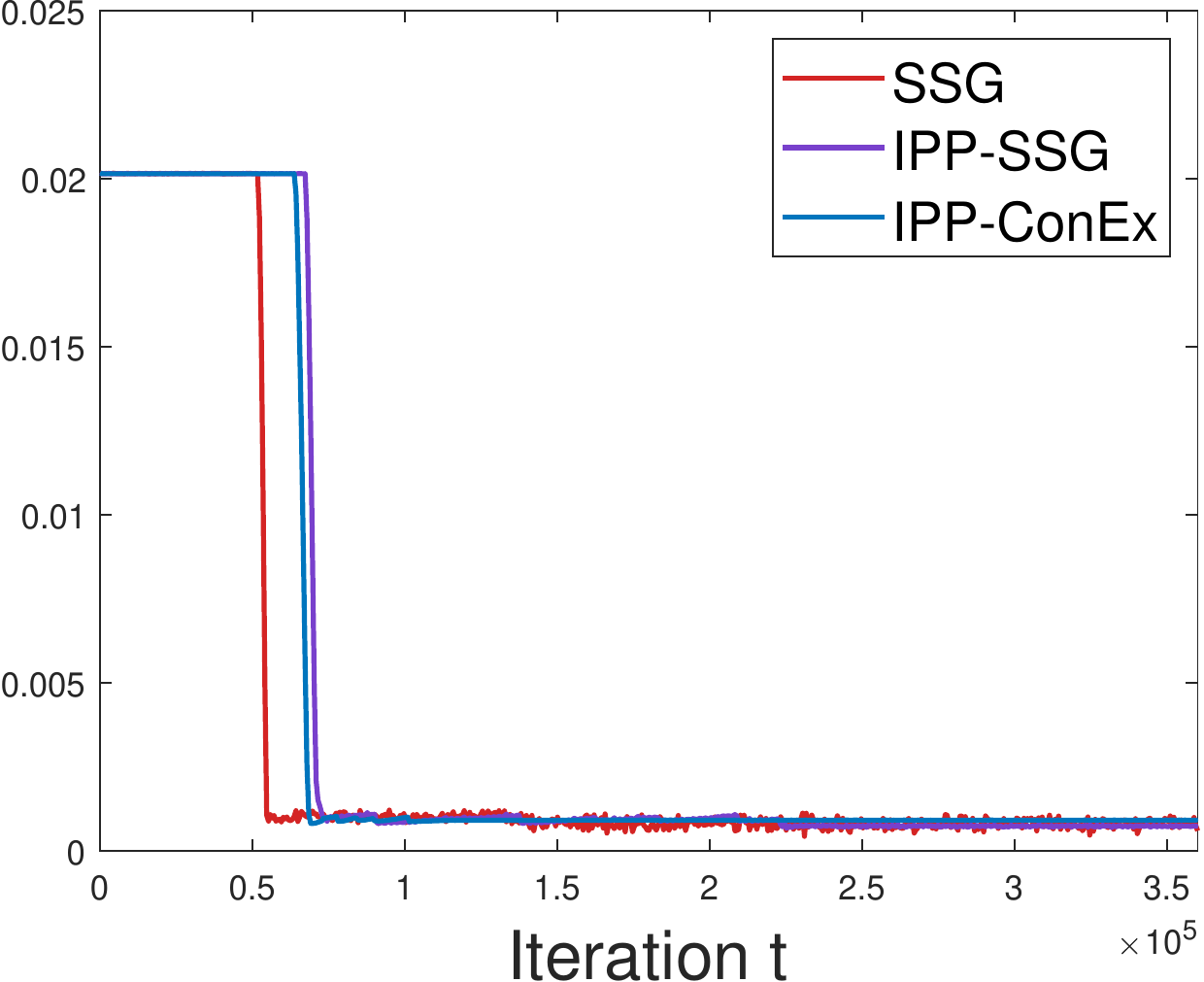}
        \end{tabular}
	\caption{Performances vs number of iterations on classification problems with demographic equity constraint.}  \label{fig:figure_weakly_convex_experiment_iteration}
	\vspace{-0.1in}
\end{figure*}

\begin{figure*}
     \begin{tabular}[h]{@{}c|ccc@{}}
      & a9a & bank & COMPAS \\
		\hline \vspace*{-0.1in}\\
		\raisebox{12ex}{\small{\rotatebox[origin=c]{90}{Objective}}}
		& \hspace*{-0.06in}\includegraphics[width=0.30\textwidth]{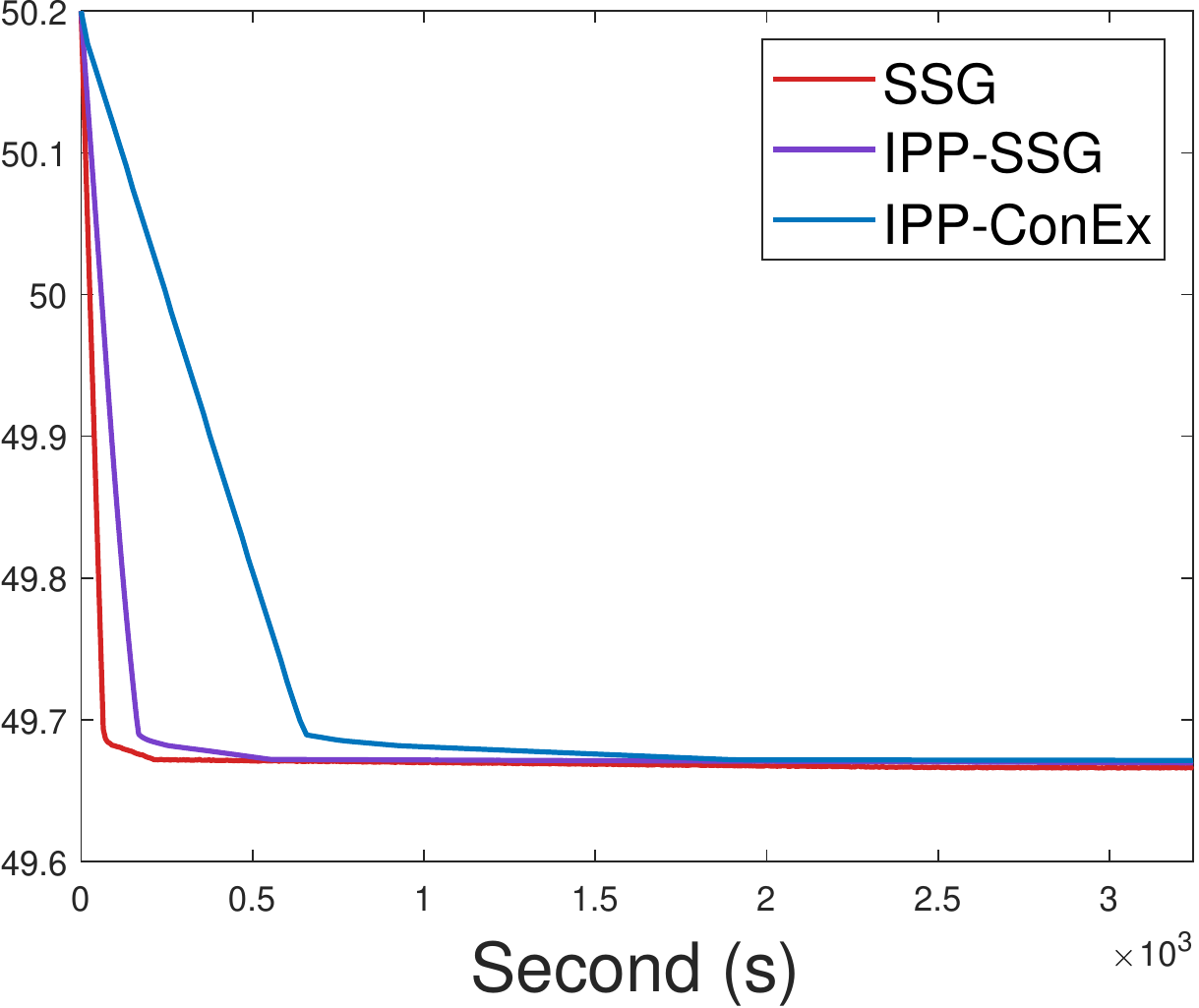}
		& \hspace*{-0.06in}\includegraphics[width=0.30\textwidth]{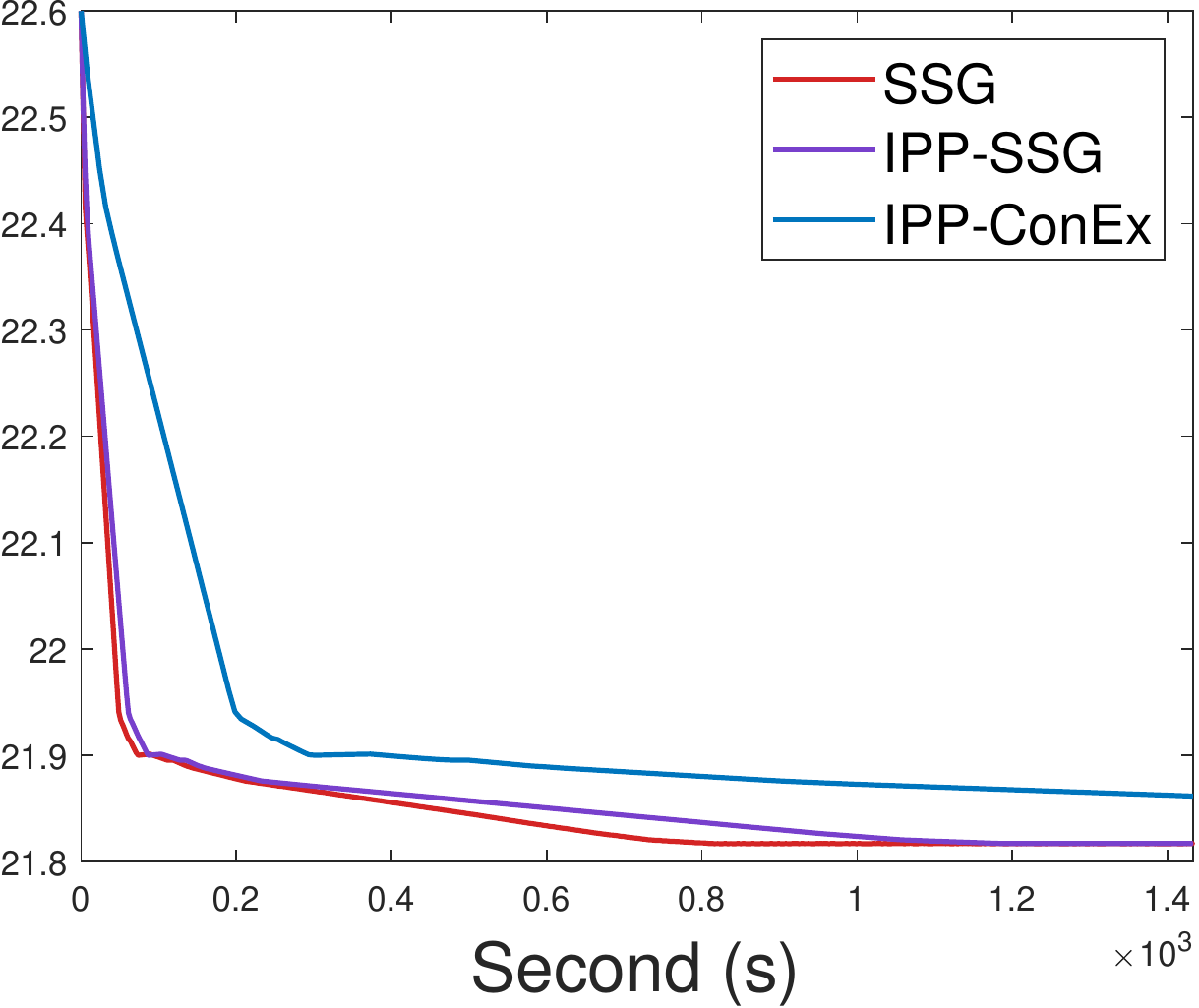}
		& \hspace*{-0.06in}\includegraphics[width=0.30\textwidth]{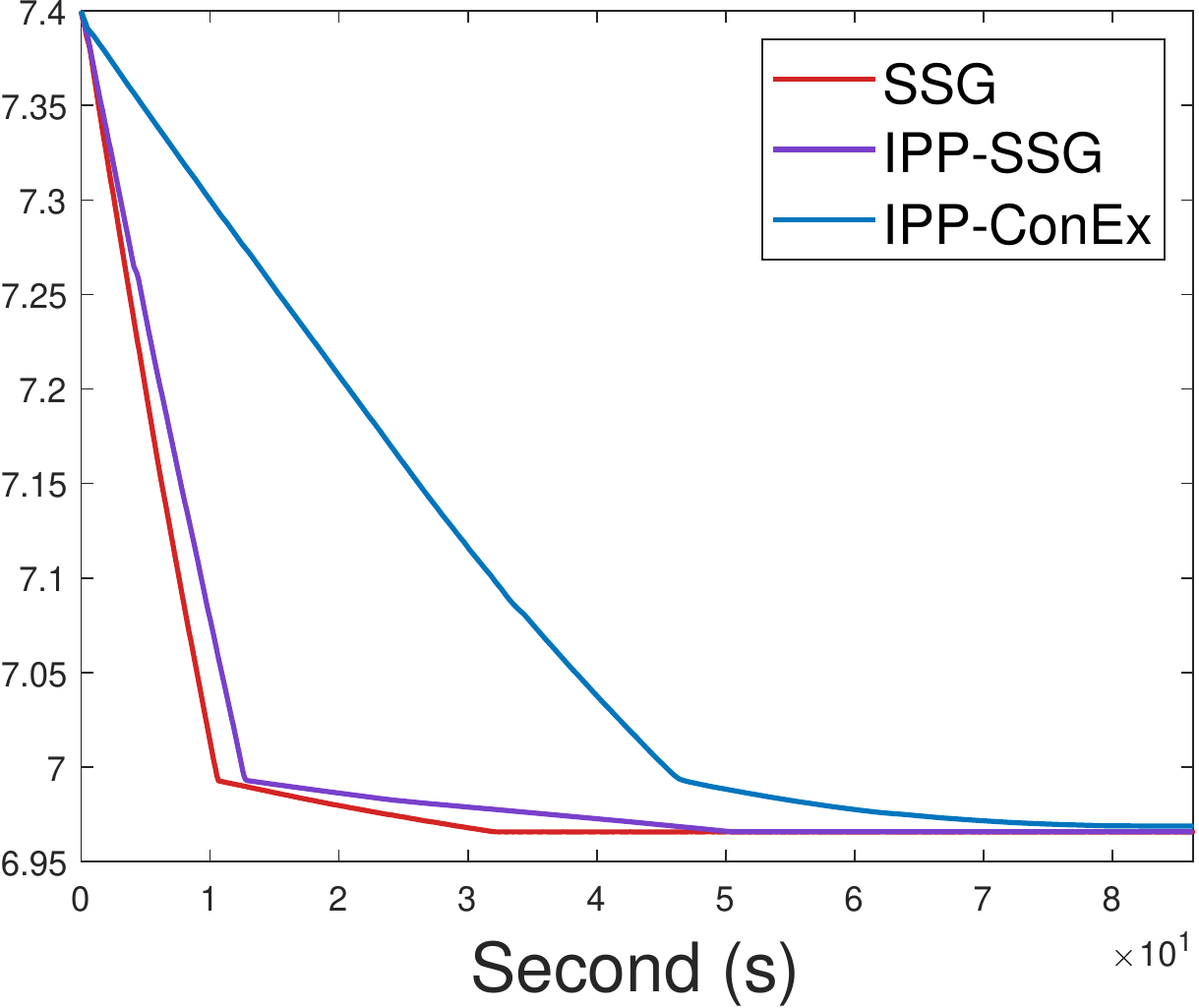}
          \\
		\raisebox{12ex}{\small{\rotatebox[origin=c]{90}{Infeasibility}}}
		& \hspace*{-0.06in}\includegraphics[width=0.30\textwidth]{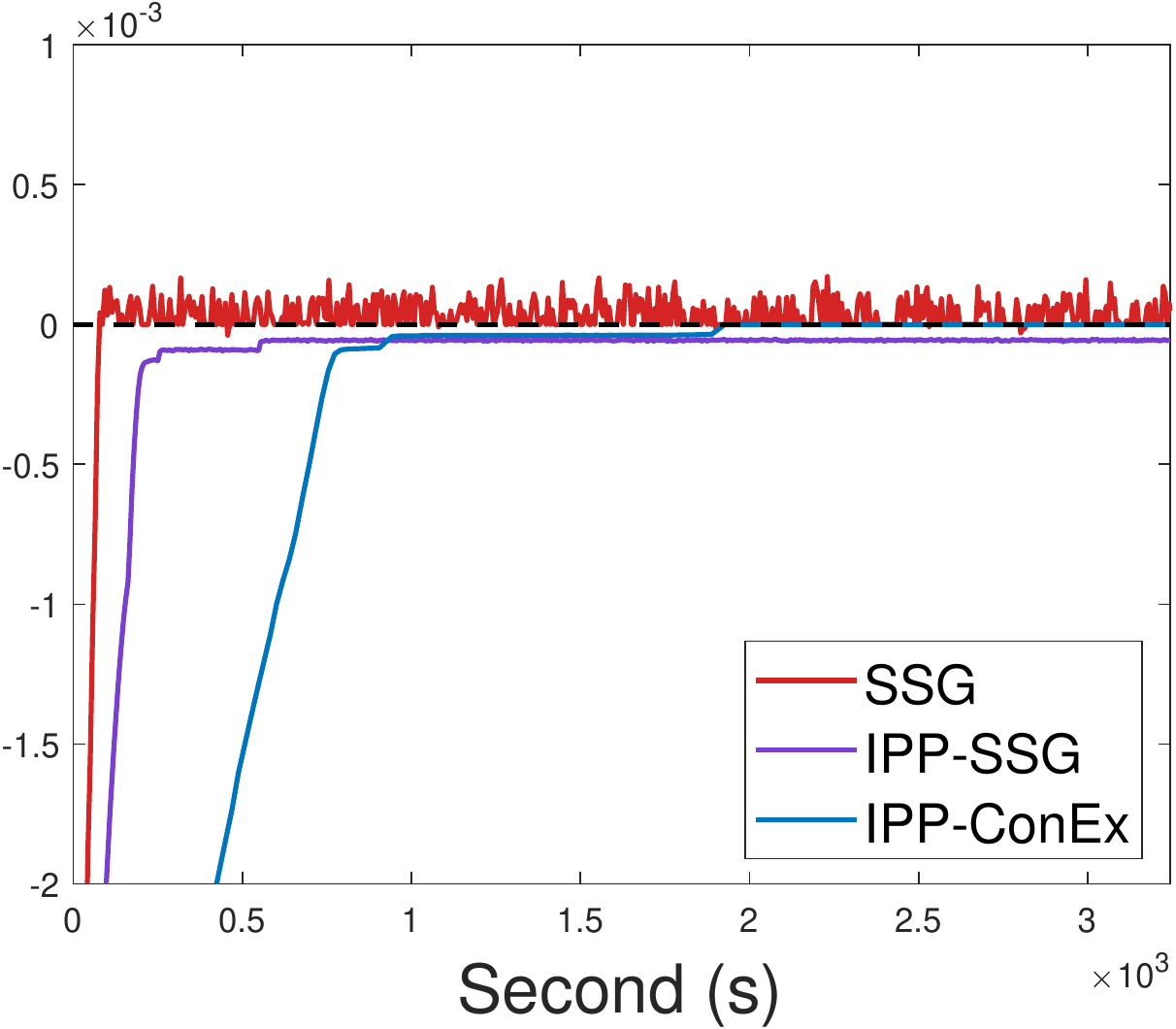}
		& \hspace*{-0.06in}\includegraphics[width=0.30\textwidth]{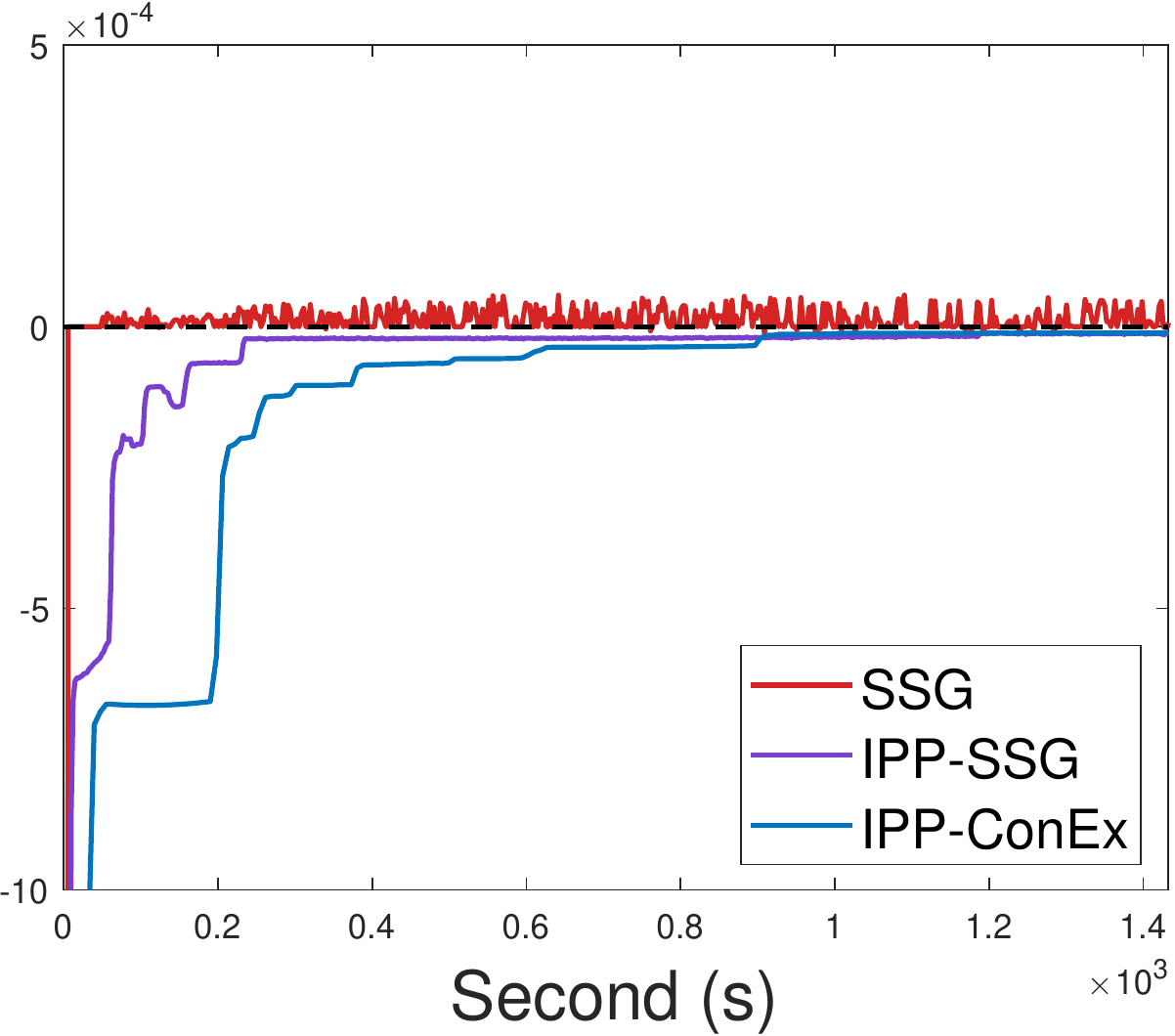}
		& \hspace*{-0.06in}\includegraphics[width=0.30\textwidth]{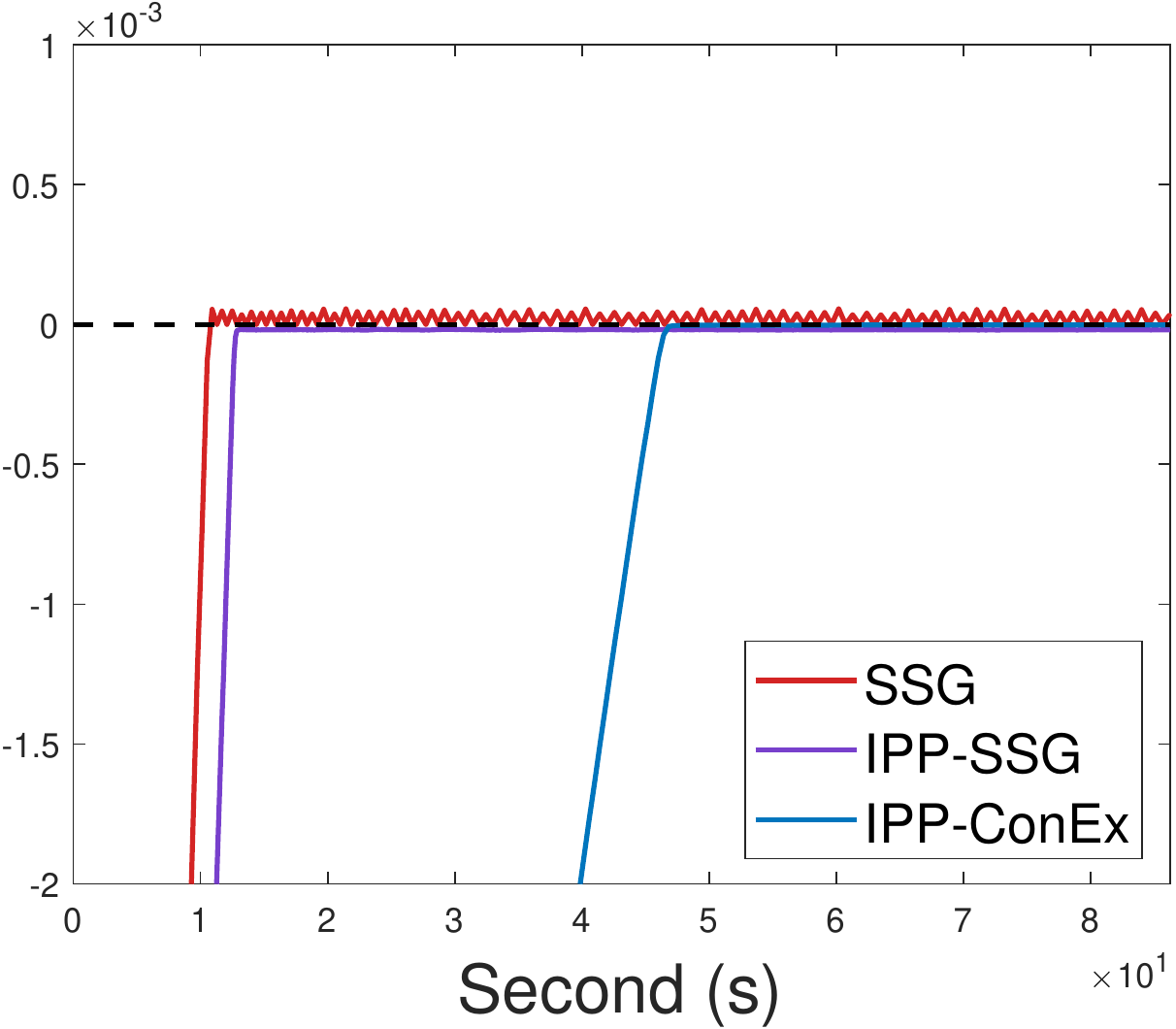}
           \\
        \raisebox{12ex}{\small{\rotatebox[origin=c]{90}{Near Stationarity}}}
		& \hspace*{-0.06in}\includegraphics[width=0.30\textwidth]{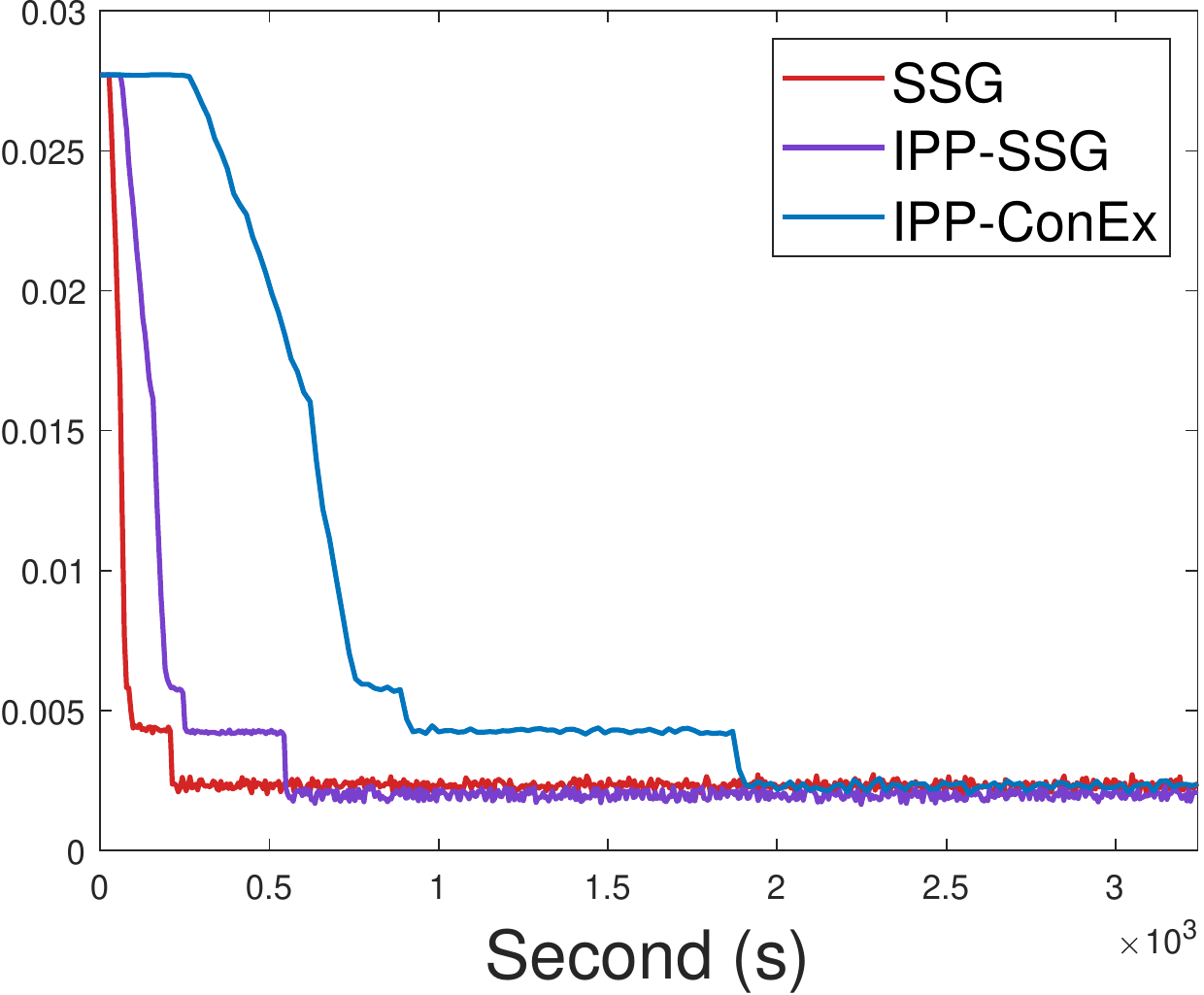}
		& \hspace*{-0.06in}\includegraphics[width=0.30\textwidth]{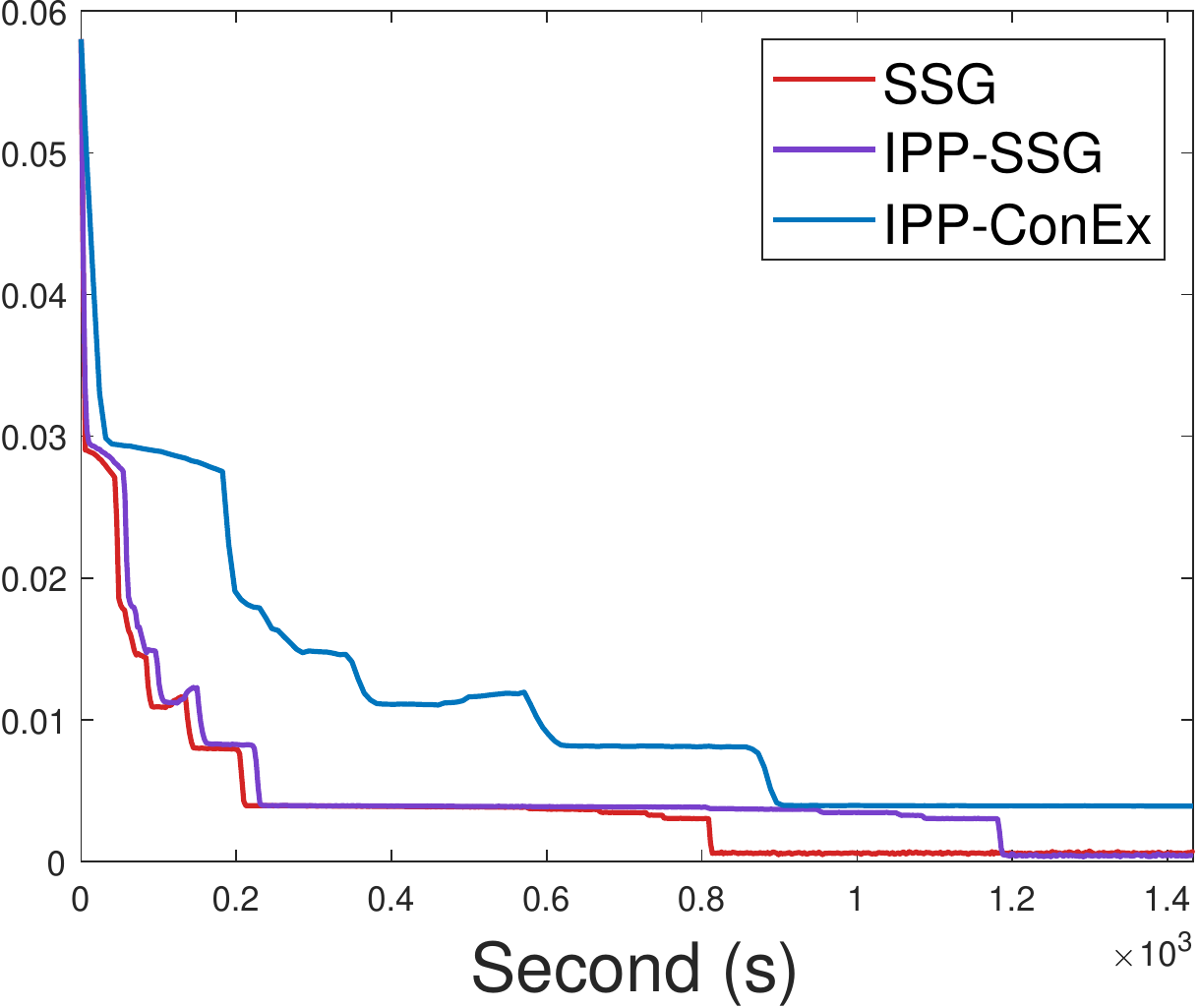}
		& \hspace*{-0.06in}\includegraphics[width=0.30\textwidth]{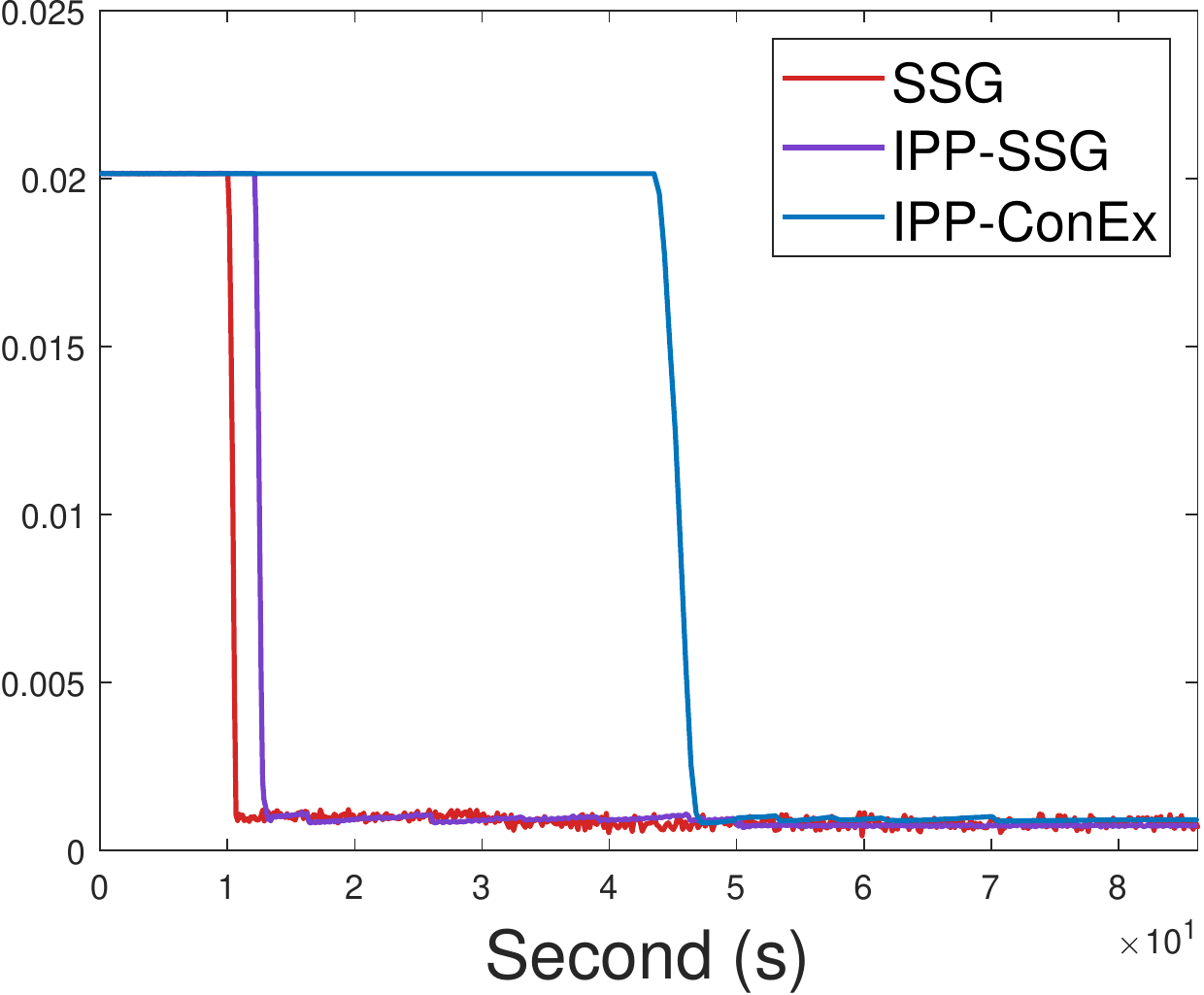}
        \end{tabular}
	\caption{Performances vs CPU time on classification problems with demographic equity constraint.}  \label{fig:figure_weakly_convex_experiment_cputime}
	\vspace{-0.1in}
\end{figure*}

Then, we present the numerical performances versus both the number of iterations and the CPU time of the methods in comparisons on problem \eqref{eq:demographicparity_fairnessClassification}. The results on each dataset are given in Figure~\ref{fig:figure_weakly_convex_experiment_iteration} and Figure~\ref{fig:figure_weakly_convex_experiment_cputime} where the meanings of axes and the procedure of calculating the (approximate) near stationarity are the same as in Figure~\ref{fig:figure_convex_experiment_iteration} and Figure~\ref{fig:figure_convex_experiment_cputime}, respectively. Since computing $\widehat\bx(\bx^{(t)})$ with a high precision for each $t$ is time-consuming, we only report near stationarity at $600$ equally spaced iterations. Figure~\ref{fig:figure_weakly_convex_experiment_iteration} shows that the SSG method is only a little faster than the other two methods in reducing the objective value and the (approximate) near stationarity while keeping the solutions nearly feasible. This is consistent with our theoretical finding that the three methods have the similar oracle complexity. Figure~\ref{fig:figure_weakly_convex_experiment_cputime} shows that the relative speeds among the two SSG methods and IPP-SSG are nearly the same as in Figure~\ref{fig:figure_weakly_convex_experiment_iteration} while IPP-ConEx is about twice time-consuming. This is again consistent with the algorithmic designs of these methods and can be explained by the same reasons we gave above for  Figure~\ref{fig:figure_convex_experiment_cputime}.

\end{document}